\documentclass[11pt, reqno]{amsart}
\usepackage[utf8]{inputenc}	
\usepackage{amssymb,amscd}
\usepackage{multirow,booktabs}
\usepackage[table]{xcolor}
\usepackage{fullpage}
\usepackage{lastpage}
\usepackage{fancyhdr}
\usepackage{mathrsfs}
\usepackage{wrapfig}
\usepackage{graphicx}
\usepackage{setspace}
\usepackage{calc}
\usepackage{bm}
\usepackage{tikz}
\usepackage{enumerate}
\usepackage{subfigure}
\usepackage{multicol}
\usepackage[colorlinks,
            linkcolor=blue,
            citecolor=red,
            ]{hyperref}
\usepackage{cancel}
\usepackage[all]{xy}
\usepackage{geometry}
\usepackage{empheq}
\usepackage{framed}
\usepackage{dsfont}
\usepackage{xcolor}
\theoremstyle{plain}
\newtheorem{theorem}{Theorem}[section]
\newtheorem{lemma}[theorem]{Lemma}
\newtheorem{proposition}[theorem]{Proposition}
\newtheorem{corollary}[theorem]{Corollary}

\theoremstyle{remark}
\newtheorem{definition}{Definition}

\newtheorem{remark}{Remark}[section]

\usepackage{indentfirst}
\begin{document}

\newcommand{\QQ}{\mathbb{Q}}
\newcommand{\RR}{\mathbb{R}}
\newcommand{\ZZ}{\mathbb{Z}}
\newcommand{\NN}{\mathbb{N}}
\newcommand{\Nor}{\mathscr{N}}
\newcommand{\CC}{\mathbb{C}}
\newcommand{\HH}{\mathbb{H}}
\newcommand{\EE}{\mathbb{E}}
\newcommand{\Var}{\operatorname{Var}}
\newcommand{\PP}{\mathbb{P}}
\newcommand{\Rd}{\mathbb{R}^d}
\newcommand{\Rn}{\mathbb{R}^n}
\newcommand{\BHH}{\overline{\mathbb{H}}}
\newcommand{\system}{(\Omega,\mathcal{F},\mu,T)}
\newcommand{\FF}{\mathcal{F}}
\newcommand{\GG}{\mathcal{G}}
\newcommand{\MBS}{(\Omega,\mathcal{F})}
\newcommand{\MS}{(\Omega,\mathcal{F},\mu)}
\newcommand{\PS}{(\Omega,\mathcal{F},\mathbb{P})}
\newcommand{\Def}{\overset{\text{def}}{=}}
\newcommand{\Ser}[2]{#1_1,#1_2,\dots,#1_#2}
\newcommand{\independent}{\perp\mkern-9.5mu\perp}
\def\avint{\mathop{\,\rlap{-}\!\!\int\!\!\llap{-}}\nolimits}

\author[Shuo Qin]{Shuo Qin}
\address[Shuo Qin]{Beijing Institute of Mathematical Sciences and Applications}
\address[Shuo Qin]{NYU-ECNU Institute of Mathematical Sciences at NYU Shanghai and Courant Institute of Mathematical Sciences}
\email{qinshuo@bimsa.cn}

\title{Step-reinforced random walks and one-half}
\date{}

  \begin{abstract}
    Under suitable moment assumptions, we show that a genuinely d-dimensional step-reinforced random walk undergoes a phase transition between recurrence and transience in dimensions $d=1,2$, and that it is transient for all reinforcement parameters in dimensions $d\geq 3$, which solves a conjecture of Bertoin.
  \end{abstract}
\maketitle

\section{General introduction}

\subsection{The model and background} 

Step-reinforced random walks form a class of random processes with reinforcement, whose step sequences are generated by an algorithm introduced by Simon \cite{MR0073085} in 1955. These walks have garnered significant attention in recent years, see e.g. \cite{MR4237267,MR4490996,MR4794980,bertenghi2021asymptotic, MR4116724,MR4221659, hu2025berry, aguech2025class, MR4168396, MR3827299}.

A canonical example is the elephant random walk (ERW), originally formulated by Schütz and Trimper \cite{schutz2004elephants} in dimension one and later extended to higher dimensions by Bercu and Laulin \cite{MR3962977}. The d-dimensional ERW starts at the origin at time $0$, with its first step selected uniformly from the set $\{\pm\bm{e}_1, \pm\bm{e}_2, \ldots, \pm\bm{e}_d\}$, where $\bm{e}_1, \ldots, \bm{e}_d$ denote the standard basis for $\mathbb{R}^d$. At each subsequent time step $n\geq 2$, the elephant uniformly samples a step from the past, and then, with probability $p\in [0,1]$, moves exactly in the same direction as the chosen step. Otherwise, it moves uniformly in one of the $2d - 1$ remaining directions. The parameter $p$ is commonly referred to as the memory parameter. For $p\geq 1/(2d)$, the following equivalent parametrization \cite{bertenghi2021asymptotic, MR3652690} 
\begin{equation}
    \label{alphapdef}
    \alpha(p):=\frac{2dp-1}{2d-1}
\end{equation}
simplifies the analysis. For any direction $\bm{v}\in \{\pm\bm{e}_1, \pm\bm{e}_2, \ldots, \pm\bm{e}_d\}$, let $N_n(\bm{v})$ count the number of steps in direction $\bm{v}$ by time $n$. Given the walk until time $n\geq 1$, the conditional probability that the next step is $\bm{v}$ then equals
    \begin{equation}
    \label{ERWpalpha}
\frac{pN_n(\bm{v})}{n}+\left(1-\frac{N_n(\bm{v})}{n}\right)\frac{1-p}{2d-1}=\frac{\alpha(p) N_n(\bm{v})}{n}+\frac{1-\alpha(p)}{2d}.
    \end{equation}
    In other words, at time step $n\geq 2$, the elephant either replicates a uniformly random historical step (with probability $\alpha(p)$) or takes a fresh step (with probability $1-\alpha(p)$), analogous to a simple symmetric random walk. This structure naturally generalizes to broader step-reinforced random walks by modifying the uniform distribution on $\{\pm\bm{e}_1, \pm\bm{e}_2, \ldots, \pm\bm{e}_d\}$ to a general measure on $\RR^d$.

\begin{definition}[SRRW]
  \label{defSRRW}
  Let $\mu$ be a probability measure on the d-dimensional Euclidean space $\RR^d$. Let $(\xi_n)_{n\geq 2}$ be i.i.d. Bernoulli random variables with success parameter $\alpha\in [0,1]$, and let $(U[n])_{n\geq 1}$ be independent random variables where each $U[n]$ is uniformly distributed on $\{1,2,\ldots,n\}$. Define the step sequence $(\bm{X}_n)_{n\geq 1}$ and the walk recursively as follows:: 
\begin{enumerate}[(i)]
  \item At time $n=1$, sample $\bm{X}_1$ from $\mu$, set $\bm{S}_0:=\bm{0}$ and $\bm{S}_1:=\bm{X}_1$;
  \item For $n\geq 1$, given $\bm{X}_1,\bm{X}_2,\dots,\bm{X}_n$: 
  \begin{itemize}
    \item If $\xi_{n+1}=1$, set $\bm{X}_{n+1}=\bm{X}_{U[n]}$;
    \item If $\xi_{n+1}=0$, sample $\bm{X}_{n+1}$ independently from $\mu$.
  \end{itemize}
 Update $\bm{S}_{n+1}:=\bm{S}_n+\bm{X}_{n+1}$.
\end{enumerate}
The process $\bm{S}=(\bm{S}_n)_{n\in \NN}$ is called a step-reinforced random walk (SRRW) in $\RR^d$ with reinforcement parameter $\alpha$ and step distribution $\mu$. For $d=1$, we often denote the walk by $S$.
\end{definition}
\begin{remark}[Invariance under measurable maps]
\label{imagesrrw}
Let $(\bm{X}_n)_{n\geq 1}$ be as in Definition \ref{defSRRW}. For any Borel measurable function $g: \RR^d \to \RR^k$, by definition, the walk  $\bm{S}^{(g)}=(\bm{S}^{(g)}_n)_{n\in \NN}$ defined by
    $$
\bm{S}^{(g)}_0:=\bm{0}, \quad \quad \bm{S}^{(g)}_n:=\sum_{i=1}^n g(\bm{X}_i), \quad  n\geq 1,
    $$
 is an SRRW in $\RR^k$ with the same parameter $\alpha$ and step distribution $\mu \circ g^{-1}$. In particular, each coordinate of $\bm{S}$ is a one-dimensional SRRW.
\end{remark}

When $\alpha=0$, $\bm{S}$ reduces to a classical random walk with i.i.d. steps distributed as $\mu$. When $\alpha=1$, then $\bm{S}_n=n\bm{X}_1$ for all $n\geq 1$ (all steps replicate $\bm{X}_1$). Note that for discrete $\mu$, one may allow $\alpha$ to be negative as long as $\alpha+(1-\alpha)\mu(\bm{v})\geq 0$ for each possible direction $\bm{v}$, i.e.,
\begin{equation}
  \label{discmuanega}
\alpha\geq \max_{\bm{v}:\mu(\bm{v})>0}\frac{-\mu(\bm{v})}{1-\mu(\bm{v})}.
\end{equation}
For example, a $d$-dimensional ERW with memory parameter $p<1/(2d)$ can also be formulated by parameter $\alpha(p)$ in (\ref{alphapdef}) via the right-hand side of (\ref{ERWpalpha}) where $\alpha(p) \in [1/(2d-1),0)$, representing negative reinforcement. However, we restrict to the case $\alpha\geq 0$ in this work, as negative $\alpha$ lacks an intuitive probabilistic interpretation despite potential extensibility of our results to the case (\ref{discmuanega}).

This paper aims to prove that \emph{under mild moment assumptions, both the radial and angular components of an SRRW exhibit phase transitions at $\alpha=1/2$ in all dimensions, independently of $\mu$}. Detailed statements appear in Section \ref{secmainre}.

\subsection{Main results}
\label{secmainre}

To state our main results, we first formalize key concepts: recurrence/transience and the genuine dimension of an SRRW.

\begin{definition}
\label{defrectran}
    An SRRW $\bm{S}=(\bm{S}_n)_{n\in \NN}$ is said to be recurrent, resp. transient, if  
  $$\mathbb{P}\left(\liminf_{n\to \infty}\|\bm{S}_n\|\leq r\right)=1 \ \text{for some}\ r>0, \quad \text{resp.}\quad \mathbb{P}\left(\lim_{n \to\infty}\|\bm{S}_n\|=\infty\right)=1,$$ 
  where $\|\cdot\|$ denotes the usual Euclidean norm.
\end{definition}
\begin{remark}
If $\alpha=0$, the Hewitt-Savage zero-one law implies that if the first probability is less than $1$, then it is $0$. However, if $\alpha>0$, the two behaviors are not a priori since the SRRW could be non-Markovian. 
\end{remark}

\begin{definition}
  Let $\mu$ be a step distribution on $\RR^d$, if the dimension of $\operatorname{span}\langle\operatorname{supp}(\mu)\rangle$, the span of the support of $\mu$, equals $k\leq d$, then we say $\mu$ (and its associated SRRWS $\bm{S}$) is genuinely $k$-dimensional.
\end{definition}

Without loss of generality, we may assume $\bm{S}$ is genuinely $d$-dimensional. To see this, suppose $\mu$ is genuinely $k$-dimensional. Let $A=(\bm{v}_1,\bm{v}_2,\dots,\bm{v}_k)$ be a $d \times k$ matrix whose columns form a basis for $\operatorname{span}\langle\operatorname{supp}(\mu)\rangle$. There exists a linear map $g: \RR^d \to \RR^k$ such that any $\bm{x} \in \operatorname{span}\langle\operatorname{supp}(\mu)\rangle$ decomposes uniquely as $\bm{x}=Ag(\bm{x})$. Indeed, $g(\bm{x})=(A^TA)^{-1}A^T\bm{x}$ since $A^TA$ is invertible, and thus,
\begin{equation}
    \label{genuinekspace}
    \|g(\bm{x})\|\leq \|(A^TA)^{-1}\|\|A^T\|\|\bm{x}\|, \quad \text{and}\quad \|\bm{x}\|\leq \|A\| \|g(\bm{x})\|,
\end{equation}
where we also use $\|\cdot\|$ for the Frobenius norm of a matrix. By Remark \ref{imagesrrw}, $(g(\bm{S}_n))_{n\in \NN}$ is a genuinely $k$-dimensional step-reinforced random walk in $\RR^k$. By (\ref{genuinekspace}), $\bm{S}$ is recurrent (resp. transient) if and only if $(g(\bm{S}_n))_{n\in \NN}$ is recurrent (resp. transient). 

From now on, we assume that $\alpha<1$ and $\mu \neq \delta_{\bm{0}}$ (Dirac measure at $\bm{0}$). Recall that $\bm{X}_1$ has distribution $\mu$. For $s\geq 1$, we let $A(s)$ be the assumption that 
\begin{equation}
\label{assmps}
  \EE\|\bm{X}_1\|^s<\infty, \quad \EE \bm{X}_1=\bm{0}.
\end{equation}
Notice that if $\mu$ satisfies Assumption $A(2)$, its covariance matrix $\EE (\bm{X}_1\bm{X}_1^{T})$ is well-defined where $\bm{X}_1$ is viewed as a column vector, and that $\mu$ is $d$-dimensional if and only if the rank of $\EE (\bm{X}_1\bm{X}_1^T)$ equals $d$.

The following result shows that the SRRW undergoes a phase transition between recurrence and transience in lower dimensions and is transient for all parameters in higher dimensions.  

\begin{theorem}
\label{mainthm}
  Let $\bm{S}=(\bm{S}_n)_{n\in \NN}$ be a genuinely d-dimensional SRRW in $\RR^d$ with parameter $\alpha$ and step distribution $\mu$. One has:
  \begin{enumerate}[(i)]
      \item Under Assumption $A(2)$,
      \begin{itemize}
      \item  if $d=1$ and $\alpha\leq 1/2$, then $S$ is recurrent;
      \item if $d\geq 1$ and $\alpha>1/2$, then $\bm{S}$ is transient. 
  \end{itemize}
  \item Under Assumption $A(2+\delta)$ where $\delta>0$, 
  \begin{itemize}
      \item if $d=2$ and $\alpha<1/2$, then $\bm{S}$ is recurrent;
      \item  if $d=2$ and $\alpha=1/2$, then a.s.,
$$
\lim_{n\to \infty} \frac{\log \|\bm{S}_n\|}{\log n} = \frac{1}{2};
$$  
      \item if $d\geq 3$ and $\alpha\leq 1/2$, then $\bm{S}$ is transient.
  \end{itemize}
  \end{enumerate}  
\end{theorem}

$\mathbf{Discussion\ for \ d=2}$: When $\alpha=0$, the Chung-Fuchs theorem ensures recurrence under $A(2)$, and for any $\delta \in (0,1)$, transient examples exist under $A(2-\delta)$, see e.g. \cite[Section 5.4]{MR3930614}. The Assumption $A(2+\delta)$ in Theorem \ref{mainthm} (ii) then seems reasonable.
On the other hand, Bertoin conjectured in a personal communication that the recurrence for $\alpha<1/2$ in $d=2$ holds even under $A(2)$. We concur that $A(2)$ suffices. And we believe that with further efforts, the arguments in this work can be used to show that the following assumption is sufficient:
$$\EE (\|\bm{X}_1\|^{2}\max\{\log^K \|\bm{X}_1\|, 0\})<\infty,  \quad \EE \bm{X}_1=\bm{0},$$ 
where $K$ is a (large) positive number. However, the sufficiency of Assumption $A(2)$ likely requires novel techniques.

\begin{definition}
  A lattice $L$ in $\RR^d$ is a discrete additive subgroup of rank $k\leq d$, that is, there exists a collection of linearly independent vectors $\bm{v}_1,\bm{v}_2,\dots,\bm{v}_k \in \RR^d$ such that 
\begin{equation}
    \label{deflattice}
    L=\left\{\sum_{i=1}^k m_i \bm{v}_i: m_1, \ldots, m_k \in \mathbb{Z}\right\}.
\end{equation}
\end{definition}

\begin{corollary}
\label{srrwlattice}
     Let $\bm{S}$ be a genuinely d-dimensional SRRW in $\RR^d$ with parameter $\alpha$ and step distribution $\mu$ supported on a lattice $L$. Define the reachable set
     $$
\mathcal{R}:=\left\{\bm{x} \in L: \mathbb{P}\left(\bm{S}_n=\bm{x}\right)>0 \text { for some } n \geq 0\right\}.
$$
Suppose that one of the following conditions for recurrence in Theorem \ref{mainthm} holds, that is,
\begin{enumerate}[(i)]
    \item $d=1$, $\alpha\leq 1/2$ and $\mu$ satisfies Assumption $A(2)$.
    \item $d=2$, $\alpha <1/2$ and $\mu$ satisfies Assumption $A(2+\delta)$ for some $\delta>0$.
\end{enumerate} 
Then, every site in $\mathcal{R}$ is visited by $\bm{S}$ infinitely often a.s..
\end{corollary}

 Theorem \ref{mainthm} and Corollary \ref{srrwlattice} generalize the prior results on elephant random walks \cite{MR4248774, MR4824917, qin2023recurrence}, and can be used to study elephant random walks on lattices beyond $\ZZ^d$. For example, consider an SRRW in $\RR^2$ with parameter $\alpha$ and step distribution $\mu$ uniform on $\{\pm(2,0), \pm (1,\sqrt{3}), \pm (1,-\sqrt{3})\}$. As in (\ref{ERWpalpha}), the walk corresponds to an ERW on the planar equilateral triangular lattice with memory parameter $p=(5\alpha+1)/6$. Theorem \ref{mainthm} and Corollary \ref{srrwlattice} imply a recurrence-transience transition at $p=7/12$, demonstrating that critical parameters for ERWs depends on both dimension and lattice structure. On the other hand, as shown in Proposition \ref{highdspeed} below, SRRWs have a uniform critical parameter $\alpha=1/2$ even in higher dimensions in terms of the escape speed.

 \begin{proposition}
\label{highdspeed}
    Let $\bm{S}=(\bm{S}_n)_{n\in \NN}$ be a genuinely d-dimensional SRRW in $\RR^d$ $(d\geq 3)$ with parameter $\alpha$ and step distribution $\mu$ satisfying Assumption $A(2+\delta)$ for some $\delta>0$. Then,
$$
\lim_{n\to \infty} \frac{\log \|\bm{S}_n\|}{\log n} = \max\{\alpha,\frac{1}{2}\} \quad \text{a.s..}
$$    
\end{proposition}

A simpler version of Theorem \ref{mainthm} and Proposition \ref{highdspeed} is summarized in Fig. \ref{FigSRRW}. Note that the recurrence, transience, and escape speed described above are properties of the radial component of the SRRW. We show in Corollary \ref{cortransang} below that there is also a phase transition of the angular component at $\alpha=1/2$. 

For $\bm{x}\in \RR^d$, we write $\widehat{\bm{x}}:=\bm{x}/\|\bm{x}\|$ if $\bm{x}\neq \bm{0}$ and write $\widehat{\bm{0}}:=\bm{0}$. For a set of real numbers $E$, we denote its indicator function by $\mathds{1}_E$.

\begin{corollary}[Phase transition of the angular component]
\label{cortransang}
    Let $\bm{S}=(\bm{S}_n)_{n\in \NN}$ be a genuinely d-dimensional SRRW in $\RR^d$ with parameter $\alpha$ and step distribution $\mu$. Then under Assumption $A(2)$, in any dimension $d\geq 1$, one has 
$$
\mathbb{P}\left(\lim _{n \rightarrow \infty} \widehat{\bm{S}}_n  \text{ exists}\right)=\mathds{1}_{\left(1/2, 1\right]}(\alpha).
$$
\end{corollary}

	\begin{figure}[t]
			\centering
			\subfigure[$d=1,2$]{
		  \begin{tikzpicture}[scale=3.5]
			\draw[->] (0,0) -- (1.1,0) node[right] {$\alpha$};
			\draw[->] (0,0) -- (0,1.1) node[above ] {$\displaystyle \lim_{n\to \infty} \frac{\log \|\bm{S}_n\|}{\log n}$};
			\draw[blue,dotted]  (0,1/2) -- (0.5,1/2);
			\draw[red,thick]  (0.5,1/2) -- (1,1);
			\node [below] at (0,0) {$0$};
			\node [below] at (0.5,0) {$\frac{1}{2}$};
			\draw[-] (0.5,0) -- (0.5,0.02);
			\node [below] at (1,0) {$1$};
			\draw[-] (1,0) -- (1,0.02);
			\node [left] at (0,1/2) {$\frac{1}{2}$};
			\node [left] at (0,1) {$1$};
			\draw[-] (0,1) -- (0.02,1);
			\node[above] at (0.25,0) {$\text{\small Recurrent}$};
			\node[above] at (0.75,0) {$\text{\small Transient}$};
		  \end{tikzpicture}
		  \label{Fig1-1}
		  }
		  \subfigure[$d\geq 3$]{
		  \begin{tikzpicture}[scale=3.5]
			\draw[->] (0,0) -- (1.1,0) node[right] {$\alpha$};
			\draw[->] (0,0) -- (0,1.1) node[above] {$\displaystyle \lim_{n\to \infty} \frac{\log \|\bm{S}_n\|}{\log n}$};
			\draw[red,thick]  (0,1/2) -- (0.5,1/2);
			\draw[red,thick]  (0.5,1/2) -- (1,1);
			\node [below] at (0,0) {$0$};
			\node [below] at (0.5,0) {$\frac{1}{2}$};
			\draw[-] (0.5,0) -- (0.5,0.02);
			\node [below] at (1,0) {$1$};
			\draw[-] (1,0) -- (1,0.02);
			\node [left] at (0,1/2) {$\frac{1}{2}$};
			\node [left] at (0,1) {$1$};
			\draw[-] (0,1) -- (0.02,1);
			\node[above] at (0.5,0) {$\text{\small Transient}$};
		  \end{tikzpicture}
		  \label{Fig1-2}
		  }
		  \caption{Recurrence and transience of the SRRW under $2+\delta$-th moment assumption (second moment assumption for $d=1$)}
		  \label{FigSRRW}
		  \end{figure}

As part of the proof of Theorem \ref{mainthm}, we establish the following Marcinkiewicz-Zygmund type strong law of large numbers for the SRRW, which may be of independent interest.

\begin{proposition}
\label{MZslln}
    Let $\bm{S}=(\bm{S}_n)_{n\in \NN}$ be an SRRW in $\RR^d$ with parameter $\alpha \leq 1/2$ and step distribution $\mu$ such that $\EE\|\bm{X}_1\|^{s}<\infty$ for some $s\in (1,2)$. Then, for any $\nu\in (0,1-s^{-1})$, 
    \begin{equation}
        \label{MZsllnequ}
        \lim_{n\to \infty}n^{\nu}\left\|\frac{\bm{S}_n}{n}-\EE \bm{X}_1\right\| = 0 \quad \text{a.s. and in}\ L^1(\PP).
    \end{equation}
\end{proposition}
\begin{remark}
   (i) If a measurable function $g: \RR^d \to \RR^k$ is such that $\EE \|g(\bm{X}_1)\|^s<\infty$ for some $s\in (1,2)$, then by Remark \ref{imagesrrw} and Proposition \ref{MZslln}, for any $\nu\in (0,1-s^{-1})$, 
    $$
\lim_{n\to \infty}n^{\nu}\left\|\frac{\sum_{i=1}^n g(\bm{X}_i)}{n}-\EE g(\bm{X}_1)\right\| = 0 \quad \text{a.s. and in}\ L^1(\PP).
    $$
 (ii) While preparing the revision, we noticed that Aguech, Hariz, Machkouri and Faouzi \cite{aguech2025class} proved independently the a.s.-convergence result in (\ref{MZsllnequ}) through a different approach.
\end{remark}

\subsection{Previous results on reinforced random walks}
\label{secprerrw}

The step-reinforced random walk can be viewed as a member of the family of reinforced random walks. We briefly review some key results on recurrence, transience and localization in related models.

The notion of reinforced random walk was introduced in 1986 by Coppersmith and Diaconis \cite{diaconis1986random} with the seminal definition of edge-reinforced random walk (ERRW). On a locally finite connected graph, the ERRW (resp. vertex-reinforced random walk, or VRRW for short)  assigns weights to edges (resp. vertices), updated linearly with each traversal (resp. each visit); transition probabilities are proportional to these weights. For non-linear reinforcement, see surveys \cite{tarres2011localization}, \cite[Section 1.3]{MR2451570}, \cite[Section 1.2]{collevecchio2024strongly}.

Coppersmith and Diaconis \cite{diaconis1986random} observed that the ERRW is recurrent on $\ZZ$. On an infinite binary tree, Pemantle \cite{MR942765} proved a recurrence/transience phase transition. Merkl and Rolles obtained recurrence criteria on graphs of the form $\ZZ \times G$ where $G$ is a finite graph and on a two-dimensional graph \cite{MR2529432, MR2561431}. Sabot and Tarr\`es \cite{MR3420510} established a link among the ERRW, the vertex-reinforced jump process (VRJP) and a supersymmetric hyperbolic sigma model, and proved positive recurrence of ERRW and VRJP on any graph of bounded degree for sufficiently small initial weights and transience of VRJP on $\ZZ^d$ $(d\geq 3)$. The recurrence results were also obtained by Angel, Crawford and Kozma \cite{MR3189433} using a different approach. Whereas the transience of ERRW on $\ZZ^d$ $(d\geq 3)$ for large initial weights was proved by Disertori, Sabot and Tarr\`es \cite{MR3366053}. Using a representation of the VRJP through a random {S}chr\"{o}dinger operator by Sabot, Tarr\`es and Zeng \cite{MR3729620, MR3904155}, Poudevigne \cite{MR4721024} proved  the uniqueness of the recurrence-transience phase transition of the ERRW and VRJP on any graph of bounded degree, and  Sabot and Zeng \cite{MR3904155, MR4218029} showed the recurrence of ERRW and VRJP on $\ZZ^2$ for any initial constant weights.

On infinite graphs, as a mixture of reversible Markov chains \cite{MR3904155}, the ERRW a.s. visits infinitely many vertices. However, the VRRW exhibits a totally different phenomenon|graph localization, i.e., it can visit only finitely many vertices.  On the integers $\ZZ$, Pemantle and Volkov \cite{MR1733153} showed that the VRRW a.s.
 localizes and, with positive probability, eventually gets stuck on five vertices, and Tarr\`es \cite{MR2078554} proved that this localization on five points is the almost sure behavior. On arbitrary graphs of bounded degree, Volkov \cite{MR1825142} showed that the VRRW localizes with positive probability, which was later generalized by Benaïm and Tarr\`es \cite{MR2932667}.

 Compared to ERRW and VRRW, SRRW has less memory: It only remembers the total number of prior steps in each direction. The time and location of the steps are irrelevant. Note that the dynamics of an SRRW is fully encoded by a measure-valued P\'{o}lya urn process whose general version was introduced and studied by Mailler and Marckert \cite{MR3629870}: Let $\mu$ be a distribution on $\RR^d$. Consider an urn with a (possibly uncountable) set of colors $\mathcal{C}$ consisting of $\operatorname{supp}(\mu)$. Initially, the urn is empty. At the first step, we add a ball whose color is sampled from $\mu$. Later, at any time step, a ball is drawn uniformly at random with replacement. If its color is $\bm{x}$, then add a ball whose color is sampled from $\alpha \delta_{\bm{x}}+(1-\alpha)\mu$ where $\alpha\in [0,1]$. By a slight abuse of notation, we denote the color of the ball added at time $i\geq 1$ by $\bm{X}_i$. Then the composition of the urn at time $n$ is given by the random measure $\sum_{i=1}^n\delta_{\bm{X}_i}$ and it is easy to check that the integrals
 $$
\bm{S}_n:= \int_{\RR^d} \bm{x} d \sum_{i=1}^n\delta_{\bm{X}_i}=\sum_{i=1}^n\bm{X}_i, \quad  n\in \NN,
 $$
 define an SRRW in $\RR^d$ with step distribution $\mu$ and parameter $\alpha$. In particular, as was observed by Baur and Bertoin \cite{baur2016elephant}, a $d$-dimensional ERW is driven by a P\'{o}lya urn process with $2d$ colors. Note that under Assumption $A(2+\delta)$, using a law of large numbers for the measure-valued P\'{o}lya urn process \cite[Theorem 2.11]{MR4609444}, one could prove that there exists a constant $\nu>0$ such that  
$$
        \lim_{n\to \infty}n^{\nu}\left\|\frac{\bm{S}_n}{n}-\EE \bm{X}_1\right\| = 0 \quad \text{a.s. and in}\ L^2(\PP).
$$
If we further assume that $\alpha\leq 1/2$, then the estimate on the rate of the a.s.-convergence could be improved, see the law of the iterated logarithm in Section \ref{secprerrw2}. Proposition \ref{MZslln} treats the situation in which $\EE \|\bm{X}_1\|^2=\infty$ but $\EE \|\bm{X}_1\|^s<\infty$ for some $s\in (1,2)$.

\subsection{Organization of the paper} The remainder of this paper is organized as follows. 

In Section \ref{secintroproof}, we introduce some preliminary notation and results, and outline the main ideas of our proofs. In Section \ref{secpro}, we prove the results stated in Section \ref{secmainre}. More precisely: 
\begin{itemize}
    \item The proof for the superdiffusive regime $(\alpha >1/2)$ is given in Section \ref{RnYnsec}. Corollary \ref{cortransang} is also proved there.
    \item Section \ref{stocappsec} is devoted to the proof of Proposition \ref{MZslln}, which is important for our analysis of the diffusive and critical regimes $(\alpha \leq 1/2)$.
    \item For $\alpha \leq 1/2$, the proof in the cases $d=1$ and $d\geq 3$ is presented in Section \ref{secdiff}; and the proof in the case $d=2$ is more complicated and is postponed to Section \ref{rec2Dproof}.
    \item Based on previous sections, we finish the proofs of Theorem \ref{mainthm}, Corollary \ref{srrwlattice}, and Proposition \ref{highdspeed} in Section \ref{secpromain}.
\end{itemize}

\section{Overview and introduction to the proofs}
\label{secintroproof}

The main techniques we use to prove Theorem \ref{mainthm} are the Lyapunov functions method and a connection of the SRRW to random recursive trees. We shall also use a convergence result of quasi-martingales. To improve readability, we outline the main ideas in this section.

Note that while some of the techniques have already been used in the author's previous work \cite{qin2023recurrence} on the ERW, substantial difficulties arise in the general setting due to the presence of long range jumps. For example, Propositions \ref{prop1dslln} and \ref{Snsio} generalize \cite[Lemma 2.1]{qin2023recurrence} and \cite[Proposition 1.15]{qin2023recurrence}, respectively, but the martingale (resp. supermartingale) used in the proof of \cite[Lemma 2.1]{qin2023recurrence} (resp. \cite[Proposition 1.15]{qin2023recurrence}) is no longer $L^2$-bounded (resp. lower-bounded) for general step distributions, and novel solutions are required.


\subsection{Notation and some preliminary results}  
\label{secnotation}

Let us introduce some preliminary notation. 

\begin{itemize}
    \item We use bold letters to denote vectors in $\RR^d$ and let $x(i)$ denote the $i$-th coordinate of a vector $\bm{x}$. The scalar product of $\bm{x}, \bm{y}\in \RR^d$ is defined to be $\bm{x} \cdot \bm{y}:=\sum_{i=1}^d x(i)y(i)$. We write $x(i)^2:=(x(i))^2$ for simplicity.
    \item  Given $x\in \RR$, we let $\lfloor x \rfloor:=\sup\{n \in \ZZ: n\leq x\}$, and write $x^{-}:=\max \{-x, 0\}$ and $x^{+}:=\max \{x, 0\}$. For $x,y\in \RR$, we write $x\wedge y:=\min\{x,y\}$.
    \item  We let $C(a_1,a_2,\dots,a_k)$ denote a positive constant depending only on real variables $ a_1, a_2, \ldots, a_k$ and let $C$ denote a positive constant. The actual values of these constants may vary from line to line.
    \item Given two sequences $\left(u_n\right)_{n \in \mathbb{N}}$ and $\left(v_n\right)_{n \in \mathbb{N}}$ taking values in $\mathbb{R}$, we write $u_n \sim v_n$ if $u_n/v_n \to 1$ as $n$ goes to infinity.
    \item Given a vector $\bm{h}$ and a $[0,\infty)$-valued function $g$ defined on $[0,\infty)$, we write $\bm{h}=O(g(x))$ as $x\to \infty$, resp. $x\to 0$, if there exist positive constants $C$ and $x_0$ such that $\|\bm{h}\| \leq C g(x)$ for all $x \geq x_0$, resp. $|x|\leq x_0$. Note that $\bm{h}$ may be a function of more variables than just $x$, and that it may depend on $x$ in a complicated way, e.g. $h(x,y) = O(x)$ as $x\to \infty$ if $h(x,y)=(x\sin y, x\cos y)$. 
    \item  We write $X \stackrel{\mathcal{L}}{=} Y$ if two random variables (or processes) $X$ and $Y$ have the same distribution. 
\end{itemize}

From now on, we let $(\FF_n)_{n\geq 1}$ denote the filtration generated by the SRRW $\bm{S}$ we are studying. For $n\geq 1$, we let $\bm{X}_n:=\bm{S}_n-\bm{S}_{n-1}$ denote the $n$-th step of $\bm{S}$. Given a Borel measurable function $g: \RR^d \to \RR^k$, let $\bm{S}^{(g)}$ denote the SRRW defined in Remark \ref{imagesrrw}. Let $(\mu_n)_{n\geq 1}$ denote the empirical measures of $(\bm{X}_n)_{n\geq 1}$, i.e., 
   \begin{equation}
        \label{defmun}
\mu_n:=\frac{\delta_{\bm{X}_1}+\delta_{\bm{X}_2}+\dots+\delta_{\bm{X}_n}}{n}, \quad n\geq 1.
    \end{equation}
For any $n\geq 1$, let $$\Delta_n:=\mu_n-\mu.$$ 
If $\mu$ is a probability measure on $\RR^d$ and $h$ is an integrable function with respect to $\mu$, we write 
$$
\Delta_n(h):= \int_{\bm{x}\in \RR^d} h(\bm{x})d \mu_n (\bm{x})- \int_{\bm{x}\in \RR^d} h(\bm{x})d \mu (\bm{x})=\frac{\bm{S}^{(h)}_n}{n}-\EE h(\bm{X}_1).
$$
If $\mu$ satisfies Assumption $A(2)$, we let $\Delta_n(\bm{x}\bm{x}^T)$ be the $d\times d$ matrix with 
$$\Delta_n(\bm{x}\bm{x}^T)_{ij}:=\Delta_n(x(i)x(j))\quad  i,j\in \{1, 2,\dots,d\}.$$

\subsubsection{Some previous results on the SRRW}
\label{secprerrw2}

We recall some previous results on the SRRW that will be used later. It is common in the literature to distinguish between three regimes: the so-called diffusive regime $\left(\alpha<1/2\right)$, the critical regime $\left(\alpha = 1/2\right)$, and the superdiffusive regime $\left(\alpha>1/2\right)$.

In dimension 1, Donsker’s invariance principle was established in the diffusive and critical regimes \cite{MR4237267, MR4490996}. Hu and Zhang \cite{MR4794980} later established the strong invariance principle, which implies the law of the iterated logarithm \cite[Corollary 1.1]{MR4794980}: If $\EE X_1^2<\infty$, then for $\alpha<1/2$, 
\begin{equation}
  \label{lildiff}
  \limsup _{n \rightarrow \infty} \frac{S_n-n\EE X_1 }{\sqrt{2 n \log \log n}}=\frac{\sqrt{\operatorname{Var}(X_1)}}{\sqrt{1-2 \alpha}} \quad \text{a.s.,}
\end{equation}
and for $\alpha=1/2$, 
\begin{equation}
  \label{lilcrit}
\limsup _{n \rightarrow \infty} \frac{S_n-n\EE X_1 }{\sqrt{2 n \log n \log \log \log n}}=\sqrt{\operatorname{Var}(X_1)} \quad \text{a.s..}
\end{equation}
In any dimension $d\geq 1$, if $\alpha \leq 1/2$ and $\mu$ satisfies Assumption $A(2)$, then by Remark \ref{imagesrrw}, (\ref{lildiff}), (\ref{lilcrit}), and using that 
$$
\max_{1\leq i \leq d} |S_n(i)| \leq \|\bm{S}_n\| \leq \sum_{1\leq i \leq d} |S_n(i)|, \quad n\in \NN,
$$
one has, 
\begin{equation}
  \label{lilSRRW}
    \limsup_{n\to \infty}\frac{\log \|\bm{S}_n\|}{\log n}=\frac{1}{2} \quad \text{a.s..}
 \end{equation}

In all dimensions $d\geq 1$, it has been proved in \cite[Theorem 1.1]{bertenghi2021asymptotic} and \cite[Theorem 2]{MR4237267} that if $\alpha>1/2$ and $\EE \|\bm{X}_1\|^2<\infty$, then 
\begin{equation}
  \label{Wdefsuperd}
  \lim _{n \rightarrow \infty} \frac{\bm{S}_n-n \mathbb{E}\bm{X}_1}{n^{\alpha}}=\bm{W}, \quad \text{a.s. and in}\ L^2,
\end{equation}
where $\bm{W}$ is a non-degenerate random vector. 

The law of large numbers for the SRRW was established in \cite[Theorem 1.1]{MR4490996}: If $\EE \|\bm{X}_1\|^2<\infty$, then for any $\alpha\in [0,1)$,
\begin{equation}
  \label{llnSRRW}
  \lim _{n \rightarrow \infty} \frac{\bm{S}_n}{n}=\EE \bm{X}_1 \quad \text{a.s. and in}\ L^2.
\end{equation} 
The second moment assumption was later relaxed by Hu and Zhang in \cite[Theorem 1.1]{MR4794980} where they showed that the a.s.-convergence in (\ref{llnSRRW}) holds under the minimal moment condition $\EE \|\bm{X}_1\|<\infty$. In particular, if $h$ is an integrable function with respect to $\mu$, then by Remark \ref{imagesrrw},
\begin{equation}
    \label{llnsh}
    \lim_{n\to \infty}\frac{\sum_{i=1}^nh(\bm{X}_i)}{n}=\EE h(\bm{X}_1), \quad \text{a.s..}
\end{equation}
It is worth mentioning that in \cite[Theorem 1.1]{MR4490996} and \cite[Theorem 1.1]{MR4794980}, the step distribution $\mu$ is assumed to be one-dimensional. But one can easily deduce the convergence of $\bm{S}$ from the convergence of its components. 

\subsubsection{Some simple results}
We state some simple results that will be used frequently in the proofs.  We state them as a separate lemma to avoid repeating similar computations.

\begin{lemma}
\label{condlemXn1}
    Let $\bm{S}=(\bm{S}_n)_{n\in \NN}$ be an SRRW in $\RR^d$ with parameter $\alpha$ and step distribution $\mu$ which satisfies Assumption $A(1)$. If $h$ is an integrable function with respect to $\mu$, then for any $n\geq 1$, one has
    \begin{equation}
        \label{condexXn1}
      \EE (\bm{X}_{n+1}\mid \FF_n)=\frac{\alpha \bm{S}_n}{n}, \quad \EE (h(\bm{X}_{n+1}) \mid \FF_n)= \EE h(\bm{X}_1)+ \alpha \Delta_n(h).     
    \end{equation}
    Moreover, for any $n\geq 1$,
    $$
\EE \Delta_n(h)=0,\quad \text{and} \quad \EE |\Delta_n(h)|\leq 2 \EE |h(\bm{X}_1)|.
    $$
\end{lemma}
\begin{proof}
    By definition, 
    $$
 \EE (h(\bm{X}_{n+1}) \mid \FF_n)= \alpha \int_{\bm{x}\in \RR^d} h(\bm{x})d \mu_n (\bm{x})+(1-\alpha) \EE h(\bm{X}_1).
    $$
    which proves the second equality in (\ref{condexXn1}). Now let $h(\bm{x})=x(j)$, $j=1,2,\dots,d$ and use that $\EE \bm{X}_1=\bm{0}$ to obtain the first equality in (\ref{condexXn1}). For $n\geq 1$, by (\ref{condexXn1}), one has 
    \begin{equation}
        \label{ExpSnequ0}
         \EE\bm{S}_n= \EE \bm{S}_1\prod_{i=1}^{n-1} (1+\frac{\alpha}{i})=\bm{0},
    \end{equation}
    with the convention that $\prod^0_{i=1}(1+\alpha/i)=1$. Now observe that $(\bm{S}_n^{(h)}-n\EE h(\bm{X}_1))_{n\in \NN}$ and $(\bm{S}_n^{(|h|)}-n\EE| h(\bm{X}_1)|)_{n\in \NN}$ are two SRRWs with step distributions both satisfying Assumption $A(1)$. Then, one obtains from (\ref{ExpSnequ0}) that
   $$
\EE \Delta_n(h)=\frac{\EE \bm{S}_n^{(h)}-n\EE h(\bm{X}_1)}{n}=0, \quad \EE |\Delta_n(h)|\leq  \frac{\EE \sum_{i=1}^n |h(\bm{X}_i)|}{n}+\EE |h(\bm{X}_1)|=2\EE |h(\bm{X}_1)|.
   $$
\end{proof}

\subsection{Lyapunov functions method}
\label{Lyapunovf}

The Lyapunov functions method was developed by Lamperti \cite{MR0126872}. This method is robust and powerful, which has been used by Menshikov, Popov and Wade \cite{MR3587911} to study non-homogeneous random walks. Their work hugely inspires our proof. 

To prove the recurrence of an SRRW $\bm{S}$, we shall choose a function $f$ which we call the Lyapunov function and estimate the quantity 
\begin{equation}
    \label{condincref}
   \EE(f(\bm{S}_{n+1})-f(\bm{S}_{n}) \mid \FF_n) \quad \text{for} \ \bm{S}_{n} \notin B(0,r), 
\end{equation}
where $B(0,r)$ is a ball of radius $r>0$. This would allow us to control the growth rate of $(f(\bm{S}_{n}))_{n\in \NN}$. To deal with long-range jumps (which are absent in the case of ERWs), we shall use truncation arguments: For $\varepsilon\in (0,1)$ and $\bm{x}\in \RR^d$, we adopt the notation used in \cite[Equation (3.19)]{MR3587911} and let
  \begin{equation}
    \label{Evarsetdef}
      E_{\varepsilon}(\bm{x}):=\left\{\bm{y} \in \RR^d:\|\bm{y}\| \leq\|\bm{x}\|^{1-\varepsilon}\right\}.
  \end{equation}
Then on the event $\{\bm{X}_{n+1}\in E_{\varepsilon}(\bm{S}_n)\}$, $\bm{X}_{n+1}$ is small compared to $\bm{S}_{n}$ and the increment $f(\bm{S}_{n+1})-f(\bm{S}_{n})$ could be estimated by using a Taylor expansion. And if the step distribution $\mu$ satisfies some suitable moment assumptions, we can control $\EE((f(\bm{S}_{n+1})-f(\bm{S}_{n})) \mathds{1}_{\{\bm{X}_{n+1}\notin E_{\varepsilon}(\bm{S}_n)\}}\mid \FF_n)$ by the first-moment method.

For example, we let $\bm{S}$ be a genuinely 2-dimensional SRRW in $\RR^2$ with parameter $\alpha<1/2$ and step distribution $\mu$ which satisfies $A(2+\delta)$ for some $\delta>0$. As will be explained in Section \ref{secpromain}, we can assume that $\EE \bm{X}_1\bm{X}_1^T$ is the identity matrix. Define a function $f$ on $\RR^2$ by
\begin{equation}
  \label{f2Ddef}
  f(\bm{x}) := \left \{ \begin{aligned} &\sqrt{\log \|\bm{x}\|}  && \text{if } \|\bm{x}\|\geq 1, \\ &0 && \text{if } \|\bm{x}\|< 1.  \end{aligned} \right.
\end{equation}

\begin{lemma}
  \label{fxdelta2dest}
  Let $f$ be as in (\ref{f2Ddef}). Then for any $\bm{x}, \bm{y}  \in \RR^2$ with $\bm{x} \neq \bm{0}$ and $\varepsilon \in (0,1)$, one has
\begin{equation}
  \label{fxdeltagener}
  f(\bm{x}+\bm{y})-f(\bm{x}) \leq 1 + \frac{\|\bm{y}\|}{\|\bm{x}\|}. 
\end{equation}
For any $\varepsilon \in (0,1)$, there exist positive constants $r$ and $C$ depending on $\varepsilon$ such that for any $\|\bm{x}\|\geq r$ and $\bm{y} \in E_{\varepsilon}(\bm{x})$, one has
  \begin{equation}
    \label{fxdeltaEvare}
    \begin{aligned}
         f(\bm{x}+\bm{y})-f(\bm{x}) &\leq  \frac{1}{2 \sqrt{\log \|\bm{x}\|}}\left(\frac{ \bm{x} \cdot \bm{y}}{\|\bm{x}\|^2}+\frac{\|\bm{y}\|^2}{2\|\bm{x}\|^2}-\frac{(\bm{x} \cdot \bm{y})^2}{\|\bm{x}\|^4}+\frac{C\|\bm{y}\|^2}{\|\bm{x}\|^{2+\varepsilon}}\right) \\
         &\quad -\frac{1}{8 \log^{3/2} \|\bm{x}\|} \frac{(\bm{x} \cdot \bm{y})^2}{\|\bm{x}\|^4}+\frac{C\|\bm{y}\|^2}{\|\bm{x}\|^{2+\varepsilon}\log^{3/2} \|\bm{x}\|}.
    \end{aligned}
  \end{equation}
\end{lemma}

The proof of Lemma \ref{fxdelta2dest} mainly relies on Taylor expansions and will be given in Section \ref{rec2Dproof}. In dimension $d=1$, a more general result was proved in \cite[Lemma 3.4.1]{MR3587911}. 

 In the Markovian case (i.e. $\alpha=0$), since we have assumed that $\EE \bm{X}_1\bm{X}_1^T$ is the $2\times 2$ identity matrix, we have
\begin{equation}
    \label{CondmarkXn1}
    \EE (\bm{X}_{n+1}\mid \FF_n)=\bm{0}, \ \EE (\| \bm{X}_{n+1}\|^2\mid \FF_n)=2,\ \EE ((\bm{S}_n\cdot \bm{X}_{n+1})^2\mid \FF_n)=\|\bm{S}_n\|^2.
\end{equation}
 Choose a small $\varepsilon$ such that $(1-\varepsilon)(2+\delta)>2$. If $\bm{S}_{n}\neq \bm{0}$, then, by the first moment method, 
$$
\begin{aligned}
    &\quad\ \EE \left(\left(1+\frac{\|\bm{X}_{n+1}\|}{\|\bm{S}_{n}\|}\right) \mathds{1}_{\{\bm{X}_{n+1}\notin E_{\varepsilon}(\bm{S}_n)\}}\mid \FF_n\right) \\
    &\leq \EE \left(\frac{\|\bm{X}_{n+1}\|^{2+\delta}}{\|\bm{S}_{n}\|^{(1-\varepsilon)(2+\delta)}}+\frac{\|\bm{X}_{n+1}\|^{2+\delta}}{\|\bm{S}_{n}\|^{1+(1-\varepsilon)(1+\delta)}}\mid \FF_n\right) \\
    &\leq \frac{\EE \|\bm{X}_{1}\|^{2+\delta} }{\|\bm{S}_{n}\|^{(1-\varepsilon)(2+\delta)}}+ \frac{\EE \|\bm{X}_{1}\|^{2+\delta} }{\|\bm{S}_{n}\|^{1+(1-\varepsilon)(1+\delta)}}.
\end{aligned}
$$
Combined with Lemma \ref{fxdelta2dest}, this would imply that there exists positive constants $\nu$ and $r$ such that if $\bm{S}_{n} \notin B(0,r)$, then
\begin{equation}
    \label{condEfminus18}
     \EE(f(\bm{S}_{n+1})-f(\bm{S}_{n}) \mid \FF_n)\leq \frac{1}{\|\bm{S}_{n}\|^{2} \log^{3/2} \|\bm{S}_n\|}\left(-\frac{1}{8}+O\left(\frac{1}{\|\bm{S}_{n}\|^{\nu}}\right)\right).
\end{equation}
For details, see Section \ref{rec2Dproof}. By possibly choosing a larger $r$, we see that if $\bm{S}_{n} \notin B(0,r)$, then
\begin{equation}
    \label{condEfnega}
    \EE(f(\bm{S}_{n+1})-f(\bm{S}_{n}) \mid \FF_n)<0.
\end{equation}
Now we can reproduce the arguments used in the proof of \cite[Theorem 2.5.2]{MR3587911} to show that $\bm{S}$ visits $B(0,r)$ infinitely often: If $\bm{S}_k \notin B(0,r)$ for some $k>0$, we can define the hitting time $\tau_k:=\inf\{n\geq k: \|\bm{S}_n\|\leq r\}$. Then $(f(\bm{S}_{n\wedge \tau_k}))_{n\geq k}$ is a non-negative supermartingale, whence it a.s. converges. We obtain that $\tau_k$ is a.s. finite since $\limsup_{n\to\infty}f(\bm{S}_n)=\infty$ a.s.. This shows that $\bm{S}$ a.s. returns to $B(0,r)$ after each exit. 

The case $\alpha\in (0,1/2)$ is more complicated. We shall also use the truncation arguments and Taylor's expansion. However, besides the reminder terms caused by truncating long-range jumps (when $\bm{X}_{n+1}\notin E_{\varepsilon}(\bm{S}_n)$) and by approximating $f$ using its Taylor expansion (when $\bm{X}_{n+1}\in E_{\varepsilon}(\bm{S}_n)$), we need to deal with the errors caused by approximating $\mu_n$ using $\mu$. For example, by Lemma \ref{condlemXn1}, (\ref{CondmarkXn1}) would become 
\begin{equation}
    \label{2dcondXn1Sn}
    \begin{aligned}
           &\EE (\bm{X}_{n+1}\mid \FF_n)=\frac{\alpha\bm{S}_n}{n}, \quad \EE (\| \bm{X}_{n+1}\|^2\mid \FF_n)=2+\alpha \Delta_n(\|\bm{x}\|^2),\\
           &\EE ((\bm{S}_n\cdot \bm{X}_{n+1})^2\mid \FF_n)=\|\bm{S}_n\|^2+\alpha \bm{S}_n^T \Delta_n(\bm{x}\bm{x}^T)\bm{S}_n, 
    \end{aligned} 
\end{equation}
    where recall that $\Delta_n(\bm{x}\bm{x}^T)$ is the $2\times 2$ matrix with $\Delta_n(\bm{x}\bm{x}^T)_{ij}=\Delta_n(x(i)x(j))$ for $i,j\in \{1, 2\}$. We shall bound these $\Delta_n's$, or more generally, estimate the convergence rate in the strong law of large numbers as in Proposition \ref{MZslln}. This would enable us to apply stopping time arguments. More specifically, in the proofs of Proposition \ref{2DSRRWprop}(i), we shall introduce four sequences of stopping times which can be non-rigorously interpreted as follows:
    \begin{itemize}
      \item $(n_k)_{k\geq 1}$: At time $n_k$, $\bm{S}$ is in a ``good position," meaning that $\|\bm{S}_{n_k}\|$ and all related $\Delta_{n_k}'s$ are sufficiently small.
      \item $(T_k)_{k\geq 1}$: $T_k$ is the first time after $n_k$ when one of the $\Delta_{n}'s$ exceeds a predefined threshold. We would deduce from Proposition \ref{MZslln} that $T_k$ is a.s. infinite for all large $k$.
      \item $(\tau_k)_{k\geq 1}$: $\tau_k$ is the first time that $\bm{S}$ returns to $B(0,r)$ after time $n_k$ where $r>0$ is some constant.
      \item $(\theta_k)_{k\geq 1}$: $\theta_k$ is the first time after $n_k$ when $\|\bm{S}\|$ exceeds a predefined threshold.
    \end{itemize}
    Then, for $n\in [n_k, \tau_k\wedge T_k \wedge \theta_k)$, we can use the second-order Taylor expansion of $f$ and truncation to estimate the conditional increment $f(\bm{S}_{n+1})-f(\bm{S}_{n})$, in particular, we can show that (\ref{condEfnega}) still holds. We then apply the optional stopping theorem to the supermartingale $(f(\bm{S}_{n\wedge \tau_k\wedge  T_k\wedge \theta_k}))_{n\geq n_k}$, and show that the stopping time $\tau_k \wedge T_k$ is a.s. finite, and thus, $\tau_k$ is a.s. finite for all large $k$. This implies that $\bm{S}$ visits $B(0,r)$ infinitely often. See Section \ref{rec2Dproof} for details. 

Notice that to obtain (\ref{condEfnega}), a first-order Taylor expansion of the Lyapunov function $f$ is not sufficient. To apply a second-order Taylor expansion, we need Assumption $A(2)$. In dimension 2, we further impose Assumption $A(2+\delta)$ to make sure that the errors from approximation and truncation are relatively small compared to the leading term in the right-hand side of (\ref{condEfminus18}), i.e.,
$$   - \frac{1}{8\|\bm{S}_{n}\|^{2} \log^{3/2} \|\bm{S}_n\|}.$$
The leading term itself is small (more precisely, its absolute value is small). In dimension 1, the corresponding leading term is relatively larger so that we can control those errors even under Assumption $A(2)$.

\subsection{Convergence of quasi-martingales}
\label{convquasimasec}

Although the Lyapunov functions method will not be used for the critical case, i.e., $d=2$ and $\alpha=1/2$, it provides a potential alternative proof strategy. Let $\bm{S}$ be a genuinely 2-dimensional SRRW in $\RR^2$ with parameter $\alpha$ and step distribution $\mu$ which satisfies $A(2+\delta)$ for some $\delta>0$. If $\alpha<1/2$, then as explained in Section \ref{Lyapunovf}:
\begin{itemize}
    \item We shall first prove a weaker result (Proposition \ref{Snsio}), that is, there exists a positive constant $\tilde{\kappa}<1/2$ such that $\|\bm{S}_n\|\leq n^{\tilde{\kappa}}$ infinitely often a.s.. This guarantees that $\bm{S}$ could start from ``good initial positions" infinitely often. 
    \item Given this weaker result, we can show that for some large $r>0$, $\|\bm{S}_n\|\leq r$ infinitely often a.s..  
\end{itemize}
In fact, one can check that the arguments used in the second part also apply to the case $\alpha=1/2$. Combined with (\ref{lilSRRW}),
this observation shows that for $\alpha=1/2$, a.s. on the event $\{\lim_{n\to \infty}\|\bm{S}_n\| =\infty \}$, one has 
\begin{equation}
\label{Sna12rate}
    \lim_{n\to \infty}\frac{\log \|\bm{S}_n\|}{\log n}=\frac{1}{2}.
\end{equation}
This motivates us to prove directly that (\ref{Sna12rate}) occurs a.s., which, at first glance, appears to be a much stronger result than transience. We shall let 
  \begin{equation}
    \label{defxnkappa}
    x_n:=\frac{\log (\|\bm{S}_n\|^2+n^{\kappa})}{\log n}, \quad n> 1,
  \end{equation}
  where $\kappa \in (0,1)$ will be chosen later and is used to control the step size $|x_{n+1}-x_n|$. Now it suffices to show that $(x_n)_{n> 1}$ is convergent a.s.. Indeed, assuming the a.s.-convergence of $(x_n)_{n> 1}$, one can deduce from (\ref{lilSRRW}) that $\lim_{n\to \infty}x_n = 1$ a.s., which implies (\ref{Sna12rate}).

  Notice that for any $n>1$,
$$
\kappa \leq x_n \leq 1+ \frac{\log (1+\frac{\|\bm{S}_n\|^2}{n})}{\log n} \leq 1+ \frac{\|\bm{S}_n\|^2}{n \log n}.
$$
  By the $L^2$-convergence in (\ref{llnSRRW}), $(x_n)_{n> 1}$ is $L^1$-bounded. One useful strategy to prove the a.s.-convergence of an $L^1$-bounded adapted process is using the following quasi-martingale convergence theorem, see e.g. \cite[Theorem 9.4]{MR0688144}. The notion of quasi-martingale was first introduced by Donald. L. Fisk \cite{MR0192542}, for its rigorous definition and basic properties, see e.g. \cite[Chapter 2]{MR0688144}.

\begin{theorem}
  \label{asconthmquasi}
  Let $(y_n)_{n\geq 1}$ be a real adapted process with respect to the filtration $(\GG_n)_{n\geq 1}$ that satisfies 
$$  \sum_{n=1}^{\infty} \EE(\left|\EE(y_{n+1} -y_n \mid \mathcal{G}_n)\right|)<\infty\quad \text{and} \quad  \sup_{n\geq 1} \EE(y_n^{-})<\infty.$$
  Then $\left(y_n\right)_{n \geq 1}$ converges a.s. towards an integrable random variable.
\end{theorem}

By Theorem \ref{asconthmquasi} with $(y_n)_{n> 1}=(x_n)_{n> 1}$ and $(\GG_n)_{n\geq 1}=(\FF_n)_{n\geq 1}$, we see that a sufficient condition for the a.s.-convergence of $(x_n)_{n> 1}$ is the following: 
 \begin{equation}
   \label{sumposibd}
   \EE\left(\sum_{n=2}^{\infty} \EE (x_{n+1}-x_n\mid \FF_n)^{-} \right)<\infty.
 \end{equation}
The proof of (\ref{sumposibd}) will be given in Section \ref{sec2dcrit12}.

\subsection{Connection to random recursive trees}
\label{secrrt}
As pointed out by  K\"ursten \cite{MR3652690} and Businger \cite{MR3827299}, the ERW and the shark random swim have a connection to Bernoulli bond percolation on random recursive trees. This connection still holds for the general SRRW and can be described as follows. 

Set $\xi_1:=0$ and let $(\xi_n)_{n\geq 2}$ be i.i.d. Bernoulli random variables with parameter $\alpha$. Let $(\bm{\Theta}_n)_{n\geq 1}$ be i.i.d. $\mu$-distributed random vectors. Given $(\xi_n)_{n\geq 1}$ and $(\bm{\Theta}_n)_{n\geq 1}$, we construct a growing random forest $(\mathscr{F}_n)_{n\geq 1}$ and assign spins on its components as follows. 

At time $n=1$, there is a vertex with label 1. We denote by $\mathscr{F}_1$ the forest with this single vertex. Later, at each time step $n\geq 2$:
\begin{enumerate}[(i)] 
  \item  We add and connect a new vertex labeled $n$ to one of the vertices of $\mathscr{F}_{n-1}$  chosen uniformly at
 random.
 \item  If $\xi_n=0$, the edge connecting the new vertex to the existing vertex is deleted; and if $\xi_n=1$, the edge is kept. We then get a forest with $n$ vertices, which we denote by $\mathscr{F}_n$. 
 \item In each connected component of $\mathscr{F}_n$, we designate the vertex with the smallest label as the root. For $i\in \{1,2,\ldots,n\}$, we denote by $c_{i, n}$ the cluster rooted at $i$ and by $\left|c_{i, n}\right|$ its size, with the convention that $c_{i, n}=\emptyset$ if there is no cluster rooted at $i$. See Fig. \ref{forestT8} for an illustration of $\mathscr{F}_8$ where $|c_{i,8}|=0$ for $i=2,5,6,7,8$. To each cluster $c_{i,n}$, we assign a spin $\bm{\Theta}_i$.  Observe that, for any $n\geq i\geq 1$, the component $c_{i, n}\neq \emptyset$ if and only if $\xi_i=0$. In particular, the root of $c_{i,n}$ and the spin assigned to $c_{i,n}$ do not change as $n$ increases, though $c_{i,n}$ may grow as $n$ increases.
\end{enumerate}

 We define a random walk $\bm{S}=(\bm{S}_n)_{n\in \NN}$ by $\bm{S}_0:=\bm{0}$ and 
 \begin{equation}
     \label{rrtconstrS}
      \bm{S}_n =\sum_{i=1}^n\left|c_{i, n}\right| \bm{\Theta}_i, \quad n\geq 1.
 \end{equation}
One can check by induction that the walk $\bm{S}$ is an SRRW with step distribution $\mu$ and parameter $\alpha$.

\begin{remark}
The procedures above construct an SRRW consecutively. For a fixed $k\geq 2$, to get the expression of $\bm{S}_k$ in the form of (\ref{rrtconstrS}), it is possible to change the order of the procedures: 
  \begin{enumerate}[(I)]
      \item (Random recursive tree). Repeat Procedure (i) until $n$ reaches $k$. One obtains a tree with $k$ vertices.
      \item (Bernoulli bond percolation). Perform Bernoulli bond percolation on the tree obtained in Procedure (I), that is, each edge is deleted with probability $1-\alpha$, independently of the other edges. One then gets a forest with 
(possibly empty) clusters $c_{1,k},c_{2,k},\dots,c_{k,k}$.
      \item (Spin system): As in Procedure (iii), to each cluster $c_{i,k}$, we assign a spin $\bm{\Theta}_i$ where $(\bm{\Theta}_i)_{1\leq i\leq k}$ are i.i.d. $\mu$-distributed random vectors.
  \end{enumerate}
Then, one has $\bm{S}_k \stackrel{\mathcal{L}}{=}\sum_{i=1}^k\left|c_{i, k}\right| \bm{\Theta}_i$.
\end{remark}

\begin{figure}
    \centering
    \begin{tikzpicture}
		\node[circle,draw, minimum size=0.5cm] (A) at  (0,0) {1};
		\node[circle,draw, minimum size=0.5cm] (B) at  (-0.5,1)  {2};
		\node[circle,draw, minimum size=0.5cm] (C) at  (0.5,1)  {3};
		\node[circle,draw, minimum size=0.5cm] (D) at  (-1,2)  {4};
		\node[circle,draw, minimum size=0.5cm] (E) at  (1,2)  {5};
		\node[circle,draw, minimum size=0.5cm] (F) at  (-1.5,3)  {6};
		\node[circle,draw, minimum size=0.5cm] (G) at  (-0.5,3)  {7};
		\node[circle,draw, minimum size=0.5cm] (H) at  (0,2)  {8};
		\draw (A) -- (B);
		\draw[dotted] (A) -- (C);
		\draw[dotted] (B) -- (D);
		\draw (C) -- (E);
		\draw (D) -- (F);
		\draw (D) -- (G);
		\draw (B) -- (H);
		\node[left] at (-0.8,1.4) {$\xi_4=0$};
		\node[right] at (0.3,0.4) {$\xi_3=0$};
		\end{tikzpicture}
    \caption{An illustration of the forest $\mathscr{F}_8$}
    \label{forestT8}
\end{figure}

Now assume that $\mu$ satisfies Assumption $A(2)$ and $\alpha>1/2$. By \cite[Theorem 3.1, Lemma 3.3]{MR3399834} and \cite[Lemma 3, Lemma 4]{MR3827299}, for any $i\geq 1$, a.s. on the event $\{\xi_i=0\}$, the limit
$$\lim_{n\to \infty}\frac{|c_{i, n}|}{n^{\alpha}}=\rho_i$$ 
exists in $(0,\infty)$, and $\sum_{i=1}^{\infty}\EE (\rho_i^2\mathds{1}_{\{\xi_i=0\}})<\infty$ with the convention that $\rho_i:=0$ on $\{\xi_i=1\}$. In particular, one has, a.s.,
$$
\lim_{n\to \infty} \frac{|c_{i, n}|}{n^{\alpha}}=(1-\xi_i)\rho_i.
$$
As already shown in the proof of \cite[Theorem 1]{MR4237267}, one can interchange the limit and infinite summation in the following sense:
$$
\lim_{n\to \infty} \frac{\bm{S}_n}{n^{\alpha}} =  \lim_{n\to \infty} \sum_{i=1}^{\infty} \frac{\left|c_{i, n}\right| }{n^{\alpha}}\bm{\Theta}_i = \sum_{i=1}^{\infty} (1-\xi_i)\rho_i\bm{\Theta}_i \quad \text{in}\ L^2,
$$
with the convention that $|c_{i, n}|=0$ for $i>n$. The random vector $\bm{W}$ in (\ref{Wdefsuperd}) then can be written as a series:
\begin{equation}
  \label{Snlimseries}
  \bm{W}=\lim_{n\to \infty} \frac{\bm{S}_n}{n^{\alpha}} = \sum_{i=1}^{\infty} (1-\xi_i)\rho_i\bm{\Theta}_i,\quad a.s..
\end{equation}

\begin{lemma}
  \label{rnYnneqx}
  Let $(\bm{\Theta}_n)_{n\geq 1}$ be i.i.d. $\mu$-distributed random vectors. Assume $(r_n)_{n\geq 1} \in \ell^2$ has infinitely many nonzero terms. If $\mu$ satisfies Assumption $A(2)$, then $\sum_{n=1}^{\infty} r_n\bm{\Theta}_n$ converges a.s., and for any $\bm{x}\in \RR^d$, $\PP(\sum_{n=1}^{\infty} r_n\bm{\Theta}_n=\bm{x})=0$. 
\end{lemma}

The proof of Lemma \ref{rnYnneqx} will be given in Section \ref{RnYnsec}. Assuming Lemma \ref{rnYnneqx}, we prove the second assertion in Theorem \ref{mainthm} (i).

\begin{proposition}
  \label{RnYnxP0}
  Suppose that $\mu$ satisfies Assumption $A(2)$ and $\alpha>1/2$. Let $\sum_{i=1}^{\infty} (1-\xi_i) \rho_i\bm{\Theta}_i$ be as in (\ref{Snlimseries}). Then,
for any $\bm{x}\in \RR^d$, 
$$\PP(\sum_{i=1}^{\infty} (1-\xi_i)\rho_i\bm{\Theta}_i=\bm{x})=0.$$ 
In particular, 
$$
\PP(\lim_{n\to \infty} \frac{\bm{S}_n}{n^{\alpha}}=\bm{0})=0.
$$
\end{proposition}
\begin{proof}
  Since $\sum_{i=1}^{\infty}\EE(\rho_i^2\mathds{1}_{\{\xi_i=0\}})<\infty$ and $\rho_i$ is positive a.s. on the event $\{\xi_i=0\}$, almost surely, $((1-\xi_i)\rho_i)_{i\geq 1}$ is in $\ell^2$, and has infinitely many nonzero terms.  Notice that $(\bm{\Theta}_i)_{i\geq 1}$ are independent of $(\xi_i)_{i\geq 1}$ and $(\rho_i)_{i\geq 1}$. Then, for any $\bm{x}\in \RR^d$, by Lemma \ref{rnYnneqx},
  $$
\PP(\sum_{i=1}^{\infty} (1-\xi_i) \rho_i\bm{\Theta}_i=\bm{x})=\EE \left(\PP\left(\sum_{i=1}^{\infty} (1-\xi_i) \rho_i\bm{\Theta}_i=\bm{x}\mid (\xi_i)_{i\geq 1}, (\rho_i)_{i\geq 1}\right)\right)=0.
  $$
  \end{proof}

Lastly, it's also worth mentioning that using the random recursive tree construction (\ref{rrtconstrS}), one can relate the SRRW to the broadcasting problem (or the root-bit reconstruction problem) studied by Addario-Berry, Devroye, Lugosi and Velona \cite{MR4386534}. The broadcasting problem consists of estimating the value of $\bm{\Theta}_1$ upon observing the unlabelled version of forest $\mathscr{F}_n$, together with the spin assigned to each component. Note that since the vertex labels are not observed, the identity of the vertex with label $1$ is unknown. In the case of 1-dimensional ERW, one may use the sign of $S_n$ (i.e. the majority) as an estimator, and an upper bound for $\limsup_{n\to \infty}\PP(\operatorname{sgn}(S_n)\neq \operatorname{sgn}(S_1))$ was given in \cite[Theorem 3]{MR4386534}. Here $\operatorname{sgn}(\cdot)$ is the usual sign function.

\section{Proofs of the main results}
\label{secpro}

\subsection{Proof for the superdiffusive regime: Lemma \ref{rnYnneqx}}
\label{RnYnsec}

In this section, we prove Lemma \ref{rnYnneqx}, which finishes the proof for the superdiffusive regime $(\alpha>1/2)$.

\begin{proof}[Proof of Lemma \ref{rnYnneqx}.]
Define $(M_n(r))_{n\geq 1}$ by 
$$
M_n(r):=\sum_{i=1}^{n} r_i\bm{\Theta}_i, \quad n\geq 1.
$$
Then, by our assumption, $(M_n(r))_{n\geq 1}$ is an $L^2$-bounded martingale, and thus converges a.s.. 

We now prove the second assertion. Suppose to the contrary that for some $\bm{x}\in \RR^d$, $\PP(\sum_{n=1}^{\infty} r_n\bm{\Theta}_n=\bm{x})=q_1>0$. Now let $(\bm{\Theta}_n^{\prime})_{n\geq 1}$ are i.i.d. $\mu$-distributed random variables independent of $(\bm{\Theta}_n)_{n\geq 1}$. Then
  $$(\bm{\Theta}_1,\bm{\Theta}_2,\dots,\bm{\Theta}_{n-1},\bm{\Theta}_n^{\prime},\bm{\Theta}_{n+1},\dots)\stackrel{\mathcal{L}}{=} (\bm{\Theta}_n)_{n\geq 1}.$$
Since $\mu$ satisfies $A(2)$ and $\mu\neq \delta_{\bm{0}}$, $\bm{\Theta}_1-\bm{\Theta}_1^{\prime}$ has positive variance, and thus, $\PP(\bm{\Theta}_1-\bm{\Theta}_1^{\prime}\neq \bm{0})>0$. By the continuity of probability measures, there exist two positive constants $\varepsilon, q_2$ such that for any $n\geq 1$, 
$$
\PP(\|\bm{\Theta}_n^{\prime}-\bm{\Theta}_n\|>2\varepsilon)\geq 2q_2.
$$
For any $\bm{y}\in \RR^d$, the union bound gives
$$
2 q_2\leq \PP(\|\bm{\Theta}_n-\bm{\Theta}_n^{\prime}\|>2\varepsilon)\leq \PP(\|\bm{\Theta}_n-\bm{y}\|>\varepsilon)+\PP(\|\bm{\Theta}_n^{\prime}-\bm{y}\|>\varepsilon), \quad n\geq 1.
$$
Therefore, using that $\bm{\Theta}_n$ and $\bm{\Theta}_n^{\prime}$ are independent and identically distributed, one has 
$$
\PP(\|\bm{\Theta}_n^{\prime}-\bm{\Theta}_n\|>\varepsilon \mid \bm{\Theta}_n) \geq q_2.
$$

  By Markov's inequality, we can choose a positive constant $C_1>1$ such that 
  $$\PP(\|\bm{\Theta}_n\|\leq C_1/2)\geq 1-\frac{\min\{q_1,q_2\}}{2}.$$ 
  In particular, for any $n\geq 1$, $\PP(\sum_{i=1}^{\infty} r_i\bm{\Theta}_i=\bm{x},\|\bm{\Theta}_n\|\leq C_1/2)\geq q_1/2$. Then, for any $n\geq 1$ with $r_n\neq 0$, we have
  $$
  \begin{aligned}
    &\quad\ \PP(\varepsilon|r_n|<\|\sum_{i=1}^{\infty} r_i\bm{\Theta}_i-\bm{x}\|\leq C_1|r_n| )\\
    &=\PP(\varepsilon |r_n|<\|r_n\bm{\Theta}_n^{\prime}+\sum_{i\neq n} r_i\bm{\Theta}_i-\bm{x}\| \leq C_1|r_n| ) \\
    &\geq \frac{q_1}{2}\PP(\varepsilon|r_n|<\|r_n\bm{\Theta}_n^{\prime}+\sum_{i\neq n} r_i\bm{\Theta}_i-\bm{x}\| \leq C_1|r_n| \ \mid \ \sum_{i=1}^{\infty} r_i\bm{\Theta}_i=\bm{x},\|\bm{\Theta}_n\|\leq \frac{C_1}{2})\\
    &= \frac{q_1}{2}\PP(\varepsilon <\|\bm{\Theta}_n^{\prime} -\bm{\Theta}_n\|\leq C_1 \mid \sum_{i=1}^{\infty} r_i\bm{\Theta}_i=\bm{x},\|\bm{\Theta}_n\|\leq \frac{C_1}{2})\\
    &\geq \frac{q_1\EE \left( \PP\left(\|\bm{\Theta}_n^{\prime} -\bm{\Theta}_n\|>\varepsilon, \|\bm{\Theta}_n^{\prime}\|\leq C_1/2 \mid \bm{\Theta}_n\right) \mathds{1}_{\{\sum_{i=1}^{\infty} r_i\bm{\Theta}_i=\bm{x},\|\bm{\Theta}_n\|\leq C_1/2\}}\right)}{2\PP(\sum_{i=1}^{\infty} r_i\bm{\Theta}_i=\bm{x},\|\bm{\Theta}_n\|\leq C_1/2)},
  \end{aligned}
  $$
  where in the last step we first conditioned on the sigma-algebra generated by $(\bm{\Theta}_i)_{i\in \NN}$ and then used that $(\bm{\Theta}_n,\bm{\Theta}_n^{\prime})$ are independent of $(\bm{\Theta}_i)_{i\neq n}$. Now, by our choice of $C_1$, the inclusion-exclusion principle gives
  $$
\PP\left(\|\bm{\Theta}_n^{\prime} -\bm{\Theta}_n\|>\varepsilon, \|\bm{\Theta}_n^{\prime}\|\leq \frac{C_1}{2} \mid \bm{\Theta}_n\right)\geq q_2 + \PP(\|\bm{\Theta}_n^{\prime}\|\leq \frac{C_1}{2} )-1 \geq \frac{q_2}{2}.
  $$
   Thus, $\PP(\varepsilon |r_n|<\|\sum_{i=1}^{\infty} r_i\bm{\Theta}_i-\bm{x}\| \leq C_1|r_n| )\geq q_1q_2/4$. As $(r_n)_{n\geq 1}$ converges to 0 and has infinitely many nonzero terms, we can find a subsequence $(r_{n_k})_{k\geq 1}$ such that $C_1|r_{n_{k+1}}|<\varepsilon|r_{n_k}|$ for all $k\geq 1$. Hence,
  $$
  1=\PP(\sum_{i=1}^{\infty} r_i\bm{\Theta}_i\in \RR^d) \geq \sum_{k=1}^{\infty}\PP(\varepsilon |r_{n_k}|<\|\sum_{i=1}^{\infty} r_i\bm{\Theta}_i-\bm{x}\|\leq C_1|r_{n_k}|) \geq  \sum_{k=1}^{\infty} \frac{q_1q_2}{4}=\infty,
  $$
  a contradiction.
\end{proof}

\begin{proof}[Proof of Corollary \ref{cortransang}.]
    For $\alpha>1/2$, by Proposition \ref{RnYnxP0}, the random vector $\bm{W}$ is a.s. not the zero vector, and thus, a.s. 
    $$
   \lim_{n\to \infty} \widehat{\bm{S}}_n= \frac{\bm{W}}{\|\bm{W}\|}.
    $$
    
  Now assume that $\alpha\leq 1/2$. By our assumption, components of $\bm{S}$ (and their opposites) are 1-dimensional SRRWs and their step distributions are non-degenerate and satisfy Assumption $A(2)$. Thus, for each $i\in \{1,2,\dots,d\}$, Equations (\ref{lildiff}) and (\ref{lilcrit}) yield that a.s.,
    $$
\liminf_{n\to \infty}S_n(i)=-\infty, \quad \limsup_{n\to \infty}S_n(i)=\infty,
    $$
whence $\lim_{n\to \infty}\widehat{\bm{S}}_n$ does not exist a.s..
    
\end{proof}

 \subsection{Marcinkiewicz-Zygmund type strong law of large numbers: Proposition \ref{MZslln}}
  \label{stocappsec}

As a preparation for the proof for the diffusive and critical regimes $(\alpha\leq 1/2)$, we prove Proposition \ref{MZslln} in this section. We first prove its one-dimensional version. 

If $h$ is an integrable function with respect to $\mu$, recall from Section \ref{secnotation} that we use the following notation:
\begin{equation}
    \label{defDeltanh}
    \Delta_n(h)= \frac{\bm{S}^{(h)}_n}{n}-\EE h(\bm{X}_1), \quad n\geq 1.
\end{equation}

\begin{proposition}
\label{prop1dslln}
 Let $\bm{S}=(\bm{S}_n)_{n\in \NN}$ be an SRRW in $\RR^d$ with parameter $\alpha \leq 1/2$ and step distribution $\mu$. Suppose that $h$ is a measurable function on $\RR^d$ such that $\EE|h(\bm{X}_1)|^{s}<\infty$ for some $s\in (1,2)$. Let $\nu\in (0,1-s^{-1})$ be a constant. Then, 
    \begin{equation}
        \label{asLicon1dim}
        \lim_{n\to \infty}n^{\nu}\Delta_n(h) = 0 \quad \text{a.s. and in}\ L^1(\PP).
    \end{equation}
    Moreover, 
  \begin{equation}
    \label{pnkIg1}
     \lim_{k\to \infty}\PP(\sup_{n> k}n^{\frac{\nu}{2}}|\Delta_n(h)|\geq 1\mid \FF_k)\mathds{1}_{\{|\Delta_k(|h|^{s})|<1, k^{\nu}|\Delta_k(h)|<1\}} =0.
  \end{equation}
\end{proposition}
\begin{remark}
    Roughly speaking, (\ref{pnkIg1}) says that if $|\Delta_k(|h|^{s})|$ is bounded and $\Delta_k(h)$ is small, then with high probability, $\Delta_n(h)$ remains small for all $n\geq k$.
\end{remark}

For the proof of Proposition \ref{prop1dslln}, we prove the conditional version of Lemma \ref{condlemXn1}. We shall use the following notation: For $n\geq 1$, write
$$ \gamma_n:=\frac{1-\alpha}{n+1}, \quad \text{and}\quad \beta_n:=\prod_{k=1}^{n-1}(1-\gamma_k),$$ 
with the convention that $\beta_1=1$. Notice that by the properties of gamma functions, one has
 \begin{equation}
   \label{betanasy}
    \lim_{n\to \infty}\beta_n n^{1-\alpha} =  \lim_{n\to \infty}\frac{\Gamma(n+\alpha) n^{1-\alpha}}{\Gamma(1+\alpha)\Gamma(n+1)}= c(\alpha),
 \end{equation}
 where $c(\alpha)$ is a positive constant.
 
\begin{lemma}
\label{lemhconrec}
 Let $\bm{S}=(\bm{S}_n)_{n\in \NN}$ be an SRRW in $\RR^d$ with parameter $\alpha$ and step distribution $\mu$ which satisfies Assumption $A(1)$. If $h$ is an integrable function with respect to $\mu$, then
  \begin{equation}
 \label{deltanhxformu}
    \Delta_n(h)=\beta_n\left(h(\bm{X}_1)+ \sum_{j=1}^{n-1}\frac{\gamma_j}{\beta_{j+1}}\epsilon_{j+1}(h)\right), \quad n\geq 1,
 \end{equation}
where $(\varepsilon_{j+1}(h))_{j\geq 1}$ is a martingale difference sequence defined by
\begin{equation}
    \label{defvareh}
    \varepsilon_{j+1}(h):=\frac{h(\bm{X}_{j+1})-\EE h(\bm{X}_{1})-
 \alpha\Delta_j(h)}{1-\alpha}, \quad j\geq 1.
\end{equation}
In particular, for any $n> k\geq 1$, 
        $$
\EE (\Delta_n(h) \mid \FF_k)=\frac{\beta_n}{\beta_k}\Delta_k(h).
    $$
\end{lemma}
\begin{proof}
By (\ref{defDeltanh}), for any $n\geq 1$,
\begin{equation}
\label{recudelth}
    \begin{aligned}
    \Delta_{n+1}(h)- \Delta_{n}(h)&=\frac{\sum_{i=1}^nh(\bm{X}_i)+h(\bm{X}_{n+1})}{n+1}-\frac{1}{n+1}\left(1+\frac{1}{n}\right)\sum_{i=1}^nh(\bm{X}_i)\\
    &=\frac{1}{n+1}\left(-\frac{\sum_{i=1}^nh(\bm{X}_i)}{n}+h(\bm{X}_{n+1})\right) \\
    &=\frac{1}{n+1}\left(-\Delta_{n}(h)+h(\bm{X}_{n+1})-\EE h(\bm{X}_{1}) \right) \\
    &= \frac{1-\alpha}{n+1}\left(-\Delta_{n}(h)+\frac{h(\bm{X}_{n+1})-\EE h(\bm{X}_{1})-\alpha\Delta_{n}(h)}{1-\alpha}\right) \\
    &=\frac{1-\alpha}{n+1}\left(-\Delta_{n}(h)+ \varepsilon_{n+1}(h)\right).
\end{aligned}
\end{equation}
That $(\varepsilon_{j+1}(h))_{j\geq 1}$ is a martingale difference sequence follows from (\ref{condexXn1}) in Lemma \ref{condlemXn1}. By induction, one can easily deduce (\ref{deltanhxformu}) from (\ref{recudelth}). For $n>k\geq 1$, (\ref{deltanhxformu}) implies that 
$$
 \Delta_n(h)=\beta_n\left(\frac{\Delta_{k}(h)}{\beta_{k}}+ \sum_{j=k}^{n-1}\frac{\gamma_j}{\beta_{j+1}}\epsilon_{j+1}(h)\right),
$$
whence the last assertion follows.
\end{proof}

The proof of Proposition \ref{prop1dslln} also needs the following known results.

\begin{lemma}
\label{lemmartvonchow}
    Let $(M_n)_{n\geq 1}=(\sum_{i=1}^n Y_i)_{n\geq 1}$ be a martingale with respect to a filtration $(\GG_n)_{n\geq 1}$. Let $s\in [1,2]$ be a constant. One has: 
  \begin{enumerate}[(i)]
      \item (von-Bahr-Esseen inequality \cite{MR170407}) There exists a positive constant $C(s)$ such that for any $n\geq 1$, 
      $$
      \EE |M_n|^s \leq C(s) \sum_{i=1}^n \EE |Y_i|^s.
      $$
     \item (Chow \cite{MR182040}, see also \cite[Theorem 2.18]{MR0624435}) Let $(a_n)_{n\geq 1}$ be an adapted non-decreasing sequence of positive random variables. Then, a.s. on the event 
$$
\left\{\lim_{n\to \infty}a_n=\infty, \sum_{i=1}^{\infty}\frac{\EE(|Y_{i+1}|^s\mid \GG_i)}{a_{i}^s}<\infty\right\},
$$
$\sum_{i=1}^n Y_i/a_i$ converges and $M_n/a_n\to 0$ as $n\to \infty$.  
\end{enumerate}  
\end{lemma}

\begin{proof}[Proof of Proposition \ref{prop1dslln}.]
Recall from Lemma \ref{lemhconrec} that 
  \begin{equation}
 \label{2nddeltanhxformu}
  \Delta_n(h)=\beta_n\left(h(\bm{X}_1)+ \sum_{j=1}^{n-1}\frac{\gamma_j}{\beta_{j+1}}\epsilon_{j+1}(h)\right), \quad n\geq 1.
  \end{equation}
where $(\varepsilon_{j+1}(h))_{j\geq 1}$ is a martingale difference sequence defined in (\ref{defvareh}). We now estimate the $s$-th absolute moments of $(\epsilon_{n+1}(h))_{n\geq 1}$ under the assumption that $\EE|h(\bm{X}_1)|^{s}<\infty$. First note that for any real numbers $y_1,y_2,\dots,y_m$ $(m\geq 1)$, one has,
\begin{equation}
    \label{Crinen2}
    (y_1+y_2+\dots+y_m)^s \leq m^{s-1} (|y_1|^s+|y_2|^s+\dots+|y_m|^s).
\end{equation}
Thus, there exists a positive constant $C(s,\alpha)$ such that for any $n\geq 1$,
\begin{equation}
    \label{eesphsbd}
    \EE |\epsilon_{n+1}(h)|^s \leq C(s,\alpha) \left(\EE |h(\bm{X}_{n+1})|^s+|\EE h(\bm{X}_{1})|^s+\EE |\Delta_n(h)+\EE h(\bm{X}_{1})|^s \right).
\end{equation}
By H\"older's inequality, 
$$
|\Delta_n(h)+\EE h(\bm{X}_{1})|^s\leq \left(\frac{\sum_{i=1}^n|h(\bm{X}_{i})|}{n}\right)^s \leq \frac{\sum_{i=1}^n|h(\bm{X}_{i})|^s}{n}.
$$
Lemma \ref{condlemXn1} (with $h=|h|^s$ in the notation there) shows that for any $n\geq 1$,
\begin{equation}
    \label{estEvars}
    \EE |\epsilon_{n+1}(h)|^s\leq 3C(s,\alpha) \EE |h(\bm{X}_{1})|^s.
\end{equation}
Define a martingale $(M_n(\nu,h))_{n\geq 1}$ by $M_1(\nu,h):=0$ and 
 $$
 M_n(\nu,h):=\sum_{j=1}^{n-1}\frac{\gamma_j}{j^{1-\alpha-\nu}\beta_{j+1}}\epsilon_{j+1}(h), \quad n>1.
 $$
 Recall from (\ref{betanasy}) that $\beta_n\sim c(\alpha)n^{\alpha-1}$, where $c(\alpha)$ is a positive constant. Then by using (\ref{estEvars}), we obtain that
$$
\sum_{n=1}^{\infty}\frac{1}{n^{s(1-\nu)}}\left(\frac{\gamma_n}{n^{-\alpha}\beta_{n+1}}\right)^s\EE|\epsilon_{n+1}(h)|^s <\infty,
$$
where we also used the assumption that $(1-\nu)s>1$. By Lemma \ref{lemmartvonchow} (ii) with $a_n=n^{1-\alpha-\nu}$ (note that $\nu<1-\alpha$), a.s., the martingale $(M_n(\nu,h))_{n\geq 1}$ converges and
$$
\lim_{n\to \infty} n^{\nu}\beta_n\sum_{j=1}^{n-1}\frac{\gamma_j}{\beta_{j+1}}\epsilon_{j+1}(h)=\lim_{n\to \infty} \frac{\beta_nn^{1-\alpha}}{n^{1-\alpha-\nu}}\sum_{j=1}^{n-1}\frac{\gamma_j}{\beta_{j+1}}\epsilon_{j+1}(h)=0,
$$
which, combined with (\ref{2nddeltanhxformu}) and that $\beta_n\sim c(\alpha)n^{\alpha-1}$, completes the proof of the a.s.-convergence in (\ref{asLicon1dim}). Applying Inequality (\ref{Crinen2}) and Lemma \ref{lemmartvonchow} (i) to (\ref{2nddeltanhxformu}), we obtain that
$$
    n^{\nu} \EE |\Delta_n(h)|^s\leq C(s)n^{\nu}\beta_n^s \EE |h(\bm{X}_1)|^s+ C(s)\beta_n^s\sum_{j=1}^{n-1}\left(\frac{\gamma_j}{\beta_{j+1}}\right)^s\EE |\epsilon_{j+1}(h)|^s,
$$
where $C(s)$ is a positive constant. By (\ref{betanasy}) and (\ref{estEvars}), the two terms on the right-hand side converge to $0$ as $n\to \infty$, which proves the $L^1(\PP)$ convergence in (\ref{asLicon1dim}). 

We now prove (\ref{pnkIg1}). Fix $k$, assume that $|\Delta_k(|h|^{s})|<1$, Lemma \ref{lemhconrec} (with $h=|h|^s$ in the notation there) implies that $|\EE (\Delta_n(|h|^{s})\mid \FF_k)|<1$ for all $n\geq k$. In particular, by (\ref{condexXn1}), for all $n>k$,
\begin{equation}
    \label{estcondhsk}
     \EE( |h(\bm{X}_{n})|^s \mid \FF_k) =\EE |h(\bm{X}_{1})|^s + \alpha\EE (\Delta_{n-1}(|h|^{s})\mid \FF_k) \leq \EE |h(\bm{X}_{1})|^s+1,
\end{equation}
and thus, by the definition of $\Delta_k(|h|^{s})$ (see (\ref{defDeltanh})),
\begin{equation}
    \label{estcondavarhsk}
\begin{aligned}
   \frac{\EE (\sum_{i=1}^n|h(\bm{X}_{i})|^s\mid \FF_k) }{n}&\leq \frac{\sum_{i=1}^k |h(\bm{X}_{i})|^s }{n}+\frac{n-k}{n}(\EE |h(\bm{X}_{1})|^s+1) \\
   &= \frac{k}{n}(\Delta_k(|h|^{s})+\EE |h(\bm{X}_{1})|^s)+\frac{n-k}{n}(\EE |h(\bm{X}_{1})|^s+1) \\
   &\leq  \EE |h(\bm{X}_{1})|^s+1.
\end{aligned}
\end{equation}
Similarly as in (\ref{estEvars}), we can use Ineqaulities (\ref{Crinen2}), (\ref{estcondhsk}), (\ref{estcondavarhsk}) and H\"older's inequality to obtain that, for any $n\geq k$,
 \begin{equation}
     \label{bdexpespsfk}
   \begin{aligned}
     &\quad\ \EE (|\epsilon_{n+1}(h)|^s \mid \FF_k) \\
     &\leq C(s,\alpha) \left(\EE (|h(\bm{X}_{n+1})|^s\mid \FF_k)+|\EE h(\bm{X}_{1})|^s+\EE (|\Delta_n(h)+\EE h(\bm{X}_{1})|^s \mid \FF_k)\right)\\
     &\leq 3 C(s,\alpha) (\EE |h(\bm{X}_{1})|^s+1),
 \end{aligned}  
 \end{equation}
where $C(s,\alpha)$ is a positive constant. 

For any $n> k\geq 1$ with $n>k$, by (\ref{2nddeltanhxformu}),
\begin{equation}
 \label{deltankformula}
 \Delta_n(x)=\beta_n\left(\frac{\Delta_{k}(x)}{\beta_{k}}+ \sum_{j=k}^{n-1}\frac{\gamma_j}{\beta_{j+1}}\epsilon_{j+1}(h)\right).
\end{equation}
If $k^{\nu}|\Delta_k(h)|<1$, by (\ref{betanasy}), there exist a positive constant $C(\alpha)$ such that for all $n>k$,
$$
\frac{\beta_n|\Delta_{k}(x)|}{\beta_{k}} \leq \frac{C(\alpha) k^{1-\alpha-\nu}}{n^{1-\alpha}} \leq \frac{C(\alpha)}{n^{\frac{\nu}{2}}} n^{-\frac{\nu}{2}}.
$$
We may choose a large $k$ such that $C(\alpha)k^{-\nu/2}<1/2$.
Therefore, to prove (\ref{pnkIg1}), it suffices to prove that on the event $\{|\Delta_k(|h|^{s})|<1, k^{\nu}|\Delta_k(h)|<1\}$,
\begin{equation}
  \label{pnkIgsec1}
   \lim_{k\to \infty} \PP\left(\sup_{n> k}n^{\frac{\nu}{2}}\beta_n\left|\sum_{j=k}^{n-1}\frac{\gamma_j}{\beta_{j+1}}\epsilon_{j+1}(h)\right|\geq \frac{1}{2} \mid \FF_k\right) =0.
\end{equation}
Define a martingale $(M_n(\frac{\nu}{2},h))_{n\geq k}$ by $M_k(\frac{\nu}{2},h):=0$ and 
 $$
 M_n(\frac{\nu}{2},h):=\sum_{j=k}^{n-1}\frac{\gamma_j}{j^{1-\alpha-\frac{\nu}{2}}\beta_{j+1}}\epsilon_{j+1}(h), \quad n>k.
 $$
Using Lemma \ref{lemmartvonchow} (i) again, we have, for all $n> k$,
\begin{equation}
  \label{MnnugLdelest}
  \begin{aligned}
  \EE\left(\left|M_n(\frac{\nu}{2},h)-M_k(\frac{\nu}{2},h)\right|^{s}\mid \FF_k\right)&\leq C(s) \sum_{j=k}^{n-1} \EE \left(\left|\frac{\gamma_j}{j^{1-\alpha-\frac{\nu}{2}}\beta_{j+1}}\epsilon_{j+1}(h)\right|^{s}\mid \FF_k\right) \\
  &\leq C(s) \sum_{j=k}^{\infty}\left(\frac{\gamma_j}{j^{1-\alpha-\frac{\nu}{2}}\beta_{j+1}}\right)^s \EE \left(|\epsilon_{j+1}(h)|^{s}\mid \FF_k\right).
\end{aligned}
\end{equation}
By (\ref{bdexpespsfk}), $\EE (|\epsilon_{j+1}|^{s}\mid \FF_k)$ is bounded for $j\geq k$. Therefore, the last sum in (\ref{MnnugLdelest}) is finite and converges to 0 as $k\to \infty$ (again, note that $(1-\nu/2)s>1$ by the choice of $\nu$). Consequently, $(M_n(\frac{\nu}{2},h))_{n\geq k}$ converges a.s. to some random variable $M_{\infty}(\frac{\nu}{2},h)$. Write
$$
\delta_n:=M_n(\frac{\nu}{2},h)-M_{\infty}(\frac{\nu}{2},h), \quad n\geq k.
$$
Doob's $L^p$ inequality for martingales and (\ref{MnnugLdelest}) imply that for any $\varepsilon>0$, 
\begin{equation}
    \label{dooblpdeltan}
    \lim_{k\to \infty}\PP(\sup_{n\geq k}|\delta_n|\geq \varepsilon \mid \FF_k) = 0.
\end{equation}
 For any $n> k$, summation by parts gives
$$
\begin{aligned}
   &\quad\ n^{\frac{\nu}{2}}\beta_n\sum_{j=k}^{n-1}\frac{\gamma_j}{\beta_{j+1}}\epsilon_{j+1}\\
  &= n^{\frac{\nu}{2}}\beta_n\sum_{j=k}^{n-1}j^{1-\alpha-\frac{\nu}{2}}\left(\delta_{j+1}-\delta_j\right)\\
  &=\beta_nn^{1-\alpha} \delta_n-n^{\frac{\nu}{2}}\beta_n\left(k^{1-\alpha-\frac{\nu}{2}}\delta_k+\sum_{j=k}^{n-1}\left((j+1)^{1-\alpha-\frac{\nu}{2}}-j^{1-\alpha-\frac{\nu}{2}}\right)\delta_{j+1}\right).
\end{aligned}
$$
By (\ref{betanasy}), we can choose $\varepsilon>0$ such that $\beta_nn^{1-\alpha} \varepsilon<1/6$ for all $n\geq 1$. Then, on the event $\{\sup_{n\geq k}|\delta_n|\geq \varepsilon\}$, for all $n>k$,
\begin{equation}
    \label{bdnnu2betasumk}
    \begin{aligned}
    \left|n^{\frac{\nu}{2}}\beta_n\sum_{j=k}^{n-1}\frac{\gamma_j}{\beta_{j+1}}\epsilon_{j+1}\right|&\leq \beta_nn^{1-\alpha} \varepsilon+ n^{\frac{\nu}{2}}\beta_n\left(k^{1-\alpha-\frac{\nu}{2}}\varepsilon+\sum_{j=k}^{n-1}((j+1)^{1-\alpha-\frac{\nu}{2}}-j^{1-\alpha-\frac{\nu}{2}})\varepsilon\right) \\
    &\leq 3\beta_nn^{1-\alpha} \varepsilon<\frac{1}{2}.
\end{aligned}
\end{equation}
Then (\ref{pnkIgsec1}) follows from (\ref{dooblpdeltan}).
\end{proof}

We are now ready to prove Proposition \ref{MZslln}. Recall that $x(i)$ denotes the $i$-th component of a vector $\bm{x}$.

\begin{proof}[Proof of Proposition \ref{MZslln}.]
Observe that
   $$ n^{\nu}\left\|\frac{\bm{S}_n}{n}-\EE \bm{X}_1\right\| \leq n^{\nu}\sum_{i=1}^d \left|\frac{S_n(i)}{n}-\EE X_1(i)\right|=n^{\nu}\sum_{i=1}^d |\Delta_n(x(i))|.$$
It remains to apply Proposition \ref{prop1dslln} with $h(\bm{x})=x(i)$, $i=1,2,\dots,d$.
\end{proof}

We end this section with two results that will be used frequently in Section \ref{secdiff}.

Recall from Section \ref{secnotation} that if the step distribution $\mu$ is a probability measure on $\RR^d$ satisfying Assumption $A(2)$, then $\Delta_n(\bm{x}\bm{x}^T)$ is a $d\times d$ matrix given by 
$$\Delta_n(\bm{x}\bm{x}^T)_{ij}:=\Delta_n(x(i)x(j))\quad  i,j\in \{1, 2,\dots,d\}.$$
We use $\|\Delta_n(\bm{x}\bm{x}^T)\|$ to denote its Frobenius norm. The following result will be used in Section \ref{secdiff} for second-order Taylor expansions.

\begin{corollary}
  \label{limx2Delta}
Let $\bm{S}=(\bm{S}_n)_{n\in \NN}$ be an SRRW in $\RR^d$ with parameter $\alpha \leq 1/2$ and step distribution $\mu$ which satisfies Assumption $A(2+\delta)$ for some $\delta>0$. Then, for any $\nu\in (0,\delta/(4+2\delta))$, a.s. and in $L^1(\PP)$,
  $$\lim_{n\to \infty}n^{\nu}\max\{|\Delta_n(\|\bm{x}\|^2)|,|\Delta_n(\|\bm{x}\|^{2+\frac{\delta}{2}})|,\|\Delta_n(\bm{x}\bm{x}^T)\| \}=0.$$
\end{corollary}
\begin{proof}
First consider $h(\bm{x}):=\|\bm{x}\|^{2+\frac{\delta}{2}}$. Then by our assumption, $\EE |h(\bm{X}_1)|^s<\infty$ where 
$$s=\frac{2+\delta}{2+\frac{\delta}{2}}=\frac{4+2\delta}{4+\delta}.$$
  Apply Proposition \ref{prop1dslln} to conclude that $n^{\nu} |\Delta_n(h)| \to 0$ a.s. and in $L^1(\PP)$. The proof for $|\Delta_n(\|\bm{x}\|^2)|$ and $\|\Delta_n(\bm{x}\bm{x}^T)\|$ is similar where one can consider $h(\bm{x}):=\|\bm{x}\|^{2}$ and $h(\bm{x}):=x(i)x(j)$ for $i,j\in \{1,2,\dots,d\}$.
\end{proof}

On the other hand, given that $|\Delta_n(\|\bm{x}\|^2)|$, $|\Delta_n(\|\bm{x}\|^{2+\frac{\delta}{2}})|$ and $\left\|\Delta_n(\bm{x}\bm{x}^T)\right\|$ are all small, we can apply the following lemma to control the distribution of the $n+1$-st step. It will be used in Section \ref{secdiff} to control unusually large increments which make Taylor expansions break down.

Recall from (\ref{Evarsetdef}) that for $\varepsilon\in (0,1)$, we write
$$
      E_{\varepsilon}(\bm{x}):=\left\{\bm{y} \in \RR^d:\|\bm{y}\| \leq\|\bm{x}\|^{1-\varepsilon}\right\}.
$$

\begin{lemma}
  \label{condEYq}
  Let $\bm{S}=(\bm{S}_n)_{n\in \NN}$ be an SRRW in $\RR^d$ with parameter $\alpha \leq 1/2$ and step distribution $\mu$ which satisfies Assumption $A(2+\delta)$ for some $\delta>0$. Let $n$ be a positive integer and $\varepsilon\in (0,1)$ be a constant. Suppose that $\bm{S}_n\neq 0$ and 
  \begin{equation}
    \label{bnnmuxdelta}
    n^{\frac{\nu}{2}}\max\{|\Delta_n(\|\bm{x}\|^2)|,|\Delta_n(\|\bm{x}\|^{2+\frac{\delta}{2}})|,\|\Delta_n(\bm{x}\bm{x}^T)\|\}\leq 1,
  \end{equation}
  where $\nu>0$ is a constant. Then for any $q\in [0, 2+\frac{\delta}{2}]$, there exists a constant $C$ independent of $\varepsilon, q$ and $n$, such that
  \begin{equation}
    \label{Xqest}
    \EE(\|\bm{X}_{n+1}\|^q\mathds{1}_{\{X_{n+1}\notin E_{\varepsilon}(\bm{S}_n)\}}\mid \FF_n) \leq \frac{C}{\|\bm{S}_n\|^{(1-\varepsilon)(2+\frac{\delta}{2}-q)}} .
  \end{equation}
  In particular, as $\|\bm{S}_n\|\to \infty$,
  \begin{equation}
    \label{approfirmo}
      \EE(\bm{X}_{n+1}\mathds{1}_{\{\bm{X}_{n+1} \in E_{\varepsilon}(\bm{S}_n)\}}\mid \FF_n)= \frac{\alpha\bm{S}_n}{n}+O\left(\frac{1}{\|\bm{S}_n\|^{(1-\varepsilon)(1+\frac{\delta}{2})}}\right),
  \end{equation}
  and as $\|\bm{S}_n\|\to \infty, n \to \infty$,
  \begin{equation}
    \label{approsecmo}
  \EE(\|\bm{X}_{n+1}\|^2\mathds{1}_{\{\bm{X}_{n+1}\in E_{\varepsilon}(\bm{S}_n)\}}\mid \FF_n)= \EE \|\bm{X}_1\|^2 +\frac{1}{n^{\frac{\nu}{2}}}+ O\left(\frac{1}{\|\bm{S}_n\|^{\frac{\delta}{2}(1-\varepsilon)}}\right).
\end{equation}
 If we further assume that $\EE \bm{X}_1\bm{X}_1^T$ is the identity matrix, then as $\|\bm{S}_n\|\to \infty, n \to \infty$,
 \begin{equation}
  \label{approSnXnmo}
  \EE\left(\frac{\|\bm{S}_n\cdot \bm{X}_{n+1}\|^2}{\|\bm{S}_n\|^2}\mathds{1}_{\{\bm{X}_{n+1}\in E_{\varepsilon}(\bm{S}_n)\}}\mid \FF_n\right)= 1  +\frac{1}{n^{\frac{\nu}{2}}}+O\left(\frac{1}{\|\bm{S}_n\|^{\frac{\delta}{2}(1-\varepsilon)}}\right).
\end{equation}
\end{lemma}
\begin{proof}
  By Lemma \ref{condlemXn1} and (\ref{bnnmuxdelta}), 
$$
  \EE(\|\bm{X}_{n+1}\|^{2+\frac{\delta}{2}}\mid \FF_n) =\EE\|\bm{X}_1\|^{2+\frac{\delta}{2}}+ \alpha \Delta_n(\|\bm{x}\|^{2+\frac{\delta}{2}})\leq \EE\|\bm{X}_1\|^{2+\frac{\delta}{2}}+1.
$$
For any $q\in [0, 2+\frac{\delta}{2}]$, 
  \begin{equation}
      \label{XqoutvarSbd}
      \|\bm{X}_{n+1}\|^q\mathds{1}_{\{\bm{X}_{n+1} \notin E_{\varepsilon}(\bm{S}_n)\}}=\frac{\|\bm{X}_{n+1}\|^{2+\frac{\delta}{2}}}{\|\bm{X}_{n+1}\|^{2+\frac{\delta}{2}-q}}\mathds{1}_{\{\bm{X}_{n+1} \notin E_{\varepsilon}(\bm{S}_n)\}}\leq \frac{\|\bm{X}_{n+1}\|^{2+\frac{\delta}{2}}}{\|\bm{S}_n\|^{(1-\varepsilon)(2+\frac{\delta}{2}-q)}}.
  \end{equation}
One obtains (\ref{Xqest}) by taking the conditional expectation on both sides of (\ref{XqoutvarSbd}).

Recall from Lemma \ref{condlemXn1} that 
$$\EE(\bm{X}_{n+1}\mid \FF_n)=\frac{\alpha\bm{S}_n}{n}, \quad \EE(\|\bm{X}_{n+1}\|^2\mid \FF_n)= \EE\|\bm{X}_1\|^{2}+ \alpha \Delta_n(\|\bm{x}\|^2).$$ 
Then (\ref{approfirmo}) and (\ref{approsecmo}) follow from (\ref{bnnmuxdelta}) and (\ref{Xqest}) with $q=1,2$. If $\EE \bm{X}_1\bm{X}_1^T$ is the identity matrix, then Lemma \ref{condlemXn1} (with $h(\bm{x})=x(i)x(j)$, $i,j\in \{1,2,\dots,d\}$) shows that 
$$
\EE\left(\frac{\|\bm{S}_n\cdot \bm{X}_{n+1}\|^2}{\|\bm{S}_n\|^2}\mid \FF_n\right)=\bm{S}_n^T \EE\left(\frac{ \bm{X}_{n+1}\bm{X}_{n+1}^T}{\|\bm{S}_n\|^2}\mid \FF_n\right)\bm{S}_n=1+\frac{\alpha\bm{S}_n^T \Delta_n(\bm{x}\bm{x}^T) \bm{S}_n}{\|\bm{S}_n\|^2}.
$$
Observe that $\|\bm{S}_n^T \Delta_n(\bm{x}\bm{x}^T) \bm{S}_n\|\leq \|\bm{S}_n\|^2\|\Delta_n(\bm{x}\bm{x}^T)\|$ and $\|\bm{S}_n\cdot \bm{X}_{n+1}\|^2\leq \|\bm{S}_n\|^2\|\bm{X}_{n+1}\|^2$. The rest of the proof for (\ref{approSnXnmo}) is then similar to that of (\ref{approsecmo}).
\end{proof}

\subsection{Proof for dimensions not equal to 2}
\label{secdiff}

\subsubsection{SRRW in dimension 1}
\label{proof1dsec}

This section is devoted to the proof of the following proposition.
\begin{proposition}
  \label{1DSRRWprop}
  Let $S$ be an SRRW in $\RR$ with parameter $\alpha\leq 1/2$ and step distribution $\mu$ such that $\EE X_1^2=1$. Then $S$ is recurrent. 
\end{proposition}

The following two auxiliary lemmas will be used in the proof of Proposition \ref{1DSRRWprop}.

\begin{lemma}
  \label{Lxdelta1dest}
  Define a function $L$ on $\RR$ by
\begin{equation}
  \label{L1Ddef}
  L(x) := \sqrt{|x|}, \quad x\in \RR.
\end{equation}
 There exist two positive constants $\varepsilon<1$ and $C(\varepsilon)>1$ such that for any $x\neq 0$ and $y \in \RR$,
  \begin{equation}
    \label{Lest}
    L(x+y)-L(x) \leq \sqrt{|x|} \left(\frac{y}{2x} -\frac{y^2}{10x^2}+\frac{C(\varepsilon)y^2\mathds{1}_{\{|y| > \varepsilon|x|\}}}{x^2}\right),
  \end{equation} 
  where, by a slight abuse of notation, we used $\mathds{1}_{\{|y| > \varepsilon|x|\}}$ for the indicator function of the set $\{(x,y)\in \RR^2: |y| > \varepsilon|x|\}$.
\end{lemma}
\begin{proof}
  By Taylor's expansion: As $x\to 0$,
  $$
  \sqrt{1+x}=1+\frac{x}{2}-\frac{x^2}{8}+O(|x|^3), 
  $$
there exists a constant $\varepsilon\in (0,1)$ such that for all $t\in (-\varepsilon,\varepsilon)$, one has 
$$
\sqrt{1+t}-1 \leq \frac{t}{2}-\frac{t^2}{10}.
$$
Note that if $|y| \leq \varepsilon |x| \neq 0$, then $x+y$ and $x$ have the same sign, and in particular, $|x+y|/|x|=(x+y)/x$. Therefore,
\begin{equation}
  \label{LxdeltaEvare}
  \begin{aligned}
   &\quad \ (L(x+y)-L(x))\mathds{1}_{\{|y| \leq \varepsilon|x|\}}  \\
   &=\sqrt{|x|}\left(\sqrt{1+\frac{y}{x}}-1\right)\mathds{1}_{\{|y| \leq \varepsilon|x|\}} \\
    &\leq\sqrt{|x|}\left(  \frac{y}{2x} -\frac{y^2}{10x^2}\right) \mathds{1}_{\{|y| \leq \varepsilon|x|\}} \\
   &\leq \sqrt{|x|}\left(  \frac{y}{2x} -\frac{y^2}{10 x^2}\right)+\sqrt{|x|}\left(  \frac{y^2}{2\varepsilon x^2} +\frac{y^2}{10 x^2}\right)\mathds{1}_{\{|y| >\varepsilon|x|\}}. 
\end{aligned}
\end{equation}
For $|y| > \varepsilon |x|$, one has,
$$
\sqrt{1+\frac{|y|}{|x|}}-1\leq \sqrt{\frac{|y|}{|x|}} \leq \frac{y^2}{\varepsilon^{3/2} x^2}.
$$
Thus,
\begin{equation}
  \label{Lxdeltagener}
  (L(x+y)-L(x))\mathds{1}_{\{|y| >\varepsilon|x|\}}  \leq \sqrt{|x|} \frac{y^2}{\varepsilon^{3/2} x^2} \mathds{1}_{\{|y| > \varepsilon|x|\}}.
\end{equation}
Now (\ref{Lest}) follows from (\ref{LxdeltaEvare}) and (\ref{Lxdeltagener}).
\end{proof}

\begin{lemma}
  \label{Svarsqrtn}
  Let $S$ be an SRRW in $\RR$ with parameter $\alpha\leq 1/2$ and step distribution $\mu$ which satisfies Assumption $A(2)$. Then, for any $\varepsilon>0$, almost surely, $|S_n|\leq \varepsilon \sqrt{n}$ infinitely often.
\end{lemma}
\begin{proof}
 Recall $(\xi_n)_{n\geq 2}$ in Definition \ref{defSRRW}. Given that $\xi_n=0$, $X_n$ is sample from $\mu$ independently of $\FF_{n-1}$. For any $\varepsilon>0$, we have
  $$
      \sum_{n=2}^{\infty}\PP(|X_{n}|\mathds{1}_{\{\xi_{n}=0\}}>\varepsilon \sqrt{n}) \leq \sum_{n=2}^{\infty} \PP(\frac{X_1^2}{\varepsilon^2}\geq  n)\leq \sum_{n=2}^{\infty}n\PP(n\leq \frac{X_1^2}{\varepsilon^2}< n+1)\leq \frac{\EE X_1^2}{\varepsilon^2}<\infty.
  $$
 The Borel-Cantelli lemma implies that a.s. there exists a (random) $k<\infty$ such that if $n\geq k$ and $\xi_{n}=0$, then $|X_n|<\varepsilon \sqrt{n}$. We prove by induction that for any $n\geq 1$,
 \begin{equation}
     \label{XninducvarsqrtnK}
     |X_n|\leq \max\{\varepsilon\sqrt{n},K\}, \quad \text{where}\ K:=\max_{1\leq j\leq k}|X_j|.
 \end{equation}
By definition, (\ref{XninducvarsqrtnK}) holds for all $n\leq k$. Assume that (\ref{XninducvarsqrtnK}) holds for all $n\leq m$ ($m\geq k$). For $n=m+1$, if $\xi_{n}=0$, then $|X_n|<\varepsilon \sqrt{n}$ by the choice of $k$; if $\xi_{n}=1$, then $X_n=X_{j}$ for some $j\leq m$, and thus, by our assumption, 
$$|X_n|\leq \max\{\varepsilon\sqrt{j},K\}\leq \max\{\varepsilon\sqrt{n},K\},$$
 which completes the proof of (\ref{XninducvarsqrtnK}).

  On the other hand, (\ref{lildiff}) (for $\alpha<1/2$) and (\ref{lilcrit}) (for $\alpha=1/2$) yield that 
  \begin{equation}
    \label{liminfsupinfity1d}
    \liminf_{n\to \infty} S_n=-\infty, \quad \limsup_{n\to \infty} S_n=\infty \quad \text{a.s.} 
  \end{equation} 
  This proves the desired results since otherwise one can find an infinite sequence $(n_j)_{j\geq 1}$ such that $S_{n_j}>\varepsilon \sqrt{n_j}$ but $S_{n_j+1}<-\varepsilon \sqrt{n_j}$, and, in particular, 
  $$2\varepsilon \sqrt{n_j}< |X_{n_j+1}|\leq \max\{\varepsilon\sqrt{n_j+1},K\},$$ 
  which is impossible for large $n_j$.
\end{proof}

Under the setting of Proposition \ref{1DSRRWprop}, for a function $h$ on $\RR$, recall that 
$$
\Delta_n(h)=\frac{S_n^{(h)}}{n}-\EE h(X_1), \quad n\geq 1.
$$
If $\EE |h(X_1)|<\infty$, then the law of large numbers (\ref{llnsh}) implies that, almost surely, 
\begin{equation}
    \label{edltanhslln1d}
    \lim_{n\to \infty}\Delta_n(h)=0.
\end{equation}

\begin{proof}[Proof of Proposition \ref{1DSRRWprop}.] 
  Let $L$, $\varepsilon$ and $C(\varepsilon)$ be as in Lemma \ref{Lxdelta1dest}. We choose $K>0$ such that $\EE X_1^2 \mathds{1}_{\{|X_1|\geq K\}}<1/(60C(\varepsilon))$.  By (\ref{edltanhslln1d}) with $h(x)=x^2$ and $h(x)=x^2\mathds{1}_{\{|x|\geq K\}}$, a.s.,
  \begin{equation}
    \label{Inx2K0}
    \lim_{n\to \infty}\Delta_n(x^2)=0,\quad \lim_{n\to \infty}\Delta_n(x^2\mathds{1}_{\{|x|\geq K\}})=0.
  \end{equation}
   We define inductively a sequence of stopping times $(n_k)_{k\geq 1}$. Set $n_1:=1$. For $k>1$, set
  $$
  n_k:=\inf\left\{n>n_{k-1}: |S_{n}|\leq \sqrt{\frac{n}{20}},\  \max\{|\Delta_n(x^{2})|,|\Delta_n(x^2\mathds{1}_{\{|x|\geq K\}})|\}  < \frac{1}{60C(\varepsilon)}\right\}
  $$
  with the convention that $\inf \emptyset=\infty$. By Lemma \ref{Svarsqrtn} and (\ref{Inx2K0}), almost surely, $n_k$ is finite for all $k\geq 1$.  Applying Lemma \ref{Lxdelta1dest} with $x=S_n$ and $y=X_{n+1}$, we see that if $S_n \notin B(0,r)$ where $r>K/\varepsilon$ is a constant, then
  \begin{equation}
    \label{LSnincreine}
    L(S_{n+1})-L(S_n)\leq \sqrt{|S_n|}\left(\frac{X_{n+1}}{2S_n} -\frac{ X_{n+1}^2}{10 S_n^2}+\frac{C(\varepsilon) X_{n+1}^2\mathds{1}_{\{|X_{n+1}| \geq  K\}}}{S_n^{2}}\right).
  \end{equation}
  For each $k\geq 1$, we define stopping times $$\tau_k:=\inf\{n\geq n_k: S_n \in B(0,r) \},\quad \theta_k:=\inf\left\{n \geq n_k: |S_n| \geq \sqrt{\frac{n}{10}} \right\},$$ 
  and 
  $$T_k:=\inf\left\{n \geq n_k:  \max\{|\Delta_n(x^{2})|,|\Delta_n(x^{2}\mathds{1}_{\{|x|\geq K\}})|\}\geq \frac{1}{30C(\varepsilon)}\right\}.
  $$
 By Lemma \ref{condlemXn1}, $\EE(X_{n+1}\mid \FF_n)=\alpha S_n/n$, $\EE(X_{n+1}^2\mid \FF_n)=1+\alpha \Delta_n(x^2)$, and 
 $$  \EE( X_{n+1}^2\mathds{1}_{\{|X_{n+1}|>K\}}\mid \FF_n)=\EE X_1^2\mathds{1}_{\{|X_1|>K\}} +\alpha \Delta_n(x^2\mathds{1}_{\{|x|\geq K\}}).$$
 If $n_k\leq n < T_k$ and $S_n \notin B(0,r)$, then by (\ref{LSnincreine}) and that $\alpha\leq 1/2$ and $C(\varepsilon)>1$, one has
  \begin{equation}
   \label{LsnlogcondE}
   \begin{aligned}
  &\quad\ \mathbb{E}\left(L\left(S_{n+1}\right)-L\left(S_n\right) \mid \FF_n\right) \\
  &\leq \sqrt{|S_n|} \left( \frac{\alpha}{2n} -\frac{1}{10 S_n^2}-\frac{\alpha \Delta_n(x^2)}{10 S_n^2}+\frac{C(\varepsilon)\EE X_1^2\mathds{1}_{\{|X_1|>K\}}}{S_n^2}+\frac{\alpha C(\varepsilon) \Delta_n(x^2\mathds{1}_{\{|x|\geq K\}})}{S_n^2} \right)\\
  &\leq \sqrt{|S_n|} \left(  \frac{\alpha}{2n} -\frac{1}{10 S_n^2} +\frac{ 1}{600 C(\varepsilon) S_n^2}+\frac{1}{60 S_n^2}+\frac{1}{60 S_n^2} \right) \\
  &\leq \sqrt{|S_n|} \left(  \frac{\alpha}{2n} -\frac{1}{20 S_n^2} \right).
   \end{aligned}
  \end{equation}
   Then, by (\ref{LsnlogcondE}), $(L(S_{(n_k+j)\wedge \tau_k \wedge T_k \wedge \theta_k}))_{j\in \NN}$ is a non-negative supermartingale, and thus converges a.s.. By (\ref{lildiff}) and (\ref{lilcrit}), $\theta_k<\infty$ a.s., and hence
   $$\lim_{j\to \infty}L(S_{(n_k+j)\wedge \tau_k \wedge T_k  \wedge \theta_k}) = L(S_{\theta_k}) \quad \text{a.s. on}\ \{\tau_k=\infty\} \cap \{T_k=\infty\}.$$ 
  By the optional stopping theorem for non-negative supermartingales, see e.g. \cite[Theorem 16, Chapter V]{MR0745449}, for all $k> 1$,
   \begin{equation}
       \label{1dopsample}
         \begin{aligned}
      \left(\frac{n_k}{20}\right)^{\frac{1}{4}} &\geq \EE (L(S_{n_k  })\mid \FF_{n_k}) \geq \EE (L(S_{\theta_k})\mathds{1}_{\{\tau_k=\infty, T_k=\infty\}}\mid \FF_{n_k}) \\
      &\geq \EE \left( \left(\frac{\theta_k}{10}\right)^{\frac{1}{4}}\mathds{1}_{\{\tau_k=\infty, T_k=\infty\}}\mid \FF_{n_k}\right) \\
      &\geq \left(\frac{n_k}{10}\right)^{\frac{1}{4}} \PP(\{\tau_k=\infty\} \cap \{T_k=\infty\}\mid \FF_{n_k}),
   \end{aligned}
   \end{equation}
   in particular, $\PP(\{\tau_k=\infty\} \cap \{T_k=\infty\}\mid \FF_{n_k}) \leq (1/2)^{1/4}$. Indeed, we have shown that for any $m\geq k>1$,
  $$
  \PP(\{\tau_k=\infty\} \cap \{T_k=\infty\}\mid \FF_{n_m}) \leq \PP(\{\tau_m=\infty\} \cap \{T_m=\infty\}\mid \FF_{n_m}) \leq (\frac{1}{2})^{\frac{1}{4}}<1.
  $$
  Thus, for any $k>1$, $\PP(\{\tau_k=\infty\} \cap \{T_k=\infty\})=0$ by Lévy's 0-1 law, and hence,
  $$
  \begin{aligned}
    \PP(\tau_k<\infty )&= \PP(\bigcup_{m\geq k} \{\tau_m<\infty\})\\
    &\geq \PP(\bigcup_{m\geq k}\left(\{\tau_m<\infty\}\cap \{T_m=\infty\}  \right) ) \\
    &=\PP(\bigcup_{m\geq k} \{T_m=\infty\})=1,
  \end{aligned}
  $$
  where we used (\ref{Inx2K0}) in the last equality. This shows that almost surely for each $k>1$, $S$ will visit $B(0,r)$ after time $n_k$, whence $S$ is recurrent.
 \end{proof}

\subsubsection{SRRW in dimensions 3 and higher}
\label{proofhighd}

This section is devoted to the proof of the following proposition.

\begin{proposition}
  \label{highDSRRWprop}
  Let $\bm{S}$ be an SRRW in $\RR^d$ $(d\geq 3)$ with parameter $\alpha\leq 1/2$ and step distribution $\mu$. Assume that $\mu$ satisfies Assumption $A(2+\delta)$ for some $\delta>0$ and $\EE \bm{X}_1\bm{X}_1^T$ is the identity matrix, then 
  \begin{equation}
    \label{rateestescap}
    \lim_{n\to \infty}\frac{\|\bm{S}_n\| (\log n)^{1+\max\{1,\frac{4}{\delta}\}}}{\sqrt{n}}=\infty, \quad \text{a.s.}
  \end{equation}
\end{proposition}

Throughout Section \ref{proofhighd}, we let $\delta>0$ be as in Proposition \ref{highDSRRWprop}. Define a function $h$ on $\RR^d$ by
\begin{equation}
  \label{h3Ddef}
  h(\bm{x}) := \left \{ \begin{aligned} &\frac{1}{\|\bm{x}\|^{\frac{\delta}{4}}}  && \text{if } \|\bm{x}\|\geq 1, \\ & 1 && \text{if } \|\bm{x}\|< 1.  \end{aligned} \right.
\end{equation}

\begin{lemma}
  \label{hxdelta2dest}
 Assume that $\bm{x},\bm{y}\in \RR^d$ with $\bm{x} \neq \bm{0}$ and $\varepsilon \in (0,1)$. There exist positive constants $r$ and $C$ such that for all $\|\bm{x}\|\geq r$ and $\bm{y} \in E_{\varepsilon}(\bm{x})$,
  \begin{equation}
    \label{hxdeltaEvare}
   h(\bm{x}+\bm{y})-h(\bm{x}) \leq \frac{-\delta}{4\|\bm{x}\|^{2+\frac{\delta}{2}}}\left(\bm{x} \cdot \bm{y} +\frac{\|\bm{y}\|^2}{2}-(1+\frac{\delta}{8})\frac{(\bm{x} \cdot \bm{y})^2}{\|\bm{x}\|^2}-\frac{C\|\bm{y}\|^2}{\|\bm{x}\|^{\varepsilon}}\right).
  \end{equation}
\end{lemma}
\begin{proof}
  Using Taylor expansions, one can show that for $\bm{x} \in \RR^d$ with $\bm{x} \neq \bm{0}$ and $\bm{y} \in E_{\varepsilon}(\bm{x})$, 
$$
   \begin{aligned}
    &\quad\ \frac{1}{ \|\bm{x}+\bm{y}\|^{\frac{\delta}{4}}}-\frac{1}{\|\bm{x}\|^{\frac{\delta}{4}}}=\frac{1}{\|\bm{x}\|^{\frac{\delta}{4}}}\left(\left(1+\frac{2 \bm{x} \cdot \bm{y} +\|\bm{y}\|^2}{\|\bm{x}\|^2}\right)^{-\frac{\delta}{8}}-1\right) \\
    &=\frac{-\delta}{4\|\bm{x}\|^{\frac{\delta}{2}}}\left(\frac{ \bm{x} \cdot \bm{y}}{\|\bm{x}\|^2}+\frac{\|\bm{y}\|^2}{2\|\bm{x}\|^2}-\frac{1+\frac{\delta}{8}}{4}\frac{(2\bm{x} \cdot \bm{y}+ \|\bm{y}\|^2)^2}{\|\bm{x}\|^4}+O\left(\frac{\|\bm{y}\|^3}{\|\bm{x}\|^{3}}\right)\right) \\ 
  &=\frac{-\delta}{4\|\bm{x}\|^{2+\frac{\delta}{2}}}\left(\bm{x} \cdot \bm{y} +\frac{\|\bm{y}\|^2}{2}-(1+\frac{\delta}{8})\frac{(\bm{x} \cdot \bm{y})^2}{\|\bm{x}\|^2}+O(\frac{\|\bm{y}\|^2}{\|\bm{x}\|^{\varepsilon}})\right).
\end{aligned}
$$
\end{proof}

We first prove the transience, which will be used in the proof of Proposition \ref{highDSRRWprop}, more specifically, Inequality (\ref{nSnlogSn5}).

\begin{proposition}
  \label{hightranprop}
 Under the setting of Proposition \ref{highDSRRWprop}, the SRRW $\bm{S}$ is transient.
\end{proposition}
\begin{proof}
Let $\nu \in (0,\delta/(4+2\delta))$ and $\varepsilon \in (0,\delta/(8+2\delta))$ be two constants. Without loss of generality, we assume that $\delta<4$. We define inductively a sequence of stopping times $(n_k)_{k\geq 1}$. Set $n_1=1$. For $k>1$, let
  $$
  \begin{aligned}
    n_k:=\inf&\left\{n>n_{k-1}: \|\bm{S}_{n} \|\geq \sqrt{n},\ |\Delta_n(\|\bm{x}\|^{2+\delta})|<1,\right. \\ 
    &\left.\max\{|\Delta_n(\|\bm{x}\|^{2})|,|\Delta_n(\|\bm{x}\|^{2+\frac{\delta}{2}})|,\|\Delta_n(\bm{x}\bm{x}^T)\|\}< \frac{1}{n^{\nu}}\right\},
  \end{aligned}
  $$
  with the convention that $\inf \emptyset=\infty$. By Corollary \ref{limx2Delta} and the strong law of large numbers (\ref{llnSRRW}) and the law of the iterated logarithm (\ref{lilSRRW}), almost surely, $n_k<\infty$ for all $k\geq 1$. For $k\geq 1$, define 
  \begin{equation}
    \label{defTkhighD}
      T_k:=\inf\left\{n \geq n_k:  \max\{|\Delta_n(\|\bm{x}\|^{2})|,|\Delta_n(\|\bm{x}\|^{2+\frac{\delta}{2}})|,\|\Delta_n(\bm{x}\bm{x}^T)\|\}\geq \frac{1}{n^{\frac{\nu}{2}}} \right\}.
  \end{equation}
  Let $h$ be as in (\ref{h3Ddef}). Then, by Lemma \ref{hxdelta2dest} with $\bm{x}=\bm{S}_n,\bm{y}=\bm{X}_{n+1}$ and Equations (\ref{approfirmo}), (\ref{approsecmo}) (\ref{approSnXnmo}) in Lemma \ref{condEYq}, there exists positive constants $C_1$ and $r$ such that if $n_k\leq n<T_k$ and $\|\bm{S}_n\|\geq r$,
  \begin{equation}
    \label{hhighDest}
      \begin{aligned}
   &\quad\ \mathbb{E}\left((h\left(\bm{S}_{n+1}\right)-h\left(\bm{S}_n\right))\mathds{1}_{\{\bm{X}_{n+1} \in E_{\varepsilon}(\bm{S}_n)\}} \mid \FF_n\right) \\
   &\leq \frac{-\delta}{4\|\bm{S}_n\|^{2+\frac{\delta}{4}}}\left(\frac{\alpha\|\bm{S}_n\|^2}{n}+\frac{d}{2}-(1+\frac{\delta}{8})-\frac{C_1}{\|\bm{S}_n\|^{\varepsilon}}-\frac{C_1}{\|\bm{S}_n\|^{\frac{\delta}{2}-\varepsilon(1+\frac{\delta}{2})}}-\frac{C_1}{n^{\frac{\nu}{2}}}\right)  \\
    &\leq \frac{-\delta}{4\|\bm{S}_n\|^{2+\frac{\delta}{4}}}\left(\frac{3}{2}-(1+\frac{\delta}{8})-\frac{C_1}{\|\bm{S}_n\|^{\varepsilon}}-\frac{C_1}{\|\bm{S}_n\|^{\frac{\delta}{2}-\varepsilon(1+\frac{\delta}{2})}}-\frac{C_1}{n^{\frac{\nu}{2}}}\right),
  \end{aligned}
  \end{equation}
  where we used the assumption that $\EE \bm{X}_1\bm{X}_1^T$ is the identity matrix and that $d\geq 3$. Since we have assumed that $\delta<4$, by possibly choosing larger $r$ and $k$, we see that the expression in the parenthesis of the last line in (\ref{hhighDest}) is lower bounded by $(4-\delta)/16$, and thus,
\begin{equation}
  \label{hcondEEvar}
  \mathbb{E}\left((h\left(\bm{S}_{n+1}\right)-h\left(\bm{S}_n\right))\mathds{1}_{\{\bm{X}_{n+1} \in E_{\varepsilon}(\bm{S}_n)\}} \mid \FF_n\right)\leq \frac{-\delta(4-\delta)}{64\|\bm{S}_n\|^{2+\frac{\delta}{4}}}<0.
\end{equation}
Note that $h(\bm{x})\in (0,1]$ for all $\bm{x}\in \RR^d$. If $n_k\leq n<T_k$, by Lemma \ref{condEYq} with $q=0$, for some constant $C_2$,
\begin{equation}
  \label{hcondEEnotvar}
  \begin{aligned}
     &\quad\ \mathbb{E}\left((h\left(\bm{S}_{n+1}\right)-h\left(\bm{S}_n\right))\mathds{1}_{\{\bm{X}_{n+1} \notin E_{\varepsilon}(\bm{S}_n)\}} \mid \FF_n\right)\\
     &\leq \PP(\bm{X}_{n+1} \notin E_{\varepsilon}(\bm{S}_n)\mid \FF_n) \leq \frac{C_2}{\|\bm{S}_n\|^{(1-\varepsilon)(2+\frac{\delta}{2})}}.
  \end{aligned}
\end{equation}
By the choice of $\varepsilon$, we see that $2+\frac{\delta}{4}<(1-\varepsilon)(2+\frac{\delta}{2})$. Then, (\ref{hcondEEvar}) and (\ref{hcondEEnotvar}) imply that if $n_k\leq n<T_k$ and $\|\bm{S}_n\|\geq r$ with possibly a larger $r$,
\begin{equation}
  \label{hSsuper}
  \mathbb{E}\left(h\left(\bm{S}_{n+1}\right)-h\left(\bm{S}_n\right) \mid \FF_n\right)<0.
\end{equation}
Now define $\tau_k:=\inf\{n\geq n_k: \bm{S}_n \in B(0,r) \}$. By (\ref{hSsuper}), for all large $k$, $(h(\bm{S}_{(n_k+j)\wedge  T_k \wedge \tau_k}))_{j\in \NN}$ is a non-negative supermartingale, and thus converges a.s.. Observe that 
$$\lim_{n\to \infty}h(\bm{S}_{(n_k+j)\wedge  T_k  \wedge \tau_k}) = h(\bm{S}_{\tau_k}) \quad \text{a.s. on}\ \{\tau_k<\infty\} \cap \{T_k=\infty\}.$$ 
Then by the optional stopping theorem for non-negative supermartingales, 
$$
\begin{aligned}
    n_k^{-\frac{\delta}{8}} &\geq \EE (h(\bm{S}_{n_k})\mid \FF_{n_k}) \\
    &\geq \EE (h(\bm{S}_{\tau_k})\mathds{1}_{\{\tau_k<\infty, T_k=\infty\}}\mid \FF_{n_k}) \\
   &\geq \frac{\PP(\{\tau_k<\infty\} \cap \{T_k=\infty\}\mid \FF_{n_k})}{r^{\frac{\delta}{4}}}. 
\end{aligned}
$$
Therefore, for all large $k$, $\PP(\{\tau_k<\infty\} \cap \{T_k=\infty\}\mid \FF_{n_k}) \leq 1/2$. By (\ref{pnkIg1}) in Proposition \ref{prop1dslln}, $\PP(T_k=\infty\mid \FF_{n_k}) \to 1$ as $k\to \infty$. Thus, for all large $k$, 
$$
\PP(\{\tau_k=\infty\} \cap \{T_k=\infty\}\mid \FF_{n_k}) \geq \PP(T_k=\infty\mid \FF_{n_k}) -\frac{1}{2}\geq \frac{1}{3}.
$$ 
By (\ref{lilSRRW}), $\lim_{n\to \infty}h(\bm{S}_n)=0$ (i.e. $\lim_{n\to \infty}\|\bm{S}_n\| = \infty$) a.s. on $\{\tau_k=\infty\} \cap \{T_k=\infty\}$. In particular,
$$
\PP(\lim_{n\to \infty}\|\bm{S}_n\| = \infty\mid \FF_{n_k})\geq \frac{1}{3},
$$
which implies that $\lim_{n\to \infty}\|\bm{S}_n\| = \infty$ a.s. by Lévy's 0-1 law.  
\end{proof}

 For an SRRW $\bm{S}$ in $\RR^d$ $(d\geq 1)$ and a positive real number $R$, we denote by $\zeta_R$ the exit time of $\bm{S}$ from the ball $B(0,R)$, i.e.
\begin{equation}
  \label{exittimedef}
  \zeta_R:=\inf\{n\in \NN: \|\bm{S}_n\|\geq R\}.
\end{equation}

\begin{lemma}
\label{expzetaRfinite}
    Let $\bm{S}$ be an SRRW in $\RR^d$ with parameter $\alpha$ and step distribution $\mu$. For any $R>0$, $\EE \zeta_R<\infty$.
\end{lemma}
\begin{proof}
Let $\bm{X}^{\prime}_1,\bm{X}^{\prime}_2,\dots$ be i.i.d. $\mu$-distributed random vectors. Sine we have assumed that $\mu(\bm{0})<1$, there exists a positive integer $m$ such that 
$$
\PP(\|\sum_{i=1}^m \bm{X}^{\prime}_i\|>2R)>0.
$$
Let $(\xi_n)_{n\geq 1}$ be as in Definition \ref{defSRRW} with the convention that $\xi_1:=0$. For any $n\in \NN$, conditionally on that $\xi_{n+1}=0,\xi_{n+2}=0,\dots,\xi_{n+m}=0$, one has
$$\bm{S}_{n+m}-\bm{S}_{n} \stackrel{\mathcal{L}}{=} \sum_{i=1}^m \bm{X}^{\prime}_i.$$
Thus,
$$
\PP(\zeta_R\leq n+m  \mid \zeta_R >n) \geq \PP(\|\bm{S}_{n+m}-\bm{S}_{n}\|>2R  \mid \zeta_R >n) \geq (1-\alpha)^m \PP(\|\sum_{i=1}^m \bm{X}^{\prime}_i\|>2R).
$$
We write $c_1:=(1-\alpha)^m \PP(\|\sum_{i=1}^m \bm{X}^{\prime}_i\|>2R)>0$. Then, one can prove by induction that
$$
\PP(\zeta_R>km)\leq (1-c_1)^k, \quad k=1,2,\dots,
$$
 whence we have $\EE \zeta_R<\infty$.
\end{proof}

The proof of Proposition \ref{highDSRRWprop} also needs the following lemma, which provides some estimates for the exit times.
 
\begin{lemma}
  Let $\bm{S}$ be an SRRW in $\RR^d$ with parameter $\alpha$ and step distribution $\mu$ such that Assumption $A(2)$ holds. Then, there exists a constant $C(\alpha,\mu)$ depending on $\alpha$ and $\mu$ such that for all $R>0$,
  \begin{equation}
    \label{Ezetamupbd}
      \EE \zeta_R \leq C(\alpha,\mu)R^2.
  \end{equation}
  Moreover, for any positive integer $K$, there exists a constant $C(K,\mu)$ such that for all $R>0$,
  \begin{equation}
    \label{zetamnotsmall}
    \PP(\zeta_R \leq K) \leq \frac{C(K,\mu)}{R^2}.
  \end{equation}
\end{lemma} 
\begin{proof}
 Define a process $(M_n)_{n\geq 1}$ by $M_0:=0$ and $M_n:=\|\bm{S}_{n}\|^2 -n(1-\alpha)\EE \|\bm{X}_1\|^2$ for $n\geq 1$. By definition, 
 $$
\EE (\|\bm{S}_{n+1}\|^2 \mid \FF_n)= \|\bm{S}_n\|^2+\frac{2\alpha\|\bm{S}_n\|^2}{n}+\EE(\|\bm{X}_{n+1}\|^2\mid \FF_n)\geq \|\bm{S}_n\|^2+(1-\alpha)\EE \|\bm{X}_1\|^2,
 $$
 which shows that $(M_n)_{n\geq 1}$ is a submartingale. Thus, by the optional stopping theorem, 
 \begin{equation}
  \label{zetamopstop}
(1-\alpha)\EE \|\bm{X}_1\|^2\EE (\zeta_R\wedge n) \leq \EE \|\bm{S}_{\zeta_R \wedge n}\|^2 \leq 3 R^2+2\EE \|\bm{X}_{\zeta_R}\|^2
 \end{equation}
 where we used that 
 $$
 \begin{aligned}
   \EE \|\bm{S}_{\zeta_R \wedge n}\|^2 &=  \EE \|\bm{S}_n\|^2\mathds{1}_{\{n<\zeta_R\}}+\EE (\|\bm{S}_{\zeta_R-1}+\bm{X}_{\zeta_R}\|^2\mathds{1}_{\{\zeta_R\leq n\}}) \\
   &\leq  \EE \|\bm{S}_n\|^2 \mathds{1}_{\{n<\zeta_R\}} + 2 \EE (\|\bm{S}_{\zeta_R-1}\|^2\mathds{1}_{\{\zeta_R\leq n\}}) + 2 \EE \|\bm{X}_{\zeta_R}\|^2.
 \end{aligned}
 $$
We define a stopping time 
$$
T_R:=\inf\{n\geq 1: \|\bm{X}_n\|\geq 2R\}.
$$
Observe that $\zeta_R\leq T_R$. Then,
\begin{equation}
  \label{equEXtaum}
  \begin{aligned}
   \EE \|\bm{X}_{\zeta_R}\|^2&=\EE (\|\bm{X}_{\zeta_R}\|^2\mathds{1}_{\{\zeta_R<T_R\}})+\EE(\|\bm{X}_{T_R}\|^2\mathds{1}_{\{\zeta_R=T_R\}})\\
   &\leq 4R^2+\sum_{n=0}^{\infty}\EE(\|\bm{X}_{n+1}\|^2\mathds{1}_{\{\zeta_R=T_R=n+1\}}) .
  \end{aligned}
\end{equation}
Note that if $T_R>n$, then $\mu_n((B(0,2R))^c)=0$, and thus, $T_R=n+1$ only if $\bm{X}_{n+1}$ is sampled according to $\mu$ independently of $\FF_n$. Therefore,
$$
\begin{aligned}
  \EE(\|\bm{X}_{n+1}\|^2\mathds{1}_{\{\zeta_R=T_R=n+1\}}\mid \FF_n)&\leq \EE (\|\bm{X}_{n+1}\|^2\mathds{1}_{\{\|\bm{X}_{n+1}\|\geq 2R\}}\mid \FF_n)\mathds{1}_{\{\zeta_R>n,T_R>n\}} \\
  &\leq (1-\alpha) \mathds{1}_{\{\zeta_R>n\}}  \EE 
(\|\bm{X}_1\|^2 \mathds{1}_{\{\|\bm{X}_1\|\geq 2R\}}) .
\end{aligned}
$$
Then, by (\ref{equEXtaum}), 
$$
\begin{aligned}
  \EE \|\bm{X}_{\zeta_R}\|^2&\leq 4R^2+ (1-\alpha)\EE( 
\|\bm{X}_1\|^2 \mathds{1}_{\{\|\bm{X}_1\|\geq 2R\}})\sum_{n=0}^{\infty}\PP(\zeta_R>n)\\
  &=4R^2+(1-\alpha)\EE 
(\|\bm{X}_1\|^2 \mathds{1}_{\{\|\bm{X}_1\|\geq 2R\}})\EE \zeta_R
\end{aligned}
$$
For large $R>0$, we have $4\EE 
(\|\bm{X}_1\|^2 \mathds{1}_{\{\|\bm{X}_1\|\geq 2R\}})\leq \EE \|\bm{X}_1\|^2$. Thus, by Lemma \ref{expzetaRfinite} and (\ref{zetamopstop}), one has,
$$
(1-\alpha)\EE \|\bm{X}_1\|^2(2\EE \zeta_R\wedge n -\EE \zeta_R)\leq 22 R^2.
$$
We let $n\to \infty$ to obtain that for all large $R$,
$$
\EE \zeta_R\leq \frac{22R^2}{(1-\alpha)\EE \|\bm{X}_1\|^2},
$$
which completes the proof of (\ref{Ezetamupbd}). 

Recall from Lemma \ref{condlemXn1} that $\EE \|\bm{X}_n\|^2=\EE \|\bm{X}_1\|^2$ for all $n\geq 1$. Then, 
$$
\PP(\zeta_R\leq K) \leq \sum_{i=1}^K \PP(\|\bm{X}_i\|\geq \frac{R}{K}) \leq  \frac{K^2}{R^2}\sum_{i=1}^K \EE \|\bm{X}_i\|^2 =\frac{K^3 \EE \|\bm{X}_1\|^2}{R^2},
$$
which implies (\ref{zetamnotsmall})
\end{proof}

Now we are ready to prove Proposition \ref{highDSRRWprop}. We shall adapt the proof of \cite[Theorem 3.10.1]{MR3587911} to our setting.
\begin{proof}[Proof of Proposition \ref{highDSRRWprop}.]
Again, we first assume that $\delta<4$. Recall $\zeta_R$ defined in (\ref{exittimedef}). For $m_2>m_1 >0$, let $\lambda_{m_2, m_1}:=\inf \left\{n \geq \zeta_{m_2}: \|\bm{S}_n\| \leq m_1\right\}$ and 
$$
\tilde{T}_{m_2}:=\inf\{n \geq \zeta_{m_2}:  \ \max\{|\Delta_n(\|\bm{x}\|^{2})|,|\Delta_n(\|\bm{x}\|^{2+\frac{\delta}{2}})|,\|\Delta_n(\bm{x}\bm{x}^T)\|\}\geq \frac{1}{n^{\frac{\nu}{2}}}\}.
$$ 
Similarly as in (\ref{hSsuper}), given that $\zeta_{m_2}\geq K$ for some large $K$ and $m_1>r$ for some large $r$, one can deduce from (\ref{hhighDest}) and (\ref{hcondEEnotvar}) that  $(h(S_{(\zeta_{m_2}+j) \wedge \tilde{T}_{m_2}\wedge \lambda_{m_2, m_1}}))_{j\in \NN}$ is a non-negative supermartingale (here we used that $\delta<4$). Then, similarly as in the proof of (\ref{inequtkcapTk}), by the optional stopping theorem for non-negative supermartingales, for any $m_1>r$, we have
\begin{equation}
  \label{lammim2bd}
  \begin{aligned}
   \frac{\mathds{1}_{\{\zeta_{m_2}\geq K\}}}{m_2^{\frac{\delta}{4}}} &\geq h(\bm{S}_{\zeta_{m_2}}) \mathds{1}_{\{\zeta_{m_2}\geq K\}} \geq \EE (h(\bm{S}_{\lambda_{m_2, m_1}})\mathds{1}_{\{\lambda_{m_2, m_1}<\infty, \tilde{T}_{m_2}=\infty,\zeta_{m_2}\geq K\}}\mid \FF_{\zeta_{m_2}})  \\
   &\geq  \frac{1}{m_1^{\frac{\delta}{4}}}\mathbb{P}\left(\lambda_{m_2, m_1}<\infty, \tilde{T}_{m_2}=\infty, \zeta_{m_2}\geq K\mid \mathcal{F}_{\zeta_{m_2}}\right).
\end{aligned}
\end{equation}
For any $x\geq 0$, define
$$
\eta_x:=\sup \left\{n \geq 0: \|\bm{S}_n\| \leq x\right\}.
$$
Note that $\eta_x$ is not a stopping time. By (\ref{Ezetamupbd}), (\ref{zetamnotsmall}) and (\ref{lammim2bd}), there exists a large $r$ such that if $m_1>r$, then
\begin{equation}
  \label{3101etaT}
  \begin{aligned}
  &\quad \ \mathbb{P}\left(\eta_{m_1}>n, \tilde{T}_{m_2}=\infty\right) \\
  &\leq \PP(\zeta_{m_2}\leq K)+ \mathbb{P}(\eta_{m_1}>n, K\leq \zeta_{m_2}\leq n, \tilde{T}_{m_2}=\infty)+\mathbb{P}(\zeta_{m_2}> n, \tilde{T}_{m_2}=\infty) \\
  &\leq \PP(\zeta_{m_2}\leq K)+ \mathbb{P}(\lambda_{m_2, m_1}<\infty, \tilde{T}_{m_2}=\infty,\zeta_{m_2}\geq K)+\mathbb{P}(\zeta_{m_2}>n)\\ 
  &\leq \frac{C(K)}{m_2^2}+\left(\frac{m_1}{m_2}\right)^{\frac{\delta}{4}}+\frac{C(\alpha)m_2^2}{n},
\end{aligned}
\end{equation}
where we also used Markov's inequality in the last inequality. Fix $\varepsilon\in (0,1/9)$ and set
$$
m_1(k)=2^k,\quad m_2(k)=2^k(\log 2^k)^{\frac{4}{\delta}+\varepsilon}, \quad  n(k)=4^k(\log 2^k)^{\frac{8}{\delta}+1+3 \varepsilon}.
$$
Then, for all $k$ sufficiently large, by (\ref{3101etaT})
$$
\mathbb{P}(\eta_{2^k}>n(k),\tilde{T}_{m_2(k)}=\infty) \leq \frac{C(K)}{4^k(k \log 2)^{\frac{8}{\delta}+2\varepsilon}} +\frac{1}{(k \log 2)^{1+\frac{\delta\varepsilon}{4} }} +\frac{C(\alpha)}{(k \log 2)^{1+\varepsilon}}.
$$
In particular,
 $$\sum_{k \in \NN} \mathbb{P}\left(\eta_{2^k}>n(k),\tilde{T}_{m_2(k)}=\infty\right)<\infty.$$ 
 The Borel-Cantelli lemma shows that, a.s., $\{\eta_{2^k}>n(k),\tilde{T}_{m_2(k)}=\infty\}$, $k=1,2,3,\dots $, occur finitely often. By Corollary \ref{limx2Delta}, a.s. $\{\tilde{T}_{m_2(k)}<\infty\}$, $k=1,2,3,\dots $, occur finitely often. Therefore, $\eta_{2^k} \leq n(k)$ for all but finitely many $k$. Thus, almost surely, for all sufficiently large $x>0$, say $x>C(\omega)$ where $C(\omega)$ is a positive number which depends on the sample path,  
$$
\eta_x \leq \eta_{2^{\lfloor \log_2 x \rfloor+1}} \leq n\left( \lfloor \log_2 x \rfloor+1 \right) \leq 4x^2 (\log (2x))^{\frac{8}{\delta}+1+3 \varepsilon}.
$$
Proposition \ref{hightranprop} shows that, a.s., for all but finitely many $n$, $\|\bm{S}_n\|>C(\omega)$, and thus,
\begin{equation}
  \label{nSnlogSn5}
  n \leq \eta_{\|\bm{S}_n\|} \leq 4\|\bm{S}_n\|^2\left(\log 2\|\bm{S}_n\|\right)^{\frac{8}{\delta}+1+3 \varepsilon} .
\end{equation}
If $\|\bm{S}_n\| \leq n^{\frac{1}{2}}(\log n)^{-\frac{4}{\delta}-\frac{2}{3}}$ for some large $n$, then,
$$
4\|\bm{S}_n\|^2\left(\log 2\|\bm{S}_n\|\right)^{\frac{8}{\delta}+1+3\varepsilon}  \leq \frac{4n }{(\log n)^{\frac{1}{3}-3\varepsilon}} ,
$$
which can occur for only finitely many $n$ in view of (\ref{nSnlogSn5}). This completes the proof of (\ref{rateestescap}) for the case $\delta<4$. If $\delta\geq 4$, we choose a positive number $\delta_1<4$ such that 
$$
\frac{4}{\delta_1}+\frac{2}{3}<2.
$$
Since the step distribution $\mu$ satisfies Assumption $A(2+\delta_1)$ with $\delta_1<4$, by the results for the case $\delta<4$, a.s. for all large $n$,
$$\|\bm{S}_n\| > n^{\frac{1}{2}}(\log n)^{-\frac{4}{\delta_1}-\frac{2}{3}}>n^{\frac{1}{2}}(\log n)^{-2},$$ 
whence (\ref{rateestescap}) follows.
\end{proof}

\subsection{Proof for dimension 2}
\label{rec2Dproof}

In this section, we prove the following result.

\begin{proposition}
  \label{2DSRRWprop}
  Let $\bm{S}$ be an SRRW in $\RR^2$ with parameter $\alpha\leq 1/2$ and step distribution $\mu$. Assume that $\mu$ satisfies Assumption $A(2+\delta)$ for some $\delta>0$ and $\EE \bm{X}_1\bm{X}_1^T$ is the identity matrix. 
  \begin{enumerate}[(i)]
      \item If $\alpha<1/2$, then $\bm{S}$ is recurrent. 
\item If $\alpha=1/2$, then almost surely,
   $$
  \lim_{n\to \infty} \frac{\log \|\bm{S}_n\|^2}{\log n}=1.
  $$
  \end{enumerate}
\end{proposition}

Let us introduce some notation that will be used throughout this section, i.e., Section \ref{rec2Dproof}. We assume that $\bm{S},\alpha,\mu,\delta$ are as described in Proposition \ref{2DSRRWprop}. Without loss of generality, we assume that $\delta\in (0,1)$. We let 
\begin{equation}
    \label{defnueps}
  \nu:=\frac{\delta}{8+4\delta}, \quad \varepsilon:=\frac{\delta}{8},  
\end{equation}
and fix $\kappa \in (0,1)$ such that 
\begin{equation}
  \label{kappavarepass}
 \kappa+\frac{\nu}{2}>1,\quad \kappa (1+\frac{\varepsilon}{2})>1,\quad \kappa (1+\frac{\delta}{4})(1-\varepsilon)>1.
\end{equation}
If $\alpha<1/2$, we shall further assume that 
\begin{equation}
    \label{defkappatilde}
   \tilde{\kappa}:=\frac{1}{3}+\frac{\kappa}{6}=\frac{1}{2}-\frac{1-\kappa}{6}>\alpha. 
\end{equation}
Note that such a $\kappa$ exists since (\ref{kappavarepass}) and (\ref{defkappatilde}) hold when $\kappa$ is sufficiently close to $1$.

For $n>1$, let
\begin{equation}
    \label{defxndn}
    x_n:=\frac{\log (\|\bm{S}_n\|^2+n^{\kappa})}{\log n}, \quad 
D_n:=\sqrt{\|\bm{S}_n\|^2+n^{\kappa}}.
\end{equation}
Notice that $D_n\geq \|\bm{S}_n\|$ and $D_n\geq n^{\kappa/2}$. Observe that 
\begin{equation}
  \label{delxnformula}
  \begin{aligned}
       x_{n+1}-x_n&=\frac{\log D_{n+1}^2- \log D_n^2 }{\log (n+1)}-\frac{\log D_n^2 }{\log n\log (n+1)}\log (1+\frac{1}{n}) \\ 
       &=\frac{\log (1+u_{n+1}) }{\log (n+1)}-\frac{\log D_n^2  }{\log n\log (n+1)}\log (1+\frac{1}{n}) 
  \end{aligned}
\end{equation}
where
\begin{equation}
  \label{undef}
  u_{n+1}:=\frac{2\bm{S}_n\cdot \bm{X}_{n+1}+\|\bm{X}_{n+1}\|^2+(n+1)^{\kappa}-n^{\kappa}}{D_n^2 }, \quad n\geq 1.
\end{equation}
For $n> 1$, we write 
\begin{equation}
    \label{defInun}
   \begin{aligned}
 I_n^{(1)}&:=u_{n+1}-\frac{2\alpha\|\bm{S}_n\|^2}{nD_n^2 }-\frac{2}{D_n^2},   \\
 I_n^{(2)}&:=u_{n+1}^2\mathds{1}_{\{\|\bm{X}_{n+1}\|\leq D_n^{1-\varepsilon}\}} - \frac{4\|\bm{S}_n\|^2}{D_n^4},   \\
  I_n^{(3)}&:=\frac{2\alpha\|\bm{S}_n\|^2}{nD_n^2 }+\frac{2n^{\kappa}}{D_n^4}-\frac{\log D_n^2}{n \log n}.
\end{aligned} 
\end{equation}

For the reader's convenience, here we briefly explain the main idea of the proof of Proposition \ref{2DSRRWprop}. As will be shown later, by Taylor's expansion and truncation, we can write
$$
x_{n+1}-x_n=\frac{1}{\log (n+1)}\left(I_n^{(1)} - \frac{1}{2} I_n^{(2)} + I_n^{(3)}+I_n^{(err)}\right), \quad n>1,
$$
where $I_n^{(err)}$ corresponds to the error terms caused by Taylor's approximation and truncation. We shall use Lemma \ref{lemcondun} below to control (the conditional expectations of) $I_n^{(1)}, I_n^{(2)}$ and $I_n^{(err)}$. Lemma \ref{lemIn3sum} explains why the proof for the diffusive regime would fail in the critical regime, and vice versa. Non-rigorously speaking:
\begin{itemize}
    \item If $\alpha<1/2$ and $\|\bm{S}_n\|\geq n^{\tilde{\kappa}}$, then $I_n^{(3)}$ is negative and $|I_n^{(3)}|$ is sufficiently large, which forces $(x_n)_{n>1}$ to decrease until it reaches below $\tilde{\kappa}$. Then randomness (more precisely, the optional stopping theorem) ensures that $(\bm{S}_n)_{n\in \NN}$ returns to some fixed ball infinitely often a.s..
    \item If $\alpha =1/2$, then the contribution of $(I_n^{(3)})_{n>1}$ is asymptotically negligible. We can then show that $(x_n)_{n>1}$ is a.s. convergent by the quasi-martingale convergence theorem. The limit is a.s. equal to $1$ by the law of the iterated logarithm.
\end{itemize}

\begin{lemma}
  \label{lemcondun}
 Under the setting of Proposition \ref{2DSRRWprop}, let $(u_{n+1})_{n\geq 1},(I_n^{(i)})_{1\leq i\leq 3,n>1}$ be as in  (\ref{undef}) and (\ref{defInun}). For any $n\geq 1$, we have,
\begin{equation}
  \label{condEun}
  \EE (u_{n+1}\mid \FF_n)=\frac{2\alpha\|\bm{S}_n\|^2}{nD_n^2 }+\frac{2}{D_n^2}+\frac{\alpha\Delta_n(\|\bm{x}\|^2)}{D_n^2}+\frac{(n+1)^{\kappa}-n^{\kappa}}{D_n^2},
\end{equation}
and
\begin{equation}
  \label{condEun2low}
\begin{aligned}
  \EE (u_{n+1}^2\mathds{1}_{\{\|\bm{X}_{n+1}\|\leq D_n^{1-\varepsilon}\}} \mid \FF_n) &\geq \frac{4\|\bm{S}_n\|^2}{D_n^4}-\frac{4\alpha\|\Delta_n(\bm{x}\bm{x}^T)\|}{D_n^2}-\frac{4 }{n^{1-\kappa}D_n^{2+\varepsilon}} \\
  &\quad-\frac{8 \EE \|\bm{X}_1\|^{2+\frac{\delta}{2}} +8\alpha\Delta_n(\|\bm{x}\|^{2+\frac{\delta}{2}}) }{D_n^{2+\frac{\delta}{2}(1-\varepsilon)}}.
\end{aligned}
\end{equation}
If one further assumes that $\alpha=1/2$, then,
\begin{equation}
  \label{condEun2up}
\begin{aligned}
  \EE (u_{n+1}^2\mathds{1}_{\{\|\bm{X}_{n+1}\|\leq D_n^{1-\varepsilon}\}} \mid \FF_n) &\leq \frac{4\|\bm{S}_n\|^2}{D_n^4}+\frac{2\|\Delta_n(\bm{x}\bm{x}^T)\|}{D_n^2}+\frac{5+ \Delta_n(\|\bm{x}\|^2)}{D_n^4} \\
  &\quad+\frac{10 \EE \|\bm{X}_1\|^{2+\frac{\delta}{2}}+5\Delta_n(\|\bm{x}\|^{2+\frac{\delta}{2}}) }{2D_n^{2+\frac{\delta}{2}(1-\varepsilon)+\varepsilon}}+\frac{4}{n^{1-\kappa}D_n^{2+\varepsilon}},
\end{aligned}
\end{equation}
and 
\begin{equation}
  \label{condEun2remsum}
    \sum_{n=2}^{\infty}\EE \left(\EE(I_n^{(1)}\mid \FF_n)^-\right) <\infty,\quad  \sum_{n=2}^{\infty}\EE\left( \EE (I_n^{(2)}\mid \FF_n)^+ \right) <\infty,
\end{equation}   
and
\begin{equation}
  \label{un2Dvarsum}
  \sum_{n=1}^{\infty}\EE \left(\frac{u_{n+1}^2}{D_n^{\varepsilon}} \mathds{1}_{\{\|\bm{X}_{n+1}\| \leq D_n^{1-\varepsilon}\} }\right)<\infty.
\end{equation}
\end{lemma}
\begin{proof}
Recall that $\EE \bm{X}_1\bm{X}_1^T$ is the $2 \times 2$ identity matrix. Equation (\ref{condEun}) follows from (\ref{undef}) and Lemma \ref{condlemXn1}. Now using that $\|\bm{S}_n\|\leq D_n$ and $0< (n+1)^{\kappa}-n^{\kappa}\leq  n^{\kappa-1}$ for all $n\geq 1$, we have, by (\ref{undef}),
  \begin{equation}
    \label{un2inelow}
    \begin{aligned}
           &\quad\ u_{n+1}^2\mathds{1}_{\{\|\bm{X}_{n+1}\|\leq D_n^{1-\varepsilon}\}} \\
           &\geq \frac{4(\bm{S}_n\cdot \bm{X}_{n+1})^2-4\|\bm{S}_{n}\|\|\bm{X}_{n+1}\|^3-4C\|\bm{S}_{n}\|\|\bm{X}_{n+1}\|n^{\kappa-1}}{D_n^4}\mathds{1}_{\{\|\bm{X}_{n+1}\|\leq D_n^{1-\varepsilon}\}}   \\
           &\geq \frac{4(\bm{S}_n\cdot \bm{X}_{n+1})^2}{D_n^4}- \frac{4 \|\bm{X}_{n+1}\|^2}{D_n^2} \mathds{1}_{\{\|\bm{X}_{n+1}\| > D_n^{1-\varepsilon}\}}-\frac{4\|\bm{X}_{n+1}\|^{2+\frac{\delta}{2}}}{D_n^{3-(3-2-\frac{\delta}{2})(1-\varepsilon)}}   -\frac{4 }{n^{1-\kappa}D_n^{2+\varepsilon}} \\
            &\geq \frac{4(\bm{S}_n\cdot \bm{X}_{n+1})^2}{D_n^4}- \frac{4 \|\bm{X}_{n+1}\|^{2+\frac{\delta}{2}}}{D_n^{2+\frac{\delta}{2}(1-\varepsilon)}} -\frac{4\|\bm{X}_{n+1}\|^{2+\frac{\delta}{2}}}{D_n^{2+\frac{\delta}{2}(1-\varepsilon)+\varepsilon}}   -\frac{4 }{n^{1-\kappa}D_n^{2+\varepsilon}}\\
            &\geq \frac{4(\bm{S}_n\cdot \bm{X}_{n+1})^2}{D_n^4}- \frac{8 \|\bm{X}_{n+1}\|^{2+\frac{\delta}{2}}}{D_n^{2+\frac{\delta}{2}(1-\varepsilon)}} -\frac{4 }{n^{1-\kappa}D_n^{2+\varepsilon}}.
    \end{aligned}
  \end{equation}
By Lemma \ref{condlemXn1} and our assumption that $\EE \bm{X}_1\bm{X}_1^T$ is the $2 \times 2$ identity matrix, 
$$
\begin{aligned}
 \EE(\|\bm{X}_{n+1}\|^{2+\frac{\delta}{2}}\mid \FF_n)&=\EE \|\bm{X}_1\|^{2+\frac{\delta}{2}} +\alpha\Delta_n(\|\bm{x}\|^{2+\frac{\delta}{2}}), \\
\EE ((\bm{S}_n\cdot \bm{X}_{n+1})^2\mid \FF_n)&=\|\bm{S}_n\|^2+\alpha\bm{S}_n^T\Delta_n(\bm{x}\bm{x}^T) \bm{S}_n.   
\end{aligned}
$$
Then (\ref{condEun2low}) follows from (\ref{un2inelow}) by using that $|\bm{S}_n^T\Delta_n(\bm{x}\bm{x}^T) \bm{S}_n|\leq \|\Delta_n(\bm{x}\bm{x}^T)\|\|\bm{S}_n\|^2$ and $\|\bm{S}_n\|\leq D_n$. 

Now assume that $\alpha=1/2$. Using similar arguments as in (\ref{un2inelow}), one has,
  \begin{equation}
    \label{un2ineup}
    \begin{aligned}
           &\quad\ u_{n+1}^2\mathds{1}_{\{\|\bm{X}_{n+1}\|\leq D_n^{1-\varepsilon}\}} \\
           &\leq \frac{4(\bm{S}_n\cdot \bm{X}_{n+1})^2+\|\bm{X}_{n+1}\|^4+1+4\|\bm{S}_{n}\|\|\bm{X}_{n+1}\|^3}{D_n^4}\mathds{1}_{\{\|\bm{X}_{n+1}\|\leq D_n^{1-\varepsilon}\}}  \\
         &\quad+\frac{4\|\bm{S}_{n}\|\|\bm{X}_{n+1}\|n^{\kappa-1}+2\|\bm{X}_{n+1}\|^2}{D_n^4}\mathds{1}_{\{\|\bm{X}_{n+1}\|\leq D_n^{1-\varepsilon}\}} \\
           &\leq \frac{4(\bm{S}_n\cdot \bm{X}_{n+1})^2+2\|\bm{X}_{n+1}\|^{2}+1}{D_n^4}+\frac{\|\bm{X}_{n+1}\|^{2+\frac{\delta}{2}}}{D_n^{4-(4-2-\frac{\delta}{2})(1-\varepsilon)}}  +\frac{4\|\bm{X}_{n+1}\|^{2+\frac{\delta}{2}}}{D_n^{3-(3-2-\frac{\delta}{2})(1-\varepsilon)}}+\frac{4 }{n^{1-\kappa}D_n^{2+\varepsilon}} \\
             &\leq \frac{4(\bm{S}_n\cdot \bm{X}_{n+1})^2+2\|\bm{X}_{n+1}\|^{2}+1}{D_n^4}  +\frac{5\|\bm{X}_{n+1}\|^{2+\frac{\delta}{2}}}{D_n^{2+\frac{\delta}{2}(1-\varepsilon)+\varepsilon}} +\frac{4 }{n^{1-\kappa}D_n^{2+\varepsilon}}.
    \end{aligned}
  \end{equation}
Then, similarly as in the proof of (\ref{condEun2low}), one can deduce (\ref{condEun2up}) from (\ref{un2ineup}) and Lemma \ref{condlemXn1}, where we used that
$$
\EE \|\bm{X}_1\|^2=2, \quad \EE(\|\bm{X}_{n+1}\|^2\mid \FF_n)=2+\frac{\Delta_n(\|\bm{x}\|^2)}{2}.
$$
Now (\ref{condEun}) shows that, for $\alpha=1/2$,
 $$\EE (I_n^{(1)}\mid \FF_n)\geq -\frac{|\Delta_n(\|\bm{x}\|^2)|}{2D_n^2}.$$
By Corollary \ref{limx2Delta}, there exists a constant $C$ such that for all $n\geq 1$, 
\begin{equation}
  \label{EInHnest}
  \max\{\EE |\Delta_n(\|\bm{x}\|^2)|,\EE \|\Delta_n(\bm{x}\bm{x}^T)\|,\EE |\Delta_n(\|\bm{x}\|^{2+\frac{\delta}{2}})| \}\leq C n^{-\nu}.
\end{equation}
By the choice of $\kappa$, i.e. (\ref{kappavarepass}),
$$\sum_{n=2}^{\infty}\EE \left(\EE(I_n^{(1)}\mid \FF_n)^-\right)\leq \sum_{n=1}^{\infty} \EE \frac{|\Delta_n(\|\bm{x}\|^2)|}{2D_n^2}\leq \sum_{n=1}^{\infty} \frac{C}{2}n^{-\kappa-\nu}<\infty,$$
which proves the first inequality in (\ref{condEun2remsum}). The second inequality in (\ref{condEun2remsum}) and (\ref{un2Dvarsum}) are proved similarly.
\end{proof}

\begin{lemma}
  \label{lemIn3sum}
  Under the setting of Proposition \ref{2DSRRWprop}, let $(I_n^{(i)})_{1\leq i\leq 3,n>1}$ be as in (\ref{defInun}). One has: \\
  (i) If $\alpha<1/2$ and $\|\bm{S}_n\|\geq n^{\tilde{\kappa}}$ where $\tilde{\kappa}$ is defined by (\ref{defkappatilde}), then 
  $$
I_n^{(3)}< -\frac{2(\tilde{\kappa}-\alpha)}{n} +\frac{2}{n^{1+\frac{1-\kappa}{3}}}.
  $$
 (ii) If $\alpha=1/2$, then 
  \begin{equation}
    \label{sumvn-bd}
    \EE \sum_{n=2}^{\infty} \frac{(I_n^{(3)})^{-}}{\log (n+1)} <\infty.
  \end{equation}
\end{lemma}
\begin{proof}
 (i) Since $D_n^2\geq \|\bm{S}_n\|^2\geq n^{2\tilde{\kappa}}$, by definition (\ref{defInun}),
$$
I_n^{(3)} \leq \frac{2\alpha}{n}+\frac{2}{n^{4\tilde{\kappa}-\kappa}}-\frac{2\tilde{\kappa}}{n}  \leq -\frac{2(\tilde{\kappa}-\alpha)}{n} +\frac{2}{n^{1+\frac{1-\kappa}{3}}}.
$$
(ii) We shall divide the left-hand side of (\ref{sumvn-bd}) into three parts. Let 
\begin{equation}
    \label{defkappa1}
    \kappa_1:=\frac{1+\kappa}{2} \in (\kappa,1).
\end{equation}
(a). If $ n^{\kappa_1} \leq \|\bm{S}_n\|^2 \leq n \log^2 n$, then $D_n^2 \leq n( \log^2 n+1)$ and thus
$$
\begin{aligned}
  \frac{\|\bm{S}_n\|^2}{D_n^2 }-\frac{\log D_n^2}{ \log n}&\geq 1-\frac{n^{\kappa}}{\|\bm{S}_n\|^2+n^{\kappa} } - \frac{\log (n( \log^2 n+1))}{ \log n}\\
  &\geq -\frac{1}{n^{\kappa_1-\kappa}+1}-\frac{\log ( \log^2 n+1) }{\log n},
\end{aligned}
$$
whence we have
\begin{equation}
  \label{vn-sn910}
  \begin{aligned}
     &\quad\ \EE \sum_{n=2}^{\infty}\frac{(I_n^{(3)})^{-}}{\log (n+1)} \mathds{1}_{\{n^{\kappa_1} \leq \|\bm{S}_n\|^2 \leq n \log^2 n\}} \\
     &\leq \sum_{n=2}^{\infty} \frac{1}{n\log (n+1)} \left(\frac{1}{n^{\kappa_1-\kappa}+1}+\frac{\log ( \log^2 n+1) }{\log n}\right)<\infty.
  \end{aligned}
\end{equation}
(b). Using that $\log(1+x) \leq \sqrt{x}$ for sufficiently large $x$ and that $n^{\kappa-1}\leq 1$ for all $n\geq 1$, it follows that for all large $n$, say $n\geq m_1$, we have
$$
\frac{\log (\frac{\|\bm{S}_n\|^2}{n}+\frac{n^{\kappa}}{n})}{n \log n \log (n+1)}\mathds{1}_{\{\|\bm{S}_n\|^2 > n \log^2 n\}} \leq \mathds{1}_{\{\|\bm{S}_n\|^2 > n \log^2 n\}}\frac{ \|\bm{S}_n\| }{n^{3/2} \log^2 n}.
$$
In particular,
$$
\begin{aligned}
\frac{(I_n^{(3)})^{-}}{\log (n+1)} \mathds{1}_{\{\|\bm{S}_n\|^2 > n \log^2 n\}} &\leq \frac{\log D_n^2}{n \log n\log (n+1)} \mathds{1}_{\{\|\bm{S}_n\|^2 > n \log^2 n\}} \\
&= \left(\frac{ 1}{n \log (n+1)} + \frac{\log (\frac{\|\bm{S}_n\|^2}{n}+\frac{n^{\kappa}}{n})}{n \log n \log (n+1)}\right)\mathds{1}_{\{\|\bm{S}_n\|^2 > n \log^2 n\}} \\
&\leq \left(\frac{ 1 }{n \log (n+1)}+\frac{ \|\bm{S}_n\| }{n^{3/2} \log^2 n}\right)\mathds{1}_{\{\|\bm{S}_n\|^2 > n \log^2 n\}} \\
&\leq \left(\frac{ 1 }{n \log (n+1)}+\frac{ \|\bm{S}_n\|^2 }{n^{2} \log^3 n}\right)\mathds{1}_{\{\|\bm{S}_n\|^2 > n \log^2 n\}}.
\end{aligned}
$$
Consequently, by Chebyshev's inequality and  \cite[Proposition 2.1]{MR4490996} (which shows that $\EE \|\bm{S}_n\|^2 =O(n \log n)$), one has
\begin{equation}
  \label{vn-snlog2n}
  \begin{aligned}
     \EE \sum_{n=m_1}^{\infty}\frac{(I_n^{(3)})^{-}}{\log (n+1)} \mathds{1}_{\{\|\bm{S}_n\|^2 > n \log^2 n\}} &\leq \sum_{n=m_1}^{\infty}\left(\frac{\PP(\|\bm{S}_n\|^2 > n \log^2 n)}{n\log (n+1)}+\frac{\EE \|\bm{S}_n\|^2}{n^2\log^3 n} \right)\\
& \leq \sum_{n=m_1}^{\infty} \frac{2\EE \|\bm{S}_n\|^2}{n^2 \log^3 n } <\infty.
  \end{aligned}
\end{equation}
(c). By the choice of $\kappa_1$, i.e. (\ref{defkappa1}),
$$
\lim_{n\to \infty}n\left(\frac{2n^{\kappa}}{(n^{\kappa_1}+n^{\kappa})^2} -\frac{\log (n^{\kappa_1}+n^{\kappa})}{n\log n}\right)=2-\kappa_1>0.
$$
Therefore, for all large $n$ such that $\|\bm{S}_n\|^2 < n^{\kappa_1}$,
$$
  \frac{2n^{\kappa}}{(\|\bm{S}_n\|^2+n^{\kappa})^2}-\frac{\log (\|\bm{S}_n\|^2+n^{\kappa})}{n \log n} >   \frac{2n^{\kappa}}{(n^{\kappa_1}+n^{\kappa})^2} -\frac{\log (n^{\kappa_1}+n^{\kappa})}{n\log n}  > 0,
$$
which implies that there exists a positive constant $m_2$ such that for any $n\geq m_2$,
\begin{equation}
  \label{vnSnleqn910}
I_n^{(3)} \mathds{1}_{\{\|\bm{S}_n\|^2 < n^{\kappa_1}\}}\geq 0.
\end{equation}

Now (\ref{sumvn-bd}) follows from  (\ref{vn-sn910}), (\ref{vn-snlog2n}) and (\ref{vnSnleqn910}).
\end{proof}

\subsubsection{Proof for the diffusive regime}
\label{diffu2d}

In this section, we prove Proposition \ref{2DSRRWprop} (i). We first prove a weaker version.

\begin{proposition}
  \label{Snsio}
  Under the setting of Proposition \ref{2DSRRWprop} (i), let $\nu,\varepsilon,\kappa,\tilde{\kappa}, (x_n)_{n> 1}$ be as in (\ref{defnueps}), (\ref{kappavarepass}), (\ref{defkappatilde}) and (\ref{defxndn}). Then, a.s., $\|\bm{S}_n\|\leq n^{\tilde{\kappa}}$ infinitely often.
\end{proposition}
\begin{proof}
  For $k\geq 1$, define 
  $$\tau_k:=\inf\{n\geq k: \|\bm{S}_n\| \leq n^{\tilde{\kappa}} \}.$$
  It suffices to show that for any $k\geq 1$, $\tau_k<\infty$ a.s.. Indeed, since $(\tau_k)_{k\geq 1}$ is non-decreasing in $k$, it suffices to prove that $\tau_k<\infty$ a.s. for all sufficiently large $k$. 

Define 
  \begin{equation}
    \label{defTk}
      T_k:=\inf\{n \geq k:  \max\{|\Delta_n(\|\bm{x}\|^{2})|,|\Delta_n(\|\bm{x}\|^{2+\frac{\delta}{2}})|,\|\Delta_n(\bm{x}\bm{x}^T)\|\}\geq \frac{1}{n^{\nu}} \ \text{or}\ \|\bm{S}_n\|\geq n\}.
  \end{equation}
Notice that by Corollary \ref{limx2Delta} and the strong law of large numbers (\ref{llnSRRW}), for any $k>1$,
\begin{equation}
\label{Tmeveninfty}
    \PP(\bigcup_{m\geq k} \{T_m=\infty\} )=1.
\end{equation}

Now fix $k\geq 1$, define a process $(z_n)_{n\geq k}$ as follows:
\begin{equation}
    \label{defznk}
  z_{n}:=x_{n\wedge \tau_k\wedge T_k}+\sum_{j=k}^{n\wedge \tau_k\wedge T_k-1}\frac{\tilde{\kappa}-\alpha}{j\log (j+1)},\quad n\geq k,  
\end{equation}
  with the convention that $z_k:=x_k$. We shall prove in Lemma \ref{lemznsupmart} that, if $k$ is sufficiently large, then $(z_n)_{n\geq k}$ is a lower-bounded supermartingale, and thus, converges a.s.. On the other hand, on the event $\{\tau_k=\infty\} \cap \{T_k=\infty\}$, as $n\to \infty$,
  $$
 z_n\geq \kappa + \sum_{j=k}^{n-1}\frac{\tilde{\kappa}-\alpha}{j\log (j+1)}  \to \infty,
  $$
  where we used that $\tilde{\kappa}-\alpha>0$ and that
  $$
\sum_{j=k}^{n-1} \frac{1}{j\log (j+1)} \geq \sum_{j=k}^{n-1} \int_{j+1}^{j+2} \frac{1}{y\log y}dy =\log \log (n+1)- \log \log (k+1).
  $$
  Thus, $\PP(\{\tau_k=\infty\} \cap \{T_k=\infty\}) =0$ for all sufficiently large $k$. Indeed, we have proved that for any $m\geq k$,
  \begin{equation}
    \label{taukTmprob0}
      \PP(\{\tau_k=\infty\} \cap \{T_m=\infty\}) \leq \PP(\{\tau_m=\infty\} \cap \{T_m=\infty\}) =0,
  \end{equation}
  where we used the fact that $\{\tau_k =\infty\} \subset \{\tau_m=\infty\}$. One can then deduce from  (\ref{Tmeveninfty}) and (\ref{taukTmprob0}) that $\PP(\tau_k<\infty)=1$.
\end{proof}

\begin{lemma}
\label{lemznsupmart}
  There exists a positive integer $K$ such that if $k\geq K$, then $(z_n)_{n\geq k}$ defined in (\ref{defznk}) is a lower-bounded supermartingale.  
\end{lemma}
\begin{proof}
We shall give upper bounds for 
$\EE (x_{n+1}-x_n\mid \FF_n)$ where $n\in [k,\tau_k\wedge T_k)$. We consider the cases $\|\bm{X}_{n+1}\| \leq D_n^{1-\varepsilon}$ and $\|\bm{X}_{n+1}\| > D_n^{1-\varepsilon}$ separately.

    By Taylor expansion, there exists an $\varepsilon_1>0$ such that if $|x|\leq \varepsilon_1$, then
\begin{equation}
 \label{tayineln1x}
 \log (1+x) \leq x -\frac{1}{2}x^2+|x|^3.
\end{equation} 
Recall $u_{n+1}$ defined in (\ref{undef}). By definition,
\begin{equation}
  \label{unDbd}
  |u_{n+1}|\leq \frac{2\|\bm{X}_{n+1}\|}{D_n}+\frac{\|\bm{X}_{n+1}\|^2+1}{D_n^2}.
\end{equation}
If $\|\bm{X}_{n+1}\| \leq D_n^{1-\varepsilon}$, then
\begin{equation}
  \label{unDvarbd}
  |u_{n+1}|\leq \frac{4}{D_n^{\varepsilon}},
\end{equation}
and, in particular, $|u_{n+1}|<\varepsilon_1$ for all $n\geq k$ if $k$ is sufficiently large. Thus, by (\ref{delxnformula}) and (\ref{tayineln1x}), for all $n\geq k$ such that $\|\bm{X}_{n+1}\| \leq D_n^{1-\varepsilon}$, 
\begin{equation}
  \label{xnrecura12}
  x_{n+1}-x_n\leq \frac{1}{\log (n+1)}\left(u_{n+1}-\frac{1}{2}u_{n+1}^2+|u_{n+1}|^3\right)-\frac{\log  D_n^2}{(n+1) \log n \log (n+1)}.
\end{equation}
Observe that 
$$
  \log D_{n+1}^2 - \log D_n^2 = \log (1+u_{n+1}) \leq u_{n+1},
$$
and thus, by (\ref{delxnformula}), no matter whether $\|\bm{X}_{n+1}\| \leq D_n^{1-\varepsilon}$ or not, one has
\begin{equation}
  \label{delxngener}
  x_{n+1}-x_n \leq \frac{u_{n+1}}{\log (n+1)}-\frac{\log D_n^2}{(n+1) \log n \log (n+1)}.
\end{equation}
By (\ref{unDvarbd}), (\ref{xnrecura12}) and (\ref{delxngener}), we have
\begin{equation}
    \label{condxnlogupbd}
    \begin{aligned}
&\quad \ (x_{n+1}-x_n)\log(n+1)  \\
   &\leq \mathds{1}_{\{\|\bm{X}_{n+1}\| \leq D_n^{1-\varepsilon}\} }\left(u_{n+1}- \frac{1}{2}u_{n+1}^2+\frac{4u_{n+1}^2 }{D_n^{\varepsilon}}\right) -\frac{\log D_n^2}{n \log n}+u_{n+1}\mathds{1}_{\{\|\bm{X}_{n+1}\| > D_n^{1-\varepsilon}\} } \\  &=\left(u_{n+1}- \frac{1}{2}u_{n+1}^2\mathds{1}_{\{\|\bm{X}_{n+1}\| \leq D_n^{1-\varepsilon}\} }-\frac{\log D_n^2}{n \log n }\right)+\frac{4u_{n+1}^2 }{D_n^{\varepsilon}}\\
   &=I_n^{(1)} - I_n^{(2)}+\frac{4u_{n+1}^2 }{D_n^{\varepsilon}} + I_n^{(3)}.
\end{aligned}
\end{equation}
  By (\ref{condEun}), if $k\leq n<\tau_k\wedge T_k$, then 
  \begin{equation}
      \label{condeIn1upbd}
      \EE (I_n^{(1)} \mid \FF_n) \leq \frac{\alpha}{n^{\kappa+\frac{\nu}{2}}}+\frac{1}{n^{1-\kappa}\|\bm{S}_n\|^2} \leq  \frac{\alpha}{n^{\kappa+\frac{\nu}{2}}}+\frac{1}{n^{1-\kappa+2\tilde{\kappa}}}.
  \end{equation}
  Note that both $\kappa+\frac{\nu}{2}$ and $1-\kappa+2\tilde{\kappa}$ are greater than $1$ by (\ref{kappavarepass}) and (\ref{defkappatilde}). 
Similarly, by (\ref{condEun2low}), one can show that there exist positive constants $C_1,C_2>1$ such that for all $ n\in [k,\tau_k\wedge T_k)$,
 \begin{equation}
      \label{condeminusIn2upbd}
\EE (-I_n^{(2)} \mid \FF_n) \leq \frac{C_1}{n^{C_2}}, \quad \EE ( \frac{4u_{n+1}^2 }{D_n^{\varepsilon}}\mid \FF_n) \leq \frac{C_1}{n^{C_2}},
  \end{equation}
where we also used the assumption (\ref{kappavarepass}) to obtain that
$$
\kappa (1+\frac{\delta}{4}(1-\varepsilon))>1,\quad \kappa (1+\frac{\varepsilon}{2})>1.
$$
In summary, if $k \leq n<\tau_k\wedge T_k$ and , then by (\ref{condxnlogupbd}), (\ref{condeIn1upbd}), (\ref{condeminusIn2upbd}) and Lemma \ref{lemIn3sum} (i),
\begin{equation}
    \label{condupxnsupmart}
    \EE (x_{n+1}-x_n\mid \FF_n) \leq \frac{1}{\log(n+1)}\left(  \frac{\alpha}{n^{\kappa+\frac{\nu}{2}}}+\frac{1}{n^{1-\kappa+2\tilde{\kappa}}}+\frac{2C_1}{n^{C_2}}-\frac{2(\tilde{\kappa}-\alpha)}{n} +\frac{2}{n^{1+\frac{1-\kappa}{3}}}\right).
\end{equation}
The right-hand side of (\ref{condupxnsupmart}) is less than $(\tilde{\kappa}-\alpha)(n\log(n+1))^{-1}$ if $k$ is sufficiently large, say $k\geq K$. In particular, $(z_n)_{n\geq k}$ is a supermartingale if $k\geq K$. It remains to note that $(z_n)_{n\geq k}$ is lower bounded by $\kappa$.
\end{proof}

We also need Lemma \ref{fxdelta2dest} given in Section \ref{Lyapunovf}.

\begin{proof}[Proof of Lemma \ref{fxdelta2dest}.]
  We may assume that $\|\bm{x}+\bm{y}\|\geq \|\bm{x}\|\geq 1$, otherwise (\ref{fxdeltagener}) is trivial. Using that $\sqrt{a+b} -\sqrt{a} \leq \sqrt{b}$ for any $a,b\geq 0$,
  we have
  $$
    \begin{aligned}
    \sqrt{\log \|\bm{x}+\bm{y}\|}-\sqrt{\log \|\bm{x}\|} \leq \sqrt{\frac{1}{2}} \sqrt{\log \left(\frac{\|\bm{x}+\bm{y}\|^2}{\|\bm{x}\|^2}\right)}  &\leq \sqrt{\frac{1}{2}} \sqrt{\log \left(\frac{2\|\bm{x}\|^2+2\|\bm{y}\|^2}{\|\bm{x}\|^2}\right)}\\
     &\leq  \sqrt{\frac{\log 2}{2}} + \sqrt{\frac{1}{2}}\frac{\|\bm{y}\|}{\|\bm{x}\|},
  \end{aligned}
  $$
  where we used that $\log (1+t) \leq t$ for all $t\geq 0$, which completes the proof of (\ref{fxdeltagener}). 
  
  Now observe that 
  \begin{equation}
      \label{fxplusyfxexpan}
    \begin{aligned}
  \sqrt{\log \|\bm{x}+\bm{y}\|}-\sqrt{\log \|\bm{x}\|}&=\sqrt{\log \|\bm{x}\|} \left(\left(1+\frac{1}{\log \|\bm{x}\|^2}\log \sqrt{1+\frac{2 \bm{x} \cdot \bm{y}+\|\bm{y}\|^2}{\|\bm{x}\|^2}}\right)^{\frac{1}{2}}-1\right).
\end{aligned}  
  \end{equation}
  Recall the following Taylor expansions: As $x\to 0$,
 \begin{equation}
     \label{taylorlog1xsqrt1x}
      \log(1+x)= x-\frac{x^2}{2}+O(|x|^3); \quad \sqrt{1+x}=1+\frac{x}{2}-\frac{x^2}{8}+O(|x|^3).
 \end{equation}
 By (\ref{fxplusyfxexpan}) and (\ref{taylorlog1xsqrt1x}), one can show that for any $\varepsilon\in (0,1)$, there exist positive constants $r$ and $C$ such that if $\|\bm{x}\|\geq r$ and $\bm{y} \in E_{\varepsilon}(\bm{x})$, then
$$
   \begin{aligned}
    f(\bm{x}+\bm{y})-f(\bm{x}) &\leq \sqrt{\log \|\bm{x}\|}\left(\frac{1}{2 \log \|\bm{x}\|}\left(\frac{ \bm{x} \cdot \bm{y}}{\|\bm{x}\|^2}+\frac{\|\bm{y}\|^2}{2\|\bm{x}\|^2}-\frac{(\bm{x} \cdot \bm{y})^2}{\|\bm{x}\|^4}+\frac{C\|\bm{y}\|^3}{\|\bm{x}\|^{3}}\right)\right. \\ 
    &\quad \left.-\frac{1}{8 \log^2 \|\bm{x}\|} \frac{(\bm{x} \cdot \bm{y})^2}{\|\bm{x}\|^4}+\frac{C\|\bm{y}\|^3}{\|\bm{x}\|^{3}\log^2 \|\bm{x}\|} \right) \\
  & \leq \sqrt{\log \|\bm{x}\|}\left(\frac{1}{2 \log \|\bm{x}\|}\left(\frac{ \bm{x} \cdot \bm{y}}{\|\bm{x}\|^2}+\frac{\|\bm{y}\|^2}{2\|\bm{x}\|^2}-\frac{(\bm{x} \cdot \bm{y})^2}{\|\bm{x}\|^4}+\frac{C\|\bm{y}\|^2}{\|\bm{x}\|^{2+\varepsilon}}\right)\right. \\ &\quad \left.-\frac{1}{8 \log^2 \|\bm{x}\|} \frac{(\bm{x} \cdot \bm{y})^2}{\|\bm{x}\|^4}+\frac{C\|\bm{y}\|^2}{\|\bm{x}\|^{2+\varepsilon}\log^2 \|\bm{x}\|} \right),
\end{aligned}
$$
which proves (\ref{fxdeltaEvare}). 
\end{proof}

\begin{proof}[Proof of Proposition \ref{2DSRRWprop}(i)]
 Let $\nu,\varepsilon,\kappa,\tilde{\kappa}, (x_n)_{n> 1}$ be as in (\ref{defnueps}), (\ref{kappavarepass}), (\ref{defkappatilde}) and (\ref{defxndn}). We define inductively a sequence of stopping times: Set $n_1 = 1$; for $k>1$, let
$$
n_k:=\inf\left\{n>n_{k-1}: \|\bm{S}_{n} \|\leq n^{\tilde{\kappa}},\  \max\{|\Delta_n(\|\bm{x}\|^{2})|,|\Delta_n(\|\bm{x}\|^{2+\frac{\delta}{2}})|,\|\Delta_n(\bm{x}\bm{x}^T)\|\}  < \frac{1}{n^{\nu}}\right\},
$$
with the convention that $\inf \emptyset=\infty$. By Corollary \ref{limx2Delta} and Proposition \ref{Snsio}, almost surely, $n_k<\infty$ for all $k\geq 1$. Fix a constant $t\in (\tilde{\kappa},1/2)$ and let $r$ be a sufficiently large positive constant which will be specified later. For each $k\geq 1$, we define stopping times 
$$\tau_k:=\inf\{n\geq n_k: S_j \in B(0,r) \},\quad \theta_k:=\inf\{n \geq n_k: \|\bm{S}_n\| \geq n^t \},$$ 
and
$$T_k:=\inf\left\{n \geq n_k:  \max\{|\Delta_n(\|\bm{x}\|^{2})|,|\Delta_n(\|\bm{x}\|^{2+\frac{\delta}{2}})|,\|\Delta_n(\bm{x}\bm{x}^T)\|\}\geq \frac{1}{n^{\frac{\nu}{2}}}\right\}.
$$

 Let $f$ be as in (\ref{f2Ddef}). By (\ref{fxdeltaEvare}) in Lemma \ref{fxdelta2dest} with $\bm{x}=\bm{S}_n,\bm{y}=\bm{X}_{n+1}$ and Equations (\ref{approfirmo}), (\ref{approsecmo}) and (\ref{approSnXnmo}) in Lemma \ref{condEYq}, there exist positive constants $C_1$ and $r_1$ such that if $k\leq n<T_k$ and $\|\bm{S}_n\|\geq r_1$, then
   \begin{equation}
    \label{condfSE}
  \begin{aligned}
    &\quad\ \mathbb{E}\left((f\left(\bm{S}_{n+1}\right)-f\left(\bm{S}_n\right))\mathds{1}_{\{\bm{X}_{n+1} \in E_{\varepsilon}(\bm{S}_n)\}} \mid \FF_n\right)-\frac{\alpha}{2n\sqrt{\log\|\bm{S}_n\|}} \\ &\leq  \frac{1}{\|\bm{S}_n\|^2 \log^{3/2}\|\bm{S}_n\|}\left(-\frac{1}{8} +C_1\log \|\bm{S}_n\|\left(\frac{1}{\|\bm{S}_n\|^{\frac{\delta}{2}-\varepsilon(1+\frac{\delta}{2})}}+\frac{1}{\|\bm{S}_n\|^{\varepsilon}}+\frac{1}{n^{\frac{\nu}{2}}}\right)\right).
  \end{aligned}
\end{equation}
On the other hand, (\ref{fxdeltagener}) in Lemma \ref{fxdelta2dest} and (\ref{Xqest}) in Lemma \ref{condEYq} (with $q=0, 1$) imply that there exists a positive constant $C_2$ such that if $k\leq n<T_k$, then
\begin{equation}
  \label{condfSEc}
  \begin{aligned}
     &\quad \ \mathbb{E}\left((f\left(\bm{S}_{n+1}\right)-f\left(\bm{S}_n\right))\mathds{1}_{\{\bm{X}_{n+1} \notin E_{\varepsilon}(\bm{S}_n)\}} \mid \FF_n\right) \\
      &\leq  \frac{C_2}{\|\bm{S}_n\|^{(1-\varepsilon)(2+\frac{\delta}{2})}} =\frac{C_2}{\|\bm{S}_n\|^2 \log^{3/2}\|\bm{S}_n\|} \frac{\log^{3/2}\|\bm{S}_n\|}{\|\bm{S}_n\|^{\frac{\delta}{2}-\varepsilon(2+\frac{\delta}{2})}}.
  \end{aligned}
\end{equation}
Note that $\frac{\delta}{2}-\varepsilon(2+\frac{\delta}{2})>0$ by the choice of $\varepsilon$.
 By (\ref{condfSE}) and (\ref{condfSEc}), if $n_k\leq n<T_k$ and $\|\bm{S}_n\|\geq r_2>r_1$ for some large $r_2$ (note that $1/10$ below can be replaced by any number strictly less than $1/8$ if we choose a sufficiently large $r_2$), then
\begin{equation}
 \label{fsnrecsqrtlog}
 \mathbb{E}\left(f\left(\bm{S}_{n+1}\right)-f\left(\bm{S}_n\right) \mid \FF_n\right) \leq \frac{-1}{\|\bm{S}_n\|^2 \log ^{3/2}\|\bm{S}_n\|}\left(\frac{1}{10}-\frac{C_1\log \|\bm{S}_n\|}{n^{\frac{\nu}{2}}}-\frac{\alpha\|\bm{S}_n\|^2\log \|\bm{S}_n\|}{2n}\right).
\end{equation}
Observe that if $n_k\leq n< \theta_k$ and $k$ is large enough, say $k\geq K_1$, then
 $$\frac{C_1\log \|\bm{S}_n\|}{n^{\frac{\nu}{2}}} \leq \frac{C_1 t\log n}{n^{\frac{\nu}{2}}}<\frac{1}{20},\quad \frac{\|\bm{S}_{n}\|^2\log \|\bm{S}_{ n}\|}{ n} < \frac{t\log  n}{n^{1-2t}} < \frac{1}{20}.$$
Now fix $r=r_2$ and $k\geq K_1$. By (\ref{fsnrecsqrtlog}), $(f(\bm{S}_{(n_k+j)\wedge \tau_k \wedge T_k \wedge \theta_k}))_{j\in \NN}$ is a non-negative supermartingale, and thus converges a.s.. By (\ref{lilSRRW}), $\theta_k<\infty$ a.s.. Consequently, 
 $$\lim_{n\to \infty}f(\bm{S}_{(n_k+j)\wedge \tau_k \wedge T_k  \wedge \theta_k}) = f(\bm{S}_{\theta_k}) \quad \text{a.s. on}\ \{\tau_k=\infty\} \cap \{T_k=\infty\}.$$ 
 Then by the optional stopping theorem for non-negative supermartingales, 
 $$
 \begin{aligned}
    \sqrt{\tilde{\kappa}}\sqrt{\log n_k} &\geq \EE (f(\bm{S}_{n_k})\mid \FF_{n_k}) \\
    &\geq \EE (f(\bm{S}_{\theta_k})\mathds{1}_{\{\tau_k=\infty, T_k=\infty\}}\mid \FF_{n_k}) \\
    &\geq \EE (\sqrt{\log \theta_k^{t}}\mathds{1}_{\{\tau_k=\infty, T_k=\infty\}}\mid \FF_{n_k}) \\
    &\geq \sqrt{t} \sqrt{\log n_k} \PP(\{\tau_k=\infty\} \cap \{T_k=\infty\}\mid \FF_{n_k}).
 \end{aligned}
 $$
 Therefore, for all $k\geq K_1$, one has,
 \begin{equation}
   \label{inequtkcapTk}
     \PP(\{\tau_k=\infty\} \cap \{T_k=\infty\}\mid \FF_{n_k}) \leq \sqrt{\frac{\tilde{\kappa}}{t}}<1.
 \end{equation}
By Corollary \ref{limx2Delta},
 $$\PP(\bigcup_{m\geq k}\{T_m=\infty\})=1.$$ 
 It remains to use the arguments below (\ref{1dopsample}) to conclude that $\tau_k<\infty$ a.s..
\end{proof}

\subsubsection{Proof for the critical regime}
\label{sec2dcrit12}

In this section, we finish the proof of Proposition \ref{2DSRRWprop} (ii). As discussed in Section \ref{convquasimasec},  it suffices to show (\ref{sumposibd}), that is, 
$$
  \EE\left(\sum_{n=2}^{\infty} \EE (x_{n+1}-x_n\mid \FF_n)^{-} \right)<\infty,$$
  where $(x_n)_{n> 1}$ is defined in (\ref{defxndn}). The reader may find it instructive to compare its proof with the proof of Lemma \ref{lemznsupmart}.

\begin{proof}[Proof of Proposition \ref{2DSRRWprop} (ii)]
Again, to estimate $\EE (x_{n+1}-x_n\mid \FF_n)^{-}$, we consider the cases $\|\bm{X}_{n+1}\| \leq D_n^{1-\varepsilon}$ and $\|\bm{X}_{n+1}\| > D_n^{1-\varepsilon}$ separately.

 By Taylor expansion, there exists an $\varepsilon_1>0$ such that if $|x|\leq \varepsilon_1$, then
\begin{equation}
 \label{tayineln1xlow}
 \log (1+x) \geq x -\frac{1}{2}x^2-|x|^3.
\end{equation} 
Recall $u_{n+1}$ defined in (\ref{undef}). If $\|\bm{X}_{n+1}\| \leq D_n^{1-\varepsilon}$, by (\ref{unDbd}), $|u_{n+1}|<\varepsilon_1$ for large $n$ (say $n\geq m$). Thus, by (\ref{delxnformula}) and (\ref{tayineln1xlow}), for all $n\geq m$ such that $\|\bm{X}_{n+1}\| \leq D_n^{1-\varepsilon}$, 
\begin{equation}
  \label{xnrecura12low}
  x_{n+1}-x_n\geq \frac{1}{\log (n+1)}\left(u_{n+1}-\frac{1}{2}u_{n+1}^2-|u_{n+1}|^3\right)-\frac{\log D_n^2}{n \log n \log (n+1)}.
\end{equation}
Observe that 
$$
  \log D_{n+1}^2 - \log D_n^2 \geq -\log ((\frac{n}{n+1})^{\kappa}+\frac{\|\bm{S}_n\|^2}{(n+1)^{\kappa}}) 
  \geq -\log (1+\frac{\|\bm{S}_n\|^2}{n^{\kappa}})\geq -\frac{\|\bm{S}_n\|^2}{n^{\kappa}},
$$
and thus, by (\ref{delxnformula}), no matter whether $\|\bm{X}_{n+1}\| \leq D_n^{1-\varepsilon}$ or not, one has
\begin{equation}
  \label{delxngenerlow}
  x_{n+1}-x_n \geq -\frac{\|\bm{S}_n\|^2}{n^{\kappa}\log (n+1)}-\frac{\log D_n^2}{n \log n \log (n+1)}.
\end{equation}
By (\ref{unDvarbd}), (\ref{xnrecura12low}) and (\ref{delxngenerlow}),
$$
\begin{aligned}
   &\quad\ x_{n+1}-x_n \\  
   &\geq\frac{\mathds{1}_{\{\|\bm{X}_{n+1}\| \leq D_n^{1-\varepsilon}\} }}{\log (n+1)}\left(u_{n+1}- \frac{1}{2}u_{n+1}^2-\frac{4u_{n+1}^2 }{D_n^{\varepsilon}}\right) \\
   &\quad -\frac{\log D_n^2}{n \log n \log (n+1)}-\frac{\|\bm{S}_n\|^2\mathds{1}_{\{\|\bm{X}_{n+1}\| > D_n^{1-\varepsilon}\} }}{n^{\kappa}\log (n+1)} \\  &\geq\frac{1}{\log (n+1)}\left(u_{n+1}- \frac{1}{2}u_{n+1}^2\mathds{1}_{\{\|\bm{X}_{n+1}\| \leq D_n^{1-\varepsilon}\} }-\frac{\log D_n^2}{n \log n }\right)\\
  &\quad-\frac{4u_{n+1}^2 }{D_n^{\varepsilon}\log (n+1)}\mathds{1}_{\{\|\bm{X}_{n+1}\| \leq D_n^{1-\varepsilon}\} }-\frac{\|\bm{S}_n\|^2\mathds{1}_{\{\|\bm{X}_{n+1}\| > D_n^{1-\varepsilon}\} }}{n^{\kappa}\log (n+1)}-\frac{ |u_{n+1}|\mathds{1}_{\{\|\bm{X}_{n+1}\| > D_n^{1-\varepsilon}\} }}{\log (n+1)}.
\end{aligned}
$$
As a result,
\begin{equation}
  \label{delxnestequ}
\begin{aligned}
    (x_{n+1}-x_n)^-&\leq \frac{1}{\log (n+1)}\left(u_{n+1}- \frac{1}{2}u_{n+1}^2\mathds{1}_{\{\|\bm{X}_{n+1}\| \leq D_n^{1-\varepsilon}\} }-\frac{\log D_n^2}{n \log n }\right)^- \\
  &\quad +\frac{4u_{n+1}^2 }{D_n^{\varepsilon}\log (n+1)}+\frac{\|\bm{S}_n\|^2\mathds{1}_{\{\|\bm{X}_{n+1}\| > D_n^{1-\varepsilon}\} }}{n^{\kappa}\log (n+1)}+\frac{ |u_{n+1}|\mathds{1}_{\{\|\bm{X}_{n+1}\| > D_n^{1-\varepsilon}\} }}{\log (n+1)}.
\end{aligned}
\end{equation}
By Lemma \ref{condlemXn1} and that $D_n^2\geq \|\bm{S}_n\|^2$ and $D_n^2\geq n^{\kappa}$, for all $n\geq m$,
\begin{equation}
    \label{conderr1negli}
    \begin{aligned}
    \frac{\EE(\|\bm{S}_n\|^2\mathds{1}_{\{\|\bm{X}_{n+1}\| > D_n^{1-\varepsilon} \}}\mid \FF_n)}{n^{\kappa }\log (n+1)}  &\leq \frac{1}{n^{\kappa} \log (n+1) }\EE\left(\frac{ \|\bm{X}_{n+1}\|^{2+\frac{\delta}{2}}}{D_n^{(2+\frac{\delta}{2})(1-\varepsilon)-2} }\mid \FF_n\right) \\
    &\leq  \frac{2\EE \|\bm{X}_1\|^{2+\frac{\delta}{2}}+\Delta_n(\|\bm{x}\|^{2+\frac{\delta}{2}})}{2 n^{\kappa (1+\frac{\delta}{4})(1-\varepsilon)}  \log (n+1)}. 
\end{aligned}
\end{equation}
where we used that $(2+\frac{\delta}{2})(1-\varepsilon)-2>0$ by the choice of $\varepsilon$, i.e., (\ref{defnueps}). Similarly, by (\ref{unDbd}),
\begin{equation}
    \label{conderr2negli}
\begin{aligned}
 &\quad \ \EE(|u_{n+1}|\mathds{1}_{\{\|\bm{X}_{n+1}\| > D_n^{1-\varepsilon}\} }\mid \FF_n) \\
 &\leq \EE\left(\frac{2 \|\bm{X}_{n+1}\|^{2+\frac{\delta}{2}}}{D_n^{1+({2+\frac{\delta}{2}}-1)(1-\varepsilon)}}+\frac{\|\bm{X}_{n+1}\|^{2+\frac{\delta}{2}}}{D_n^{2+(2+\frac{\delta}{2}-2)(1-\varepsilon)}}+\frac{\|\bm{X}_{n+1}\|^{2+\frac{\delta}{2}}}{D_n^{2+(2+\frac{\delta}{2})(1-\varepsilon)}}\mid \FF_n\right) 
  \\ &\leq \EE\left(\frac{4 \|\bm{X}_{n+1}\|^{2+\frac{\delta}{2}}}{D_n^{\varepsilon+(2+\frac{\delta}{2})(1-\varepsilon)}}\mid \FF_n\right)\leq  \frac{4\EE \|\bm{X}_1\|^{2+\frac{\delta}{2}}+2\Delta_n(\|\bm{x}\|^{2+\frac{\delta}{2}})}{n^{\kappa ((1+\frac{\delta}{4})(1-\varepsilon)+\frac{\varepsilon}{2})}}.
\end{aligned}
\end{equation}
Note that Corollary \ref{limx2Delta}, $\EE |\Delta_n(\|\bm{x}\|^{{2+\frac{\delta}{2}}})|$ is upper-bounded for all $n\geq m$. Then, by (\ref{conderr1negli}), (\ref{conderr2negli}) and the assumption (\ref{kappavarepass}), we have
\begin{equation}
  \label{remintwosummable}
  \sum_{n=m}^{\infty} \frac{\EE(\|\bm{S}_n\|^2\mathds{1}_{\{\|\bm{X}_{n+1}\| > D_n^{1-\varepsilon} \}})}{n^{\kappa }\log (n+1)} <\infty, \quad \sum_{n=m}^{\infty}\EE(|u_{n+1}|\mathds{1}_{\{\|\bm{X}_{n+1}\| > D_n^{1-\varepsilon}\} })<\infty.
\end{equation}
Since $(x_n)_{n>1}$ is $L^1$ bounded, by (\ref{un2Dvarsum}), (\ref{delxnestequ}) and (\ref{remintwosummable}), we see that (\ref{sumposibd}) can be deduced from the following:
\begin{equation}
  \label{lognunlogS}
  \EE\left(\sum_{n=m}^{\infty}\frac{1}{\log (n+1)} \EE \left(u_{n+1}-\frac{1}{2}u_{n+1}^2\mathds{1}_{\{\|\bm{X}_{n+1}\| \leq D_n^{1-\varepsilon}\} }-\frac{\log D_n^2}{n \log n}\mid \FF_n\right)^{-} \right) <\infty.
\end{equation}
Now observe that 
$$
u_{n+1}-\frac{1}{2}u_{n+1}^2\mathds{1}_{\{\|\bm{X}_{n+1}\| \leq D_n^{1-\varepsilon}\} }-\frac{\log D_n^2}{n \log n}=I_n^{(1)}-\frac{1}{2}I_n^{(2)}+I_n^{(3)},
$$
where $I_n^{(i)}$, $i=1,2,3$, are defined in (\ref{defInun}).
On the other hand, since $(x-y)^-\leq x^-+y^+$ and $I_n^{(3)}$ is $\FF_n$-measurable, one has
$$
\begin{aligned} 
 &\quad \ \EE \left(u_{n+1}-\frac{1}{2}u_{n+1}^2\mathds{1}_{\{\|\bm{X}_{n+1}\| \leq D_n^{1-\varepsilon}\}}-\frac{\log D_n^2}{n \log n}\mid \FF_n\right)^{-} \\
 &\leq \EE (I_n^{(1)} \mid \FF_n)^-+ \frac{1}{2}\EE (I_n^{(2)}\mid \FF_n)^++ (I_n^{(3)})^{-}.
\end{aligned}
$$ 
Therefore, one can deduce (\ref{lognunlogS}) from (\ref{condEun2remsum}) in Lemma \ref{lemcondun} and (\ref{sumvn-bd}) in Lemma \ref{lemIn3sum}, which completes the proof of (\ref{sumposibd}). 
\end{proof}

\subsection{Proof of the main results}
\label{secpromain}

In this section, we finish the proofs of the results stated in Section \ref{secmainre}. Recall that we have proved Corollary \ref{cortransang} in Section \ref{RnYnsec} and proved Proposition \ref{MZslln} in Section \ref{stocappsec}.

\begin{proof}[Proof of Theorem \ref{mainthm} and Proposition \ref{highdspeed}.] Suppose that $\mu$ satisfies Assumption $A(2)$, then $\EE \bm{X}_1\bm{X}_1^T$ is a $d\times d$ positive-definite matrix since $\bm{S}$ is assumed to be genuinely $d$-dimensional.  Remark \ref{imagesrrw} shows that  
\begin{equation}
    \label{deftildes}
    \tilde{\bm{S}}:=(\EE \bm{X}_1\bm{X}_1^T)^{-\frac{1}{2}}\bm{S}
\end{equation}
is an SRRW with the same parameter $\alpha$. Its step distribution $\tilde{\mu}$ satisfies Assumption $A(2)$. Moreover:
\begin{itemize}
    \item $\tilde{\mu}$ satisfies $A(2+\delta)$ if $\mu$ satisfies Assumption $A(2+\delta)$.
    \item The covariance matrix of $\tilde{\mu}$ is the $d\times d$ identity matrix.
    \item For any $n\geq 1$, 
$$
\frac{\|\tilde{\bm{S}}_n\|}{\|(\EE \bm{X}_1\bm{X}_1^T)^{-\frac{1}{2}}\|}\leq \|\bm{S}_n\| \leq \|(\EE \bm{X}_1\bm{X}_1^T)^{\frac{1}{2}}\|\|\widehat{\bm{S}}_n\|.
$$
In particular, $\bm{S}$ is recurrent, resp. transient, if and only if $\tilde{\bm{S}}$ is recurrent, resp. transient. Moreover,
$$
\lim_{n\to \infty} \frac{\log \|\bm{S}_n\|}{\log n} = \lim_{n\to \infty} \frac{\log \|\bm{\tilde{S}}_n\|}{\log n} ,
$$
if either of the left-hand and right-hand limits exists.
\end{itemize}
Then, Theorem \ref{mainthm} and Proposition \ref{highdspeed} follow from (\ref{Snlimseries}) and Propositions \ref{RnYnxP0}, \ref{1DSRRWprop}, \ref{highDSRRWprop}, \ref{2DSRRWprop}.
\end{proof}

\begin{proof}[Proof of Corollary \ref{srrwlattice}.]
    By Theorem \ref{mainthm}, there exists $r>0$ such that $\mathbb{P}\left(\liminf_{n\to \infty}\|\bm{S}_n\|\leq r\right)=1$. Note that there are only finitely many sites in $L\cap B(0,r)$. Indeed, let $L$ be as in (\ref{deflattice}). Define a $d \times k$ matrix $A$ by $A:=(\bm{v}_1,\bm{v}_2,\dots,\bm{v}_k)$. If $\ell=\sum_{i=1}^k m_i \bm{v}_i \in L$ for some $m_1, \ldots, m_k \in \mathbb{Z}$, i.e. $\ell=A (m_1,m_2,\dots,m_k)^T$, then similarly as in (\ref{genuinekspace}), $\|\ell\|<r$ only if 
    $$
   \|(m_1,m_2,\dots,m_k)\| < \|(A^{T}A)^{-1}\|\|A^{-1}\|r,
    $$
    which holds for only finitely many vectors $(m_1,m_2,\dots,m_k)\in \mathbb{Z}^k$.

Recall that $\bm{S}$ is a Markov chain if $\alpha=0$. For $n\in \NN$, we denote by $P_n(\cdot,\cdot)$ the $n$-step transition probability for this Markov chain. Now observe that $\mathcal{R}$ does not depend on $\alpha$. Under our assumptions, by the Chung-Fuchs theorem and \cite[Theorem 5.4.1]{MR3930614}, $\mathcal{R}$ is a subgroup of $L$ under addition. In particular, for any $\bm{x},\bm{y} \in \mathcal{R}$, there exists an $m\in \NN$ such that $P_m(\bm{x},\bm{y})=P_m(0,\bm{y}-\bm{x})>0$, and thus 
$$
\PP(\bm{S}_{n+m}=\bm{y} \mid \FF_n) \mathds{1}_{\{\bm{S}_n=\bm{x}\}} \geq (1-\alpha)^mP_m(\bm{x},\bm{y}) \mathds{1}_{\{\bm{S}_n=\bm{x}\}}.
$$
The conditional Borel-Cantelli lemma (see e.g. \cite{MR0837655}) then implies that for any $\bm{x},\bm{y}\in \mathcal{R}$, $\bm{y}$ is visited infinitely often a.s. on the event that $\bm{x}$ is visited infinitely often. Since a.s. at least one site in $L\cap B(0,r)$ is visited infinitely often, a.s. all sites in $\mathcal{R}$ are visited by $\bm{S}$ infinitely often.
\end{proof}

\section{Acknowledgements}

The author would like to thank Jean Bertoin for introducing him to this question. The author also thanks Pierre Tarrès for some very helpful comments on this paper. The author is grateful to an anonymous referee for carefully reading this paper and providing many valuable comments on improving the manuscript.

\bibliographystyle{plain}
\bibliography{math_ref}

\end{document}